\documentclass[11pt]{amsart}
\usepackage[margin=1.1in]{geometry}
\numberwithin{equation}{section}
\usepackage{amsmath,amssymb,amsthm}
\usepackage{xcolor}
\usepackage{tikz}
\usetikzlibrary{decorations.pathreplacing,shapes.geometric}
\usepackage{algorithm}
\usepackage[noend]{algpseudocode}
\makeatletter
\providecommand\theHALG@line{}
\renewcommand\theHALG@line{\thealgorithm.\arabic{ALG@line}}
\makeatother
\algrenewcommand\algorithmicrequire{\textbf{Input:}}
\algnewcommand{\Keep}{\textsc{keep}}
\algnewcommand{\Discard}{\textsc{discard}}
\usepackage[colorlinks=true,linkcolor=blue!50!black,citecolor=blue!50!black,urlcolor=blue!50!black]{hyperref}

\newtheorem{theorem}{Theorem}[section]
\newtheorem{lemma}[theorem]{Lemma}
\newtheorem{proposition}[theorem]{Proposition}
\newtheorem{corollary}[theorem]{Corollary}

\DeclareMathOperator{\disc}{disc}
\DeclareMathOperator{\vol}{vol}
\DeclareMathOperator{\Var}{Var}
\DeclareMathOperator{\tr}{tr}
\DeclareMathOperator{\rank}{rank}
\DeclareMathOperator{\sgn}{sgn}
\DeclareMathOperator{\Bin}{Bin}
\DeclareMathOperator{\Poi}{Poi}
\DeclareMathOperator{\Gam}{Gamma}

\newcommand{\R}{\mathbb{R}}
\newcommand{\C}{\mathbb{C}}
\newcommand{\Z}{\mathbb{Z}}
\newcommand{\E}{\mathbb{E}}
\newcommand{\Prob}{\mathbb{P}}
\newcommand{\one}{\mathbf{1}}
\newcommand{\cA}{\mathcal{A}}
\newcommand{\cB}{\mathcal{B}}
\newcommand{\cD}{\mathcal{D}}
\newcommand{\cF}{\mathcal{F}}
\newcommand{\cG}{\mathcal{G}}
\newcommand{\cL}{\mathcal{L}}
\newcommand{\cP}{\mathcal{P}}
\newcommand{\cS}{\mathcal{S}}
\newcommand{\cV}{\mathcal{V}}
\newcommand{\eps}{\varepsilon}
\newcommand{\inner}[2]{\langle #1,#2\rangle}
\newcommand{\Dpre}{D^{\mathrm{pre}}}
\newcommand{\ii}{\mathrm{i}}

\title[Tusn\'ady's problem and online thinning]{Tusn\'ady's problem and online two-sample thinning}
\author{Yuhao Lei}
\address{Department of Mathematical Sciences, Carnegie Mellon University, Pittsburgh, PA 15213, USA}
\email{yuhaolei@andrew.cmu.edu}

\hypersetup{pdftitle={Tusnady's problem and online two-sample thinning},pdfauthor={Yuhao Lei}}

\begin{document}

\begin{abstract}We determine the sharp order of discrepancy in Tusnády's problem in every fixed dimension \(d\ge 2\). The largest combinatorial discrepancy of \(N\) points in \(\mathbb R^d\) with respect to axis-parallel boxes is \(\Theta_d((\log N)^{d-1/2})\), extending Wei's planar result and matching Nikolov's upper bound. More strongly, with high probability, this lower bound holds simultaneously for every subset obtained by deleting at most a sufficiently small fixed fraction of a uniform \(N\)-point sample.

We apply this geometric estimate to online thinning of two independent samples of size \(n\) from an unknown common distribution, where observations must be kept or discarded upon arrival. For every fixed dimension \(k\ge 1\) and expected deletion budget \(\rho n\), with \(0<\rho<2\), we prove that the minimax expected prefix discrepancy over lower orthants is \(\Theta_{k,\rho}((\log n)^{k+1/2})\). In dimension one, the optimal terminal discrepancy is instead \(\Theta_\rho((\log\log n)^{3/2})\), giving a sharp separation between the prefix and terminal objectives. For \(k\ge 2\), we obtain terminal bounds differing by a factor of \(\log n\).

The lower bounds are based on Gaussian small-ball estimates for box-incidence matrices, while the upper bounds use stationary Metropolis thinning and, in one dimension, a terminal cancellation procedure.
\end{abstract}

\maketitle

\section{Introduction}

\subsection{Tusn\'ady's problem}

Let $\cB_d$ be the family of half-open axis-parallel boxes $\prod_{a=1}^{d}[s_a,t_a)$ in $\R^d$. For a family $\cF$ of subsets of $\R^d$ and a finite set $X\subset\R^d$ put
\begin{equation}\label{eq:discdef}
\disc(X,\cF)=\min_{\chi\colon X\to\{-1,1\}}\ \sup_{F\in\cF}\Big|\sum_{x\in X\cap F}\chi(x)\Big|,\qquad\Delta_d(N)=\sup_{|X|=N}\disc(X,\cB_d).
\end{equation}
The boundary convention is immaterial on a finite point set. Tusn\'ady's problem asks for the order of $\Delta_d(N)$ as $N\to\infty$ with $d$ fixed. Beck \cite{Bec81} first studied the problem. Transference from the theorems of Roth and Schmidt gives lower bounds for combinatorial discrepancy \cite{Rot54,Sch72,ABN18}, and Matou\v{s}ek and Nikolov later obtained $\Omega_d((\log N)^{d-1})$ using factorization norms \cite{MN15,MNT20}. Following successive upper-bound improvements \cite{Bec89,Boh90,Sri97,Lar14,BG17}, Nikolov proved
\begin{equation}\label{eq:nikolov}
\Delta_d(N)\le C_d(\log N)^{d-1/2}\qquad(d\ge2,\ N\ge2).
\end{equation}
His argument uses Banaszczyk's theorem on signed series of vectors to handle one coordinate at a cost of $\sqrt{\log N}$ \cite{Nik17,Ban12}. Bansal conjectured that the lower bound $(\log N)^{d-1}$ was the correct order \cite[Conjecture 1.4.1]{Ban22}.

Wei proved the matching planar lower bound \cite{Wei26}. For every $A>0$, with probability at least $1-e^{-AN}$, every coloring of $N$ independent uniform points in the unit square has an anchored rectangle with imbalance at least $c_A(\log_2N)^{3/2}$. The proof reveals one coordinate digit by digit and accumulates the oscillation of signed counts along the other, using the opposing coordinate scales in Talagrand's Brownian-sheet argument \cite{Tal94}. Wei left the higher-dimensional order open \cite[Open Problem 2]{Wei26}.

For $X\subset(0,1)^d$ let
\[
\cA_d=\{[0,t):t\in[0,1]^d\},\qquad[0,t)=\prod_{a=1}^d[0,t_a),
\]
be the family of anchored boxes. Inclusion and exclusion express the indicator of a box in $\cB_d$ as a signed sum of at most $2^d$ indicators of anchored boxes, hence
\begin{equation}\label{eq:anchored}
\disc(X,\cA_d)\le\disc(X,\cB_d)\le2^d\disc(X,\cA_d).
\end{equation}

\begin{theorem}\label{thm:main}
For every fixed $d\ge2$ and all $N\ge2$,
\[
c_d(\log N)^{d-1/2}\le\Delta_d(N)\le C_d(\log N)^{d-1/2}.
\]
There is also a constant $\eta_d>0$ with the following property. If $X_N$ consists of $N$ independent uniform points in $(0,1)^d$, then for all sufficiently large $N$, with probability at least $1-\exp(-b_dN)$,
\[
\disc(Y,\cA_d)\ge c_d(\log N)^{d-1/2}\qquad\text{for every }Y\subseteq X_N\text{ with }|X_N\setminus Y|\le\eta_dN.
\]
\end{theorem}

The first assertion supplies the lower bound for $d\ge3$. In dimension one, alternating colors in increasing order gives $\Delta_1(N)=1$. The second assertion holds simultaneously over the indicated subsets and their colorings.

For independent uniform arrivals, every online coloring has discrepancy at least $c_d(\log T)^{d-1/2}$ at each sufficiently large fixed horizon $T$, with probability at least $1-e^{-b_dT}$. The previous lower bound in this model was $\Omega_d((\log T)^{d/2})$ \cite{BJSS20}. Wei noted the corresponding planar consequence \cite{Wei26}.

\subsection{Online two-sample thinning}\label{sec:model}

Let $X_1,\dots,X_n$ and $Y_1,\dots,Y_n$ be independent observations from a common Borel probability measure $\mu$ on $\R^k$. The queries are the lower orthants $(-\infty,b]=\prod_{a=1}^k(-\infty,b_a]$, $b\in\R^k$, and inequalities between vectors are understood coordinatewise. Give the two samples the signs $+1$ and $-1$, respectively. At each step an online algorithm may stop or select a sample with an unread observation. It selects the sample using an independent random seed and the past, then observes the next point and immediately keeps or discards it. Decisions are irrevocable. After an early stop, complete the reveal order by listing the remaining observations in a fixed order and discarding them. Write $R$ for the total number deleted, $x_i$ for the $i$-th observation in this order and $\eps_i$ for its source sign. Define
\begin{equation}\label{eq:St}
S_t(b)=\sum_{\substack{i\le t\\ x_i\text{ kept}}}\eps_i\one\{x_i\le b\},\qquad D=\sup_{b\in\R^k}|S_{2n}(b)|,\qquad\Dpre=\max_{0\le t\le2n}\sup_{b\in\R^k}|S_t(b)|.
\end{equation}
Both discrepancies are unnormalized signed counts. The budget $\E R\le\rho n$ allows an expected fraction $\rho/2$ of the $2n$ observations to be discarded. Thus $\rho=1/8$ corresponds to one sixteenth of the input, while $\rho=2$ allows deletion of both samples. The algorithms use $n$ and $\rho$ and are common to all input laws. Without deletions, the expected discrepancy on the line has order $\sqrt n$ for every fixed nondegenerate law, by the Dvoretzky--Kiefer--Wolfowitz--Massart inequality and the central limit theorem \cite{Mas90}.

Dwivedi, Feldheim, Gurel-Gurevich and Ramdas introduced online thinning for a single uniform stream, where at least one of every two observations must be kept \cite{DFGR19}. Their model is related to the power of two choices \cite{ABKU99}. Other single-sample variants allow deletion \cite{BS25} or relocation \cite{SV26c}.

Smirnov and Vershynin introduced two-sample thinning \cite{SV26a}. For arbitrary laws on $\R^k$, they obtained discrepancy $O_{k,\rho}((\log n)^{2k})$. Their algorithm represents each point by its dyadic intervals and rejects steps that would leave a cube. Their subsequent stationary Metropolis construction gives the exponent $3k/2$ \cite[Section 1.7]{SV26b}. They asked for the optimal order at a fixed deletion fraction. In this model, the source determines the sign and the algorithm chooses whether to keep the point. Related online discrepancy models allow the algorithm to choose the sign, under stochastic arrivals \cite{BJSS20,BJMSS21} or oblivious adversaries \cite{ALS21,KRR24}.

Let $\cV^{\mathrm{pre}}_{k,\rho}(n)$ and $\cV^{\mathrm{end}}_{k,\rho}(n)$ be the infima of $\sup_\mu\E_\mu\Dpre$ and $\sup_\mu\E_\mu D$, respectively, over online rules satisfying $\sup_\mu\E_\mu R\le\rho n$. Each rule may use $k,n,\rho$, an independent random seed and the observed history. The suprema range over all Borel probability measures on $\R^k$.

\begin{theorem}\label{thm:multi}
For every fixed integer $k\ge1$ and every fixed $0<\rho<2$,
\[
\cV^{\mathrm{pre}}_{k,\rho}(n)\asymp_{k,\rho}(\log n)^{k+1/2}\qquad(n\to\infty).
\]
For every fixed $k\ge2$ in the same budget range,
\[
c_{k,\rho}(\log n)^{k-1/2}\le\cV^{\mathrm{end}}_{k,\rho}(n)\le C_{k,\rho}(\log n)^{k+1/2}
\]
for all sufficiently large $n$. Both lower bounds already hold for the uniform law on $(0,1)^k$. More precisely, for $k\ge1$ there are constants $c_k,b_k>0$ and $n_k$ such that, whenever $0<\rho\le1/8$ and $n\ge n_k$, every online algorithm with $\E R\le\rho n$ satisfies
\[
\Prob\{\Dpre\ge c_k(\log n)^{k+1/2}\}\ge1-\rho-e^{-b_kn}.
\]
For $k\ge2$, the same probability bound holds for $\{D\ge c_k(\log n)^{k-1/2}\}$, including for procedures which see both samples before deleting.
\end{theorem}

For $k\ge2$, closing the factor $\log n$ between the terminal bounds remains open. On the line, later arrivals can cancel the net count in each of several intervals, giving the sharper terminal estimate below.

\begin{theorem}\label{thm:terminal}
Fix $0<\rho<1/8$.

(i) There is an online thinning algorithm, which uses $n$ and $\rho$ but not $\mu$, such that for every Borel probability measure $\mu$ on $\R$ and every $n\ge2$,
\begin{equation}\label{eq:terminal-upper}
\E R\le\rho n,\qquad\E D\le\frac C\rho\Big(\log\log(2n)+\log\frac2\rho\Big)^{3/2},
\end{equation}
where $C>0$ is an absolute constant. Whenever $\rho^2n\ge2\cdot10^4\log(2n)$, the same algorithm satisfies
\begin{equation}\label{eq:terminal-whp}
\Prob\Big(D>\frac C\rho\Big(\log\log(2n)+\log\frac2\rho\Big)^{3/2}\Big)\le4(2n)^{-9}.
\end{equation}

(ii) There are absolute constants $c>0$ and $n_0$ such that, for the uniform distribution on $[0,1]$ and every $n\ge n_0$, every online thinning algorithm with $\E R\le n/8$ satisfies
\begin{equation}\label{eq:terminal-lower}
\Prob\big(D\ge c(\log\log n)^{3/2}\big)\ge\frac14.
\end{equation}
\end{theorem}

\begin{corollary}\label{cor:terminal-allrho}
For every fixed $0<\rho<2$,
\[
\cV^{\mathrm{end}}_{1,\rho}(n)
\asymp_\rho (\log\log n)^{3/2}\qquad(n\to\infty).
\]
More precisely, there are constants $c_\rho>0$ and $n_\rho$ such that,
for the uniform law on $(0,1)$ and every $n\ge n_\rho$, every online
thinning algorithm with $\E R\le\rho n$ satisfies
\[
\Prob\bigl(D\ge c_\rho(\log\log n)^{3/2}\bigr)
\ge\frac{2-\rho}{16}.
\]
\end{corollary}

In dimension one, the sharp prefix order follows from Theorem~\ref{thm:multi}. The upper bound has the order obtained in \cite[Section 1.7]{SV26b}, while the construction below also gives a pathwise bound; see Section~\ref{sec:prefix-upper}. The matching lower and upper bounds show that requiring control at every prefix costs a factor of order $(\log n/\log\log n)^{3/2}$ compared with the terminal objective.

\paragraph{Offline comparison.}
For $0<\rho<1/8$, the optimal one-dimensional terminal discrepancy with offline deletions is of order $1/\rho$. More precisely, for the uniform law on $[0,1]$, every procedure with $\E R\le\rho n$ satisfies $\E D\ge1/(64\rho)$ whenever $n\ge\rho^{-2}/16$. Conversely, for every Borel law on $\R$ and $n\ge1$, a distribution-free offline procedure satisfies $\E R\le\rho n$ and $D\le\lceil1/\rho\rceil$ almost surely. The argument is close to the reflected-walk bounds of \cite[Proposition 1.4 and Section 3]{SV26b}, but the sorted two-sample signs form a balanced sequence rather than independent signs.

For the lower bound, partition $[0,1]$ into $m=\max\{2,\lceil16\rho^2n\rceil\}$ equal intervals, where $16\rho^2n\le m\le32\rho^2n\le n/2$, and let $\Delta_j$ be the original signed count in interval $j$. Since each retained interval count has absolute value at most $2D$,
\[
\sum_{j=1}^m|\Delta_j|\le R+2mD.
\]
Writing $V=\E\Delta_j^2=(2n/m)(1-1/m)\ge n/m\ge1$, independence within each sample gives $\E\Delta_j^4\le3V^2+V\le4V^2$. H\"older's inequality yields $\E|\Delta_j|\ge\sqrt V/2$, and therefore
\[
\frac12\sqrt{nm}\le\rho n+2m\E D.
\]
As $16\rho^2n\le m\le32\rho^2n$, this implies $\E D\ge1/(64\rho)$.

For the upper bound, sort all observations using independent tie-breaking tags, and write their source signs as $\eps_1,\ldots,\eps_{2n}$. They form a uniformly random balanced sign sequence. Use the reflected-walk rule of \cite[Section 3]{SV26b}: starting from $S_0=0$, keep $\eps_{t+1}$ precisely when $|S_t+\eps_{t+1}|\le B$, where $B=\lceil1/\rho\rceil$. Then $D\le B$ and
\[
S_{t+1}^2-S_t^2
=2S_t\eps_{t+1}+1-(2B+1)\one\{\text{step }t+1\text{ is discarded}\}.
\]
Conditional on $\eps_{t+1}=+1$, the remaining $2n-1$ positions contain one fewer positive sign than conditional on $\eps_{t+1}=-1$. Coupling these two conditional sequences by inclusion and using the monotonicity of the reflected walk gives $\E[S_t\eps_{t+1}]\le0$. Summing the preceding identity, we obtain
\[
(2B+1)\E R
=2n+2\sum_{t=0}^{2n-1}\E[S_t\eps_{t+1}]-\E S_{2n}^2
\le2n,
\]
and hence $\E R\le2n/(2B+1)<\rho n$.

Under an expected deletion budget, an algorithm may discard all observations with probability $\rho/2$, so a positive probability of zero discrepancy is possible. The thinning lower bounds use the fixed source signs. The geometric deletion statement of Theorem~\ref{thm:main} applies to every coloring of every sufficiently large subset.

\subsection{Ideas of the proofs}

The geometric lower bound rests on a Gaussian small-ball estimate that is uniform over subspaces of dimension proportional to $N$. The dyadic decomposition of an anchored box separates coarse incidence vectors, whose span has small dimension compared with $N$, from a fine part \cite{Mat99}. After the coarse span is removed, a spectral estimate bounds the integrated variance below by a constant times $(\log N)^{2d-2}$. The further factor $\sqrt{\log N}$ comes from overlap estimates through moment order $\log N$. These estimates make a normalized Gaussian moment sufficiently concentrated for the inequality of Paouris and Valettas to give the desired small-ball probability \cite{PV18}.

The overlap calculation is motivated by Neumayer and Ullrich's use of common dyadic prefixes in their Brownian-sheet small-deviation argument \cite[Section 4]{NU26}. Their estimate matches the probability lower bound of Dunker, Linde, K\"uhn and Lifshits \cite{DLKL99}. In our setting the overlaps are between finite-sample incidence vectors. The Gaussian estimate first holds on a sample event of fixed positive probability. A bounded-differences argument shows that, with exponentially high probability, a sample differs in only a small proportion of its points from one in this event. Uniformity over large subspaces allows the Gaussian estimate to survive these changes. Binary rounding then yields the coloring lower bound \cite{LSV86}.

For thinning, the same Gaussian estimate and Pisier's regular ellipsoid theorem bound the number of ternary vectors with small discrepancy \cite{Pis89}. Folding the first coordinate transforms the source signs into independent fair signs, independent of the folded locations, even under adaptive source selection. Each ternary vector with a large support has exponentially small probability of agreeing with these signs, so a union bound excludes all large retained sets. Recording arrival order in one additional coordinate gives the prefix lower bound.

The upper bounds use Smirnov and Vershynin's signed-stream reduction and stationary Metropolis rule \cite{SV26a,SV26b}. A pilot sample sets the interval boundaries for the dyadic counts. Conditional on the locations, the Metropolis coordinates are independent at each fixed time, so their sums over a query concentrate. On the line, later arrivals cancel the net count in each bin. Once all bin totals vanish, only the bin containing the threshold contributes, allowing dyadic depth of order $\log\log n$. For the matching lower bound, Wei's planar estimate forces many bins to cross a discrepancy threshold. A long interarrival gap then preserves one crossing until the terminal time.

\subsection{Notation}

Logarithms are natural unless a base is indicated. We write $[N]=\{1,\ldots,N\}$ and $\vol$ for Lebesgue measure. The dimensions $d$ and $k$ are fixed. Constants denoted $c,C$ may change from line to line and depend only on the fixed parameters specified in the text. In the geometric arguments, they are uniform in $N$, the dyadic depth and the moment order. For symmetric matrices, $K\succeq L$ means that $K-L$ is positive semidefinite.

\section{A Gaussian criterion for box discrepancy}\label{sec:gauss}

For a sample $X=(x_1,\dots,x_N)$ of points in $[0,1]^d$ and $t\in[0,1]^d$ put
\[
a(t)=a_X(t)=\big(\one\{x_i\in[0,t)\}\big)_{i=1}^N\in\{0,1\}^N,
\]
and let $A_X$ be the matrix whose rows are the distinct vectors $a(t)$, $t\in[0,1]^d$. This matrix has at most $(N+1)^d$ rows. Then $\|A_Xg\|_\infty=\sup_t|\inner{a(t)}{g}|$ for every $g\in\R^N$. For a linear subspace $V\subseteq\R^N$ let $P_V$ be the orthogonal projection onto $V$. The following finite-sample estimate is the input to the rounding argument of Section~\ref{sec:disc} and to the counting argument of Section~\ref{sec:fixed}.

\begin{proposition}\label{prop:gauss}
Fix $d\ge2$, $0<\alpha<1$ and $B>0$. There are constants $a_{d,\alpha,B}>0$ and $b_{d,\alpha}>0$ such that for every sufficiently large $N$, with probability at least $1-\exp(-b_{d,\alpha}N)$ over the uniform sample $X$, every subspace $V\subseteq\R^N$ of dimension at least $\alpha N$ satisfies
\[
\Prob_g\big(\|A_XP_Vg\|_\infty\le a_{d,\alpha,B}(\log N)^{d-1/2}\big)\le e^{-BN},\qquad g\sim N(0,I_N).
\]
\end{proposition}

For each deterministic sample $X$, the probability is over $g$. The sample event is common to all $B>0$. We begin with a small-ball estimate in terms of the moments of the correlations of a projected Gaussian process.

\subsection{Gaussian inequalities}\label{sec:ineq}

Let $\gamma_n$ denote the standard Gaussian measure on $\R^n$. The log-concavity theorem of Pr\'ekopa \cite{Pre73} and Borell \cite{Bor74} states that, for convex sets $K_1,K_2\subseteq\R^n$ and $0<\lambda<1$,
\begin{equation}\label{eq:logconcave}
\gamma_n\big(\lambda K_1+(1-\lambda)K_2\big)\ge\gamma_n(K_1)^\lambda\gamma_n(K_2)^{1-\lambda}.
\end{equation}
The same inequality holds for the law of a centered Gaussian vector with singular covariance, which has a Gaussian density on the range $E$ of its covariance matrix, for convex subsets of $E$.

Let $Y$ and $Z$ be independent random vectors in $\R^n$, let $Y$ be centered Gaussian, possibly with singular covariance, and let $K\subset\R^n$ be a compact convex set with $K=-K$. Anderson's inequality \cite{And55} gives
\begin{equation}\label{eq:anderson}
\Prob\{Y+Z\in K\}\le\Prob\{Y\in K\}.
\end{equation}
The small ball inequality of Paouris and Valettas for norms
\cite[Proposition 3.2]{PV18} rests on their small deviation
inequality for convex functions, the B-theorem of \cite{CFM04}
and the small ball parameter of Klartag and Vershynin \cite{KV07}. There is an absolute constant $c_1>0$ such that for every norm $\|\cdot\|$ on $\R^n$, with $M$ a median of $\|g\|$ for $g\sim\gamma_n$,
\begin{equation}\label{eq:pv}
\Prob\{\|g\|\le\eta M\}\le\frac12\,\eta^{\,c_1M^2/\Var\|g\|}\qquad(0<\eta<1/2).
\end{equation}
The same estimate holds for every nonzero seminorm $f$. Indeed, $f$ is constant on cosets of its kernel and restricts to a norm on $(\ker f)^\perp$. The orthogonal projection of $g$ onto this subspace is a standard Gaussian there and has the same value of $f$, hence the same median and variance.

\subsection{The normalized moment criterion}

\begin{lemma}\label{lem:criterion}
Let $T$ be a measurable space, let $a\colon T\to\R^N$ be measurable with finitely many values, and let $P$ be an orthogonal projection on $\R^N$. Let $\mu$ be a probability measure on $T$ and $s>0$ a number with $\sigma(t)=\|Pa(t)\|_2\ge s$ for $\mu$-almost every $t$. Put $u(t)=Pa(t)/\sigma(t)$ where $\sigma(t)>0$ and $r(t,v)=\inner{u(t)}{u(v)}$. Let $p$ be an even integer with $2\le p\le\log_2N$, and suppose that
\begin{equation}\label{eq:varhyp}
\sum_{\substack{2\le j\le p\\ j\text{ even}}}\frac{p^j}{j!}\iint r(t,v)^j\,d\mu(t)\,d\mu(v)\le C\frac{p^2}N.
\end{equation}
Then there is $N_0=N_0(C)$ such that for every $B>0$ there is $c=c(C,B)>0$ with
\begin{equation}\label{eq:smallball}
\Prob_g\Big\{\sup_{t\in T}|\inner{a(t)}{g}|\le c\,s\sqrt p\Big\}\le e^{-BN}\qquad(N\ge N_0),\quad g\sim N(0,I_N).
\end{equation}
\end{lemma}

\begin{proof}
Let $m_p=\E|\zeta|^p$ for a standard normal variable $\zeta$, and set
\[
U(g)=\frac1{m_p}\int_T|\inner{u(t)}{g}|^p\,d\mu(t),
\qquad f(g)=U(g)^{1/p}.
\]
Since $\|u(t)\|_2=1$ almost everywhere, $\E_gU=1$ and $U$ has
finite Gaussian moments of all orders. The function $U$ is convex,
and Minkowski's inequality shows that $f$ is a nonzero seminorm.

For two standard normal variables with correlation $r$, Isserlis'
theorem \cite{Iss18}, grouping the pairings by their number $j$ of
cross pairs, gives
\[
\frac{\E[\zeta_1^p\zeta_2^p]}{m_p^2}
=1+\sum_{\substack{2\le j\le p\\j\text{ even}}}
\frac{\big(p(p-2)\cdots(p-j+2)\big)^2}{j!}\,r^j
\le 1+\sum_{\substack{2\le j\le p\\j\text{ even}}}
\frac{p^j}{j!}\,r^j.
\]
Integrating this inequality and using \eqref{eq:varhyp}, we obtain
\begin{equation}\label{eq:v0}
v_0:=\Var_gU\le C\frac{p^2}N.
\end{equation}

Convexity of $U$ and \eqref{eq:logconcave} imply that
$y\mapsto\Prob_g\{U\le y\}$ is log-concave on $(0,\infty)$.
Since $7/8$ is the midpoint of $1/2$ and $5/4$, Chebyshev's
inequality gives, once $v_0\le1/32$,
\begin{equation}\label{eq:Uhalf}
\Prob_g\{U\le1/2\}
\le\frac{\Prob_g\{U\le7/8\}^2}{\Prob_g\{U\le5/4\}}
\le\frac{(64v_0)^2}{1-16v_0}
\le2^{13}v_0^2.
\end{equation}
The derivative of $y^{1/p}$ is at most $2/p$ on $[1/2,\infty)$.
On $\{U\le1/2\}$ we have $|f-1|\le1$, so
\begin{equation}\label{eq:fvar}
\E_g(f-1)^2
\le\frac{4v_0}{p^2}+2^{13}v_0^2
\le\frac{4C}N+2^{13}C^2\frac{p^4}{N^2}
\le\frac{5C}N
\end{equation}
for $N\ge N_0(C)$, because $p\le\log_2N$.
Thus $\Var f\le5C/N$ and
$\Prob_g\{f<1/2\}\le4\E_g(f-1)^2<1/2$ for sufficiently large $N$.
Every median $M$ of $f(g)$ therefore satisfies $M\ge1/2$.

Applying \eqref{eq:pv} to $f$, we get
\[
\Prob_g\{f\le\eta M\}
\le\frac12\eta^{c_1N/(20C)}\qquad(0<\eta<1/2).
\]
Given $B>0$, choose $\eta=\eta(C,B)\in(0,1/2)$ such that
$(c_1/(20C))\log(1/\eta)\ge B$. Then
\begin{equation}\label{eq:feta}
\Prob_g\{f\le\eta/2\}
\le\Prob_g\{f\le\eta M\}\le e^{-BN}.
\end{equation}
Since $m_p=(p-1)!!\ge(p/2)!\ge(p/(2e))^{p/2}$,
\[
f(g)\le
\frac{\sup_t|\inner{Pa(t)}{g}|}{s\,m_p^{1/p}}.
\]
With $c=\eta/(2\sqrt{2e})$, the event
$\sup_t|\inner{Pa(t)}{g}|\le cs\sqrt p$ therefore implies
$f(g)\le\eta/2$.

Let $A$ have the distinct values of $a$ as its rows.
The vectors $APg$ and $A(I-P)g$ are independent centered Gaussians,
so Anderson's inequality \eqref{eq:anderson} yields
\[
\Prob_g\{\|Ag\|_\infty\le cs\sqrt p\}
\le\Prob_g\{\|APg\|_\infty\le cs\sqrt p\}
\le e^{-BN}.
\]
The thresholds imposed on $N$ depend only on $C$.
\end{proof}

\section{Projection and dyadic overlaps}\label{sec:layers}

\subsection{Variance on large subspaces}\label{sec:spectral}

Fix $d\ge2$ and put $q=d-1$. Let $x_1,\dots,x_N$ be independent uniform points in $(0,1)^d$. Define
\begin{equation}\label{eq:gram}
K_X=\int_{[0,1]^d}a(t)a(t)^\top\,dt,
\end{equation}
so that $w^\top K_Xw=\int\inner{w}{a(t)}^2\,dt$ for $w\in\R^N$. We need a sample event on which every subspace of dimension proportional to $N$ retains a positive amount of this variance, uniformly over the subspace. The tests are complex exponentials, whose sampled values have constant modulus.

\begin{lemma}\label{lem:sampled}
Fix $0<\delta<1$. For all sufficiently large $N$, with probability at least $3/4$ over the sample, there is an orthogonal projection $E$ on $\R^N$ of rank at least $(1-\delta)N$ such that
\begin{equation}\label{eq:operator}
K_X\succeq b_{d,\delta}\frac{(\log N)^{2q}}N\,E.
\end{equation}
\end{lemma}

\begin{proof}
Write $\mathbb{N}=\{1,2,\dots\}$. For an integer $l\ge0$ consider the frequency layer
\[
\cL_l=\Big\{m\in\mathbb{N}^d:\ \sum_{a=1}^d\lfloor\log_2m_a\rfloor=l\Big\}.
\]
The number of $m_a\in\mathbb{N}$ with $\lfloor\log_2m_a\rfloor=r_a$ is $2^{r_a}$, so $|\cL_l|=\sum_{r_1+\cdots+r_d=l}2^l=2^l\binom{l+q}q$, the sum being over $r\in\Z_{\ge0}^d$. Put $M=|\cL_l|$ and let $l$ be the least integer with $M\ge64N/\delta$. Since $|\cL_0|=1<64N/\delta$, we have $l\ge1$, and since $|\cL_l|/|\cL_{l-1}|=2(l+q)/l\le2d$, minimality gives
\begin{equation}\label{eq:Mbounds}
\frac{64N}\delta\le M<\frac{128dN}\delta.
\end{equation}
Every $m\in\cL_l$ satisfies $m_a<2^{\lfloor\log_2m_a\rfloor+1}$, hence $\prod_am_a<2^{l+d}$. For $N\ge N_{d,\delta}$ we also have $l>\frac12\log_2N$. Otherwise $M\le2^l(l+q)^q\le\sqrt N(\tfrac12\log_2N+q)^q<64N/\delta$, contradicting \eqref{eq:Mbounds}. Finally $\binom{l+q}q\ge l^q/q!$ and \eqref{eq:Mbounds} give
\begin{equation}\label{eq:2l}
2^l\le\frac{q!\,M}{l^q}<\frac{128d\,q!\,N}{\delta\,l^q}.
\end{equation}

Let $F$ be the $N\times M$ matrix with entries $F_{im}=\exp(2\pi\ii\inner{m}{x_i})$, $i\in[N]$, $m\in\cL_l$, where $\ii^2=-1$, and let $Q=\overline F F^\top/M$. Thus $Q$ is a Hermitian positive semidefinite $N\times N$ matrix with
\[
Q_{ij}=\frac1M\sum_{m\in\cL_l}e^{2\pi\ii\inner{m}{x_j-x_i}},\qquad v^*Qv=\frac1M\|F^\top v\|_2^2\quad(v\in\C^N).
\]
Its diagonal entries are $Q_{ii}=1$. For $i\ne j$, the independence of $x_i$ and $x_j$ and the identity $\E e^{2\pi\ii\inner{k}{x}}=\one\{k=0\}$ for $k\in\Z^d$ and a uniform point $x$ give
\[
\E|Q_{ij}|^2=\frac1{M^2}\sum_{m,m'\in\cL_l}\E e^{2\pi\ii\inner{m-m'}{x_j}}\,\E e^{-2\pi\ii\inner{m-m'}{x_i}}=\frac1{M^2}\sum_{m,m'\in\cL_l}\one\{m=m'\}=\frac1M.
\]
Hence $\E\|Q-I\|_{\mathrm{HS}}^2=N(N-1)/M\le\delta N/64$ by \eqref{eq:Mbounds}, and Markov's inequality gives $\|Q-I\|_{\mathrm{HS}}^2\le\delta N/16$ with probability at least $3/4$. On this Borel event, which we call the spectral event, the eigenvalues $\lambda_1\ge\cdots\ge\lambda_N\ge0$ of $Q$ satisfy $\sum_i(\lambda_i-1)^2\le\delta N/16$, and each eigenvalue below $1/2$ contributes more than $1/4$ to this sum, so at most $\delta N/4$ of them are below $1/2$. Let $\cS_1\subseteq\C^N$ be the span of the eigenvectors of $Q$ with eigenvalues at least $1/2$, so that $\dim_\C\cS_1\ge(1-\delta/4)N$ and
\begin{equation}\label{eq:S1}
\|F^\top v\|_2^2=Mv^*Qv\ge\frac M2\|v\|_2^2\qquad(v\in\cS_1).
\end{equation}

For $v\in\C^N$ consider the bounded function $f_v(t)=v^*a(t)=\sum_i\overline{v_i}\,a_i(t)$ on $[0,1]^d$, whose squared norm in $L^2([0,1]^d)$ is $\int|f_v|^2\,dt=v^*K_Xv$. The function $t\mapsto a_i(t)$ is the indicator of $\prod_a(x_{ia},1]$, and for an integer $m_a\ge1$ and $x\in[0,1]$,
\[
\int_x^1e^{-2\pi\ii m_as}\,ds=\frac{e^{-2\pi\ii m_a}-e^{-2\pi\ii m_ax}}{-2\pi\ii m_a}=\frac{1-e^{-2\pi\ii m_ax}}{-2\pi\ii m_a}.
\]
By Fubini's theorem, the Fourier coefficient of $f_v$ at a frequency $m\in\cL_l$ is therefore
\[
\int f_v(t)e^{-2\pi\ii\inner{m}{t}}\,dt=\sum_i\overline{v_i}\prod_{a=1}^d\frac{1-e^{-2\pi\ii m_ax_{ia}}}{-2\pi\ii m_a},\qquad\Big|\int f_v(t)e^{-2\pi\ii\inner{m}{t}}\,dt\Big|=\frac{|(G^\top v)_m|}{\prod_a2\pi m_a},
\]
Here $G$ is the $N\times M$ matrix with entries $G_{im}=\prod_{a=1}^d(1-e^{2\pi\ii m_ax_{ia}})$, and the second identity follows by conjugating inside the modulus. The functions $t\mapsto e^{2\pi\ii\inner{m}{t}}$, $m\in\cL_l$, are orthonormal in $L^2([0,1]^d)$, so Bessel's inequality and $\prod_am_a<2^{l+d}$ give
\begin{equation}\label{eq:besselG}
v^*K_Xv=\int|f_v|^2\,dt\ge\sum_{m\in\cL_l}\frac{|(G^\top v)_m|^2}{\prod_a(2\pi m_a)^2}\ge(2\pi)^{-2d}4^{-l-d}\,\|G^\top v\|_2^2\qquad(v\in\C^N).
\end{equation}

Expanding the product defining $G_{im}$,
\[
G_{im}=\sum_{I\subseteq[d]}(-1)^{|I|}\exp\Big(2\pi\ii\sum_{a\in I}m_ax_{ia}\Big)=(-1)^dF_{im}+H_{im},
\]
where
\[
H_{im}=\sum_{I\subsetneq[d]}(-1)^{|I|}\exp\Big(2\pi\ii\sum_{a\in I}m_ax_{ia}\Big).
\]
For a fixed proper subset $I\subsetneq[d]$, the corresponding
matrix depends on $m$ only through $(m_a)_{a\in I}$.
If $|I|=j\ge1$, it has at most
\[
\sum_{r=0}^l2^r\binom{r+j-1}{j-1}
\le2^{l+1}(l+1)^{j-1}
\]
distinct columns. The term $I=\emptyset$ has rank one. Hence
$\rank H\le C_d2^l(l+1)^{d-2}$,
and by \eqref{eq:2l} and $l>\frac12\log_2N$,
\[
\rank H\le C_d\frac{128d\,q!\,N}{\delta}\,\frac{(l+1)^{d-2}}{l^{d-1}}\le\frac{C_d'N}{\delta\,l}\le\frac{3\delta N}4\qquad(N\ge N_{d,\delta}).
\]
Let $\cS=\cS_1\cap\ker H^\top$.
Then $\dim_\C\cS\ge\dim_\C\cS_1-\rank H\ge(1-\delta)N$.
For $v\in\cS$ we have $G^\top v=(-1)^dF^\top v$, and \eqref{eq:besselG} with \eqref{eq:S1} give
\begin{equation}\label{eq:quadlower}
v^*K_Xv\ge(2\pi)^{-2d}4^{-l-d}\,\|F^\top v\|_2^2\ge\vartheta\|v\|_2^2,\qquad\vartheta:=(2\pi)^{-2d}4^{-l-d}\,\frac M2\qquad(v\in\cS).
\end{equation}

Since $K_X$ is real symmetric, its real and complex spectra
coincide. By the min-max principle, \eqref{eq:quadlower} gives at
least $(1-\delta)N$ eigenvalues at least $\vartheta$.
The orthogonal projection $E$ onto their real eigenspaces therefore
satisfies $\rank E\ge(1-\delta)N$ and $K_X\succeq\vartheta E$.
Finally, by \eqref{eq:Mbounds}, \eqref{eq:2l} and $l>\frac12\log_2N$,
\[
\vartheta\ge(2\pi)^{-2d}4^{-d}\,\frac{32N}{\delta}\,\Big(\frac{\delta\,l^q}{128d\,q!\,N}\Big)^2\ge c_d\,\delta\,\frac{l^{2q}}N\ge b_{d,\delta}\frac{(\log N)^{2q}}N,
\]
which proves \eqref{eq:operator} on the spectral event.
\end{proof}

\subsection{The fine part of a query}\label{sec:fine}

A dyadic interval of order $r\in\Z_{\ge0}$ is $[m2^{-r},(m+1)2^{-r})$ with $m\in\{0,\dots,2^r-1\}$. A dyadic box is a product of these intervals, and its order is the sum of their orders. For $t\in[0,1)$ write $t=\sum_{r\ge1}\eps_r2^{-r}$ with $\eps_r\in\{0,1\}$, choosing the terminating expansion at dyadic rationals. For each $r$ with $\eps_r=1$ put
\[
I_r(t)=\Big[\sum_{s<r}\eps_s2^{-s},\ \sum_{s\le r}\eps_s2^{-s}\Big),
\]
a dyadic interval of order $r$, since its left endpoint is a multiple of $2^{-r}$ and its length is $2^{-r}$. These intervals are pairwise disjoint and their union is $[0,t)$. For $t=1$ we use the single interval $[0,1)$ of order $0$, and for $t=0$ the decomposition is empty. If $x\in I_r(t)$, then
\begin{equation}\label{eq:prefix}
0<t-x\le2^{1-r},
\end{equation}
since $t-x\le t-\sum_{s<r}\eps_s2^{-s}=\sum_{s\ge r}\eps_s2^{-s}\le2^{1-r}$, and $t-x>t-\sum_{s\le r}\eps_s2^{-s}\ge0$. When $t=1$, the bound reads $0<1-x\le2$.

For $t\in[0,1]^d$ the products $\prod_{a=1}^dI_{r_a}(t_a)$ of intervals of the decompositions of the coordinates partition $[0,t)$ into dyadic boxes. For an integer $h\ge1$ let $S_t$ be the union of those products whose order $r_1+\cdots+r_d$ exceeds $h$, and put
\[
z(t)=\big(\one\{x_i\in S_t\}\big)_{i=1}^N.
\]
Let $W$ be the span in $\R^N$ of the incidence vectors $(\one\{x_i\in Q\})_i$ of all dyadic boxes $Q$ of order at most $h$. The vector $a(t)-z(t)$ is the sum of the incidence vectors of the dyadic boxes of order at most $h$ in the decomposition of $[0,t)$, so it lies in $W$, and consequently
\begin{equation}\label{eq:Paz}
Pa(t)=Pz(t)\qquad(t\in[0,1]^d)
\end{equation}
for every orthogonal projection $P$ whose range is contained in $W^\perp$.
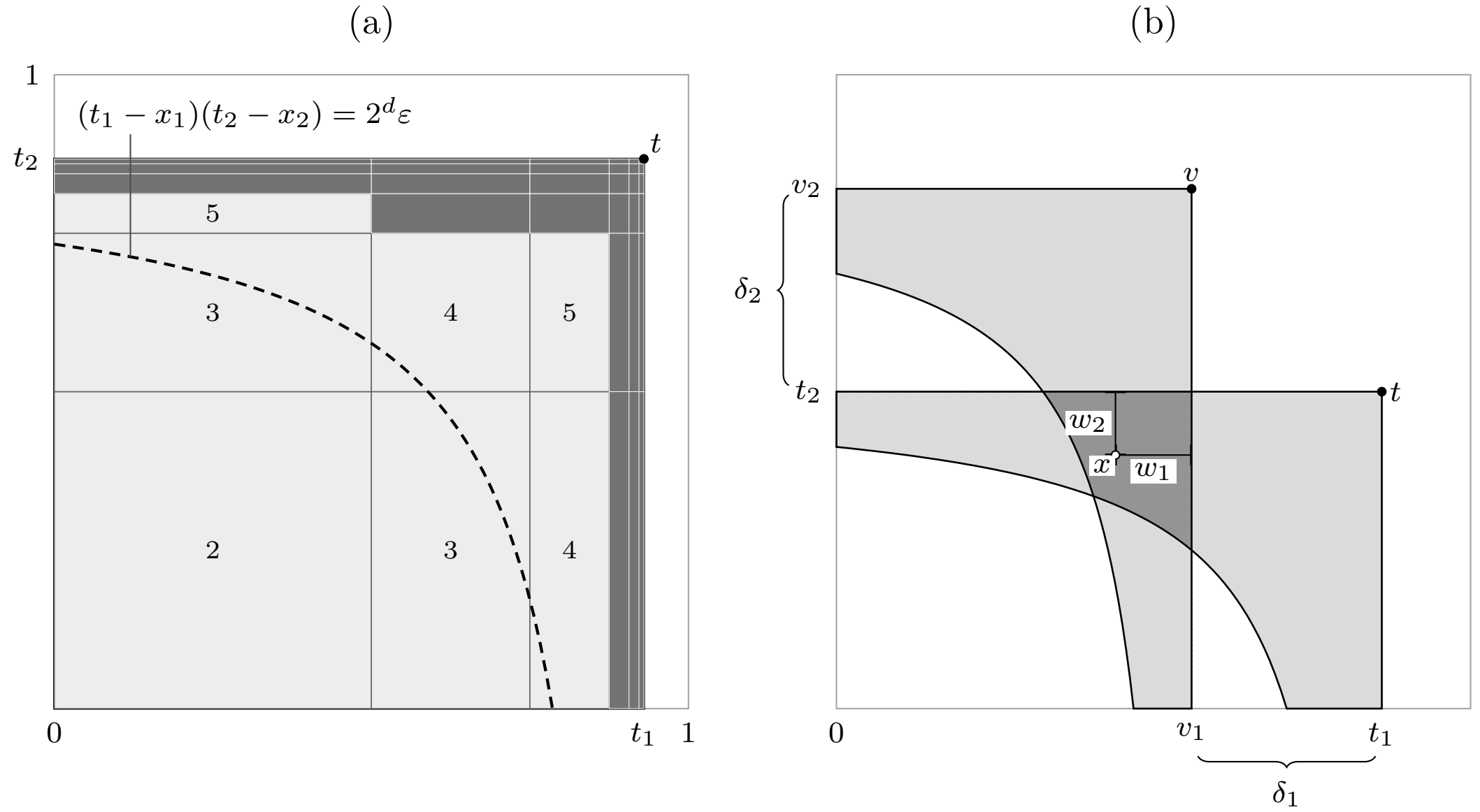
\begin{figure}[t]
\centering
\begin{tikzpicture}[x=6cm,y=6cm]
\begin{scope}
\draw[black!35] (0,0) rectangle (1,1);
\draw[semithick] (0,0) rectangle (0.9297,0.8672);
\draw[black!60, line width=0.3pt, fill=black!7] (0,0) rectangle (0.5,0.5);
\node[font=\tiny] at (0.25,0.25) {$2$};
\draw[black!60, line width=0.3pt, fill=black!7] (0,0.5) rectangle (0.5,0.75);
\node[font=\tiny] at (0.25,0.625) {$3$};
\draw[black!60, line width=0.3pt, fill=black!7] (0,0.75) rectangle (0.5,0.8125);
\node[font=\tiny] at (0.25,0.7812) {$5$};
\fill[black!55] (0,0.8125) rectangle (0.5,0.8438);
\draw[white, line width=0.15pt] (0,0.8125) rectangle (0.5,0.8438);
\fill[black!55] (0,0.8438) rectangle (0.5,0.8594);
\draw[white, line width=0.15pt] (0,0.8438) rectangle (0.5,0.8594);
\fill[black!55] (0,0.8594) rectangle (0.5,0.8672);
\draw[white, line width=0.15pt] (0,0.8594) rectangle (0.5,0.8672);
\draw[black!60, line width=0.3pt, fill=black!7] (0.5,0) rectangle (0.75,0.5);
\node[font=\tiny] at (0.625,0.25) {$3$};
\draw[black!60, line width=0.3pt, fill=black!7] (0.5,0.5) rectangle (0.75,0.75);
\node[font=\tiny] at (0.625,0.625) {$4$};
\fill[black!55] (0.5,0.75) rectangle (0.75,0.8125);
\draw[white, line width=0.15pt] (0.5,0.75) rectangle (0.75,0.8125);
\fill[black!55] (0.5,0.8125) rectangle (0.75,0.8438);
\draw[white, line width=0.15pt] (0.5,0.8125) rectangle (0.75,0.8438);
\fill[black!55] (0.5,0.8438) rectangle (0.75,0.8594);
\draw[white, line width=0.15pt] (0.5,0.8438) rectangle (0.75,0.8594);
\fill[black!55] (0.5,0.8594) rectangle (0.75,0.8672);
\draw[white, line width=0.15pt] (0.5,0.8594) rectangle (0.75,0.8672);
\draw[black!60, line width=0.3pt, fill=black!7] (0.75,0) rectangle (0.875,0.5);
\node[font=\tiny] at (0.8125,0.25) {$4$};
\draw[black!60, line width=0.3pt, fill=black!7] (0.75,0.5) rectangle (0.875,0.75);
\node[font=\tiny] at (0.8125,0.625) {$5$};
\fill[black!55] (0.75,0.75) rectangle (0.875,0.8125);
\draw[white, line width=0.15pt] (0.75,0.75) rectangle (0.875,0.8125);
\fill[black!55] (0.75,0.8125) rectangle (0.875,0.8438);
\draw[white, line width=0.15pt] (0.75,0.8125) rectangle (0.875,0.8438);
\fill[black!55] (0.75,0.8438) rectangle (0.875,0.8594);
\draw[white, line width=0.15pt] (0.75,0.8438) rectangle (0.875,0.8594);
\fill[black!55] (0.75,0.8594) rectangle (0.875,0.8672);
\draw[white, line width=0.15pt] (0.75,0.8594) rectangle (0.875,0.8672);
\fill[black!55] (0.875,0) rectangle (0.9062,0.5);
\draw[white, line width=0.15pt] (0.875,0) rectangle (0.9062,0.5);
\fill[black!55] (0.875,0.5) rectangle (0.9062,0.75);
\draw[white, line width=0.15pt] (0.875,0.5) rectangle (0.9062,0.75);
\fill[black!55] (0.875,0.75) rectangle (0.9062,0.8125);
\draw[white, line width=0.15pt] (0.875,0.75) rectangle (0.9062,0.8125);
\fill[black!55] (0.875,0.8125) rectangle (0.9062,0.8438);
\draw[white, line width=0.15pt] (0.875,0.8125) rectangle (0.9062,0.8438);
\fill[black!55] (0.875,0.8438) rectangle (0.9062,0.8594);
\draw[white, line width=0.15pt] (0.875,0.8438) rectangle (0.9062,0.8594);
\fill[black!55] (0.875,0.8594) rectangle (0.9062,0.8672);
\draw[white, line width=0.15pt] (0.875,0.8594) rectangle (0.9062,0.8672);
\fill[black!55] (0.9062,0) rectangle (0.9219,0.5);
\draw[white, line width=0.15pt] (0.9062,0) rectangle (0.9219,0.5);
\fill[black!55] (0.9062,0.5) rectangle (0.9219,0.75);
\draw[white, line width=0.15pt] (0.9062,0.5) rectangle (0.9219,0.75);
\fill[black!55] (0.9062,0.75) rectangle (0.9219,0.8125);
\draw[white, line width=0.15pt] (0.9062,0.75) rectangle (0.9219,0.8125);
\fill[black!55] (0.9062,0.8125) rectangle (0.9219,0.8438);
\draw[white, line width=0.15pt] (0.9062,0.8125) rectangle (0.9219,0.8438);
\fill[black!55] (0.9062,0.8438) rectangle (0.9219,0.8594);
\draw[white, line width=0.15pt] (0.9062,0.8438) rectangle (0.9219,0.8594);
\fill[black!55] (0.9062,0.8594) rectangle (0.9219,0.8672);
\draw[white, line width=0.15pt] (0.9062,0.8594) rectangle (0.9219,0.8672);
\fill[black!55] (0.9219,0) rectangle (0.9297,0.5);
\draw[white, line width=0.15pt] (0.9219,0) rectangle (0.9297,0.5);
\fill[black!55] (0.9219,0.5) rectangle (0.9297,0.75);
\draw[white, line width=0.15pt] (0.9219,0.5) rectangle (0.9297,0.75);
\fill[black!55] (0.9219,0.75) rectangle (0.9297,0.8125);
\draw[white, line width=0.15pt] (0.9219,0.75) rectangle (0.9297,0.8125);
\fill[black!55] (0.9219,0.8125) rectangle (0.9297,0.8438);
\draw[white, line width=0.15pt] (0.9219,0.8125) rectangle (0.9297,0.8438);
\fill[black!55] (0.9219,0.8438) rectangle (0.9297,0.8594);
\draw[white, line width=0.15pt] (0.9219,0.8438) rectangle (0.9297,0.8594);
\fill[black!55] (0.9219,0.8594) rectangle (0.9297,0.8672);
\draw[white, line width=0.15pt] (0.9219,0.8594) rectangle (0.9297,0.8672);
\draw[densely dashed, thick] (0,0.7327) -- (0.0098,0.7313) -- (0.0196,0.7298) -- (0.0295,0.7283) -- (0.0393,0.7268) -- (0.0491,0.7252) -- (0.0589,0.7236) -- (0.0687,0.722) -- (0.0786,0.7203) -- (0.0884,0.7186) -- (0.0982,0.7169) -- (0.108,0.7151) -- (0.1178,0.7132) -- (0.1277,0.7113) -- (0.1375,0.7094) -- (0.1473,0.7074) -- (0.1571,0.7054) -- (0.1669,0.7033) -- (0.1767,0.7012) -- (0.1866,0.699) -- (0.1964,0.6967) -- (0.2062,0.6944) -- (0.216,0.692) -- (0.2258,0.6896) -- (0.2357,0.6871) -- (0.2455,0.6845) -- (0.2553,0.6818) -- (0.2651,0.6791) -- (0.2749,0.6763) -- (0.2848,0.6734) -- (0.2946,0.6704) -- (0.3044,0.6673) -- (0.3142,0.6641) -- (0.324,0.6608) -- (0.3339,0.6574) -- (0.3437,0.6539) -- (0.3535,0.6502) -- (0.3633,0.6465) -- (0.3731,0.6426) -- (0.383,0.6386) -- (0.3928,0.6344) -- (0.4026,0.63) -- (0.4124,0.6255) -- (0.4222,0.6209) -- (0.432,0.616) -- (0.4419,0.6109) -- (0.4517,0.6057) -- (0.4615,0.6002) -- (0.4713,0.5945) -- (0.4811,0.5885) -- (0.491,0.5823) -- (0.5008,0.5757) -- (0.5106,0.5689) -- (0.5204,0.5618) -- (0.5302,0.5543) -- (0.5401,0.5464) -- (0.5499,0.5381) -- (0.5597,0.5293) -- (0.5695,0.5201) -- (0.5793,0.5104) -- (0.5892,0.5001) -- (0.599,0.4892) -- (0.6088,0.4776) -- (0.6186,0.4654) -- (0.6284,0.4523) -- (0.6383,0.4383) -- (0.6481,0.4233) -- (0.6579,0.4073) -- (0.6677,0.39) -- (0.6775,0.3715) -- (0.6874,0.3514) -- (0.6972,0.3296) -- (0.707,0.3059) -- (0.7168,0.28) -- (0.7266,0.2516) -- (0.7364,0.2203) -- (0.7463,0.1857) -- (0.7561,0.1471) -- (0.7659,0.104) -- (0.7757,0.0553) -- (0.7855,0);
\fill (0.9297,0.8672) circle (1.3pt);
\node[font=\scriptsize, above right, inner sep=1.5pt] at (0.9297,0.8672) {$t$};
\node[font=\scriptsize, below] at (0,0) {$0$};
\node[font=\scriptsize, below] at (0.9297,0) {$t_1$};
\node[font=\scriptsize, below] at (1,0) {$1$};
\node[font=\scriptsize, left] at (0,0.8672) {$t_2$};
\node[font=\scriptsize, left] at (0,1) {$1$};
\draw[black!70, thin] (0.12,0.7128) -- (0.12,0.915);
\node[font=\scriptsize, above, fill=white, inner sep=1pt] at (0.3,0.905) {$(t_1-x_1)(t_2-x_2)=2^d\varepsilon$};
\node[font=\small] at (0.5,1.08) {(a)};
\end{scope}
\begin{scope}[xshift=7.4cm]
\draw[black!35] (0,0) rectangle (1,1);
\fill[black!14] (0,0.5) -- (0.86,0.5) -- (0.86,0) -- (0.71,0) -- (0.7011,0.0279) -- (0.6922,0.0529) -- (0.6834,0.0754) -- (0.6745,0.0957) -- (0.6656,0.1141) -- (0.6567,0.131) -- (0.6479,0.1464) -- (0.639,0.1606) -- (0.6301,0.1737) -- (0.6212,0.1859) -- (0.6124,0.1971) -- (0.6035,0.2076) -- (0.5946,0.2174) -- (0.5857,0.2265) -- (0.5769,0.2351) -- (0.568,0.2432) -- (0.5591,0.2507) -- (0.5502,0.2579) -- (0.5414,0.2646) -- (0.5325,0.271) -- (0.5236,0.277) -- (0.5148,0.2828) -- (0.5059,0.2882) -- (0.497,0.2934) -- (0.4881,0.2983) -- (0.4792,0.303) -- (0.4704,0.3075) -- (0.4615,0.3118) -- (0.4526,0.3159) -- (0.4437,0.3198) -- (0.4349,0.3236) -- (0.426,0.3272) -- (0.4171,0.3307) -- (0.4082,0.334) -- (0.3994,0.3372) -- (0.3905,0.3403) -- (0.3816,0.3432) -- (0.3728,0.3461) -- (0.3639,0.3488) -- (0.355,0.3515) -- (0.3461,0.3541) -- (0.3372,0.3565) -- (0.3284,0.3589) -- (0.3195,0.3612) -- (0.3106,0.3635) -- (0.3018,0.3657) -- (0.2929,0.3678) -- (0.284,0.3698) -- (0.2751,0.3718) -- (0.2662,0.3737) -- (0.2574,0.3755) -- (0.2485,0.3774) -- (0.2396,0.3791) -- (0.2308,0.3808) -- (0.2219,0.3825) -- (0.213,0.3841) -- (0.2041,0.3856) -- (0.1952,0.3872) -- (0.1864,0.3887) -- (0.1775,0.3901) -- (0.1686,0.3915) -- (0.1598,0.3929) -- (0.1509,0.3942) -- (0.142,0.3955) -- (0.1331,0.3968) -- (0.1242,0.3981) -- (0.1154,0.3993) -- (0.1065,0.4005) -- (0.0976,0.4016) -- (0.0887,0.4028) -- (0.0799,0.4039) -- (0.071,0.4049) -- (0.0621,0.406) -- (0.0532,0.407) -- (0.0444,0.408) -- (0.0355,0.409) -- (0.0266,0.41) -- (0.0177,0.411) -- (0.0089,0.4119) -- (0,0.4128) -- cycle;
\fill[black!14] (0,0.82) -- (0.56,0.82) -- (0.56,0) -- (0.4685,0) -- (0.4627,0.0493) -- (0.4568,0.0931) -- (0.451,0.1321) -- (0.4451,0.1672) -- (0.4393,0.1989) -- (0.4334,0.2276) -- (0.4275,0.2538) -- (0.4217,0.2778) -- (0.4158,0.2998) -- (0.41,0.3201) -- (0.4041,0.3389) -- (0.3983,0.3563) -- (0.3924,0.3725) -- (0.3865,0.3876) -- (0.3807,0.4017) -- (0.3748,0.415) -- (0.369,0.4274) -- (0.3631,0.4391) -- (0.3573,0.4501) -- (0.3514,0.4605) -- (0.3455,0.4703) -- (0.3397,0.4796) -- (0.3338,0.4884) -- (0.328,0.4968) -- (0.3221,0.5047) -- (0.3163,0.5123) -- (0.3104,0.5195) -- (0.3045,0.5264) -- (0.2987,0.533) -- (0.2928,0.5393) -- (0.287,0.5453) -- (0.2811,0.5511) -- (0.2753,0.5566) -- (0.2694,0.5619) -- (0.2636,0.567) -- (0.2577,0.5719) -- (0.2518,0.5766) -- (0.246,0.5812) -- (0.2401,0.5855) -- (0.2343,0.5897) -- (0.2284,0.5938) -- (0.2226,0.5977) -- (0.2167,0.6015) -- (0.2108,0.6052) -- (0.205,0.6087) -- (0.1991,0.6122) -- (0.1933,0.6155) -- (0.1874,0.6187) -- (0.1816,0.6218) -- (0.1757,0.6248) -- (0.1698,0.6278) -- (0.164,0.6306) -- (0.1581,0.6334) -- (0.1523,0.6361) -- (0.1464,0.6387) -- (0.1406,0.6412) -- (0.1347,0.6437) -- (0.1288,0.646) -- (0.123,0.6484) -- (0.1171,0.6506) -- (0.1113,0.6529) -- (0.1054,0.655) -- (0.0996,0.6571) -- (0.0937,0.6592) -- (0.0879,0.6612) -- (0.082,0.6631) -- (0.0761,0.665) -- (0.0703,0.6669) -- (0.0644,0.6687) -- (0.0586,0.6704) -- (0.0527,0.6722) -- (0.0469,0.6738) -- (0.041,0.6755) -- (0.0351,0.6771) -- (0.0293,0.6787) -- (0.0234,0.6802) -- (0.0176,0.6817) -- (0.0117,0.6832) -- (0.0059,0.6847) -- (0,0.6861) -- cycle;
\begin{scope}
\clip (0,0.5) -- (0.86,0.5) -- (0.86,0) -- (0.71,0) -- (0.7011,0.0279) -- (0.6922,0.0529) -- (0.6834,0.0754) -- (0.6745,0.0957) -- (0.6656,0.1141) -- (0.6567,0.131) -- (0.6479,0.1464) -- (0.639,0.1606) -- (0.6301,0.1737) -- (0.6212,0.1859) -- (0.6124,0.1971) -- (0.6035,0.2076) -- (0.5946,0.2174) -- (0.5857,0.2265) -- (0.5769,0.2351) -- (0.568,0.2432) -- (0.5591,0.2507) -- (0.5502,0.2579) -- (0.5414,0.2646) -- (0.5325,0.271) -- (0.5236,0.277) -- (0.5148,0.2828) -- (0.5059,0.2882) -- (0.497,0.2934) -- (0.4881,0.2983) -- (0.4792,0.303) -- (0.4704,0.3075) -- (0.4615,0.3118) -- (0.4526,0.3159) -- (0.4437,0.3198) -- (0.4349,0.3236) -- (0.426,0.3272) -- (0.4171,0.3307) -- (0.4082,0.334) -- (0.3994,0.3372) -- (0.3905,0.3403) -- (0.3816,0.3432) -- (0.3728,0.3461) -- (0.3639,0.3488) -- (0.355,0.3515) -- (0.3461,0.3541) -- (0.3372,0.3565) -- (0.3284,0.3589) -- (0.3195,0.3612) -- (0.3106,0.3635) -- (0.3018,0.3657) -- (0.2929,0.3678) -- (0.284,0.3698) -- (0.2751,0.3718) -- (0.2662,0.3737) -- (0.2574,0.3755) -- (0.2485,0.3774) -- (0.2396,0.3791) -- (0.2308,0.3808) -- (0.2219,0.3825) -- (0.213,0.3841) -- (0.2041,0.3856) -- (0.1952,0.3872) -- (0.1864,0.3887) -- (0.1775,0.3901) -- (0.1686,0.3915) -- (0.1598,0.3929) -- (0.1509,0.3942) -- (0.142,0.3955) -- (0.1331,0.3968) -- (0.1242,0.3981) -- (0.1154,0.3993) -- (0.1065,0.4005) -- (0.0976,0.4016) -- (0.0887,0.4028) -- (0.0799,0.4039) -- (0.071,0.4049) -- (0.0621,0.406) -- (0.0532,0.407) -- (0.0444,0.408) -- (0.0355,0.409) -- (0.0266,0.41) -- (0.0177,0.411) -- (0.0089,0.4119) -- (0,0.4128) -- cycle;
\fill[black!42] (0,0.82) -- (0.56,0.82) -- (0.56,0) -- (0.4685,0) -- (0.4627,0.0493) -- (0.4568,0.0931) -- (0.451,0.1321) -- (0.4451,0.1672) -- (0.4393,0.1989) -- (0.4334,0.2276) -- (0.4275,0.2538) -- (0.4217,0.2778) -- (0.4158,0.2998) -- (0.41,0.3201) -- (0.4041,0.3389) -- (0.3983,0.3563) -- (0.3924,0.3725) -- (0.3865,0.3876) -- (0.3807,0.4017) -- (0.3748,0.415) -- (0.369,0.4274) -- (0.3631,0.4391) -- (0.3573,0.4501) -- (0.3514,0.4605) -- (0.3455,0.4703) -- (0.3397,0.4796) -- (0.3338,0.4884) -- (0.328,0.4968) -- (0.3221,0.5047) -- (0.3163,0.5123) -- (0.3104,0.5195) -- (0.3045,0.5264) -- (0.2987,0.533) -- (0.2928,0.5393) -- (0.287,0.5453) -- (0.2811,0.5511) -- (0.2753,0.5566) -- (0.2694,0.5619) -- (0.2636,0.567) -- (0.2577,0.5719) -- (0.2518,0.5766) -- (0.246,0.5812) -- (0.2401,0.5855) -- (0.2343,0.5897) -- (0.2284,0.5938) -- (0.2226,0.5977) -- (0.2167,0.6015) -- (0.2108,0.6052) -- (0.205,0.6087) -- (0.1991,0.6122) -- (0.1933,0.6155) -- (0.1874,0.6187) -- (0.1816,0.6218) -- (0.1757,0.6248) -- (0.1698,0.6278) -- (0.164,0.6306) -- (0.1581,0.6334) -- (0.1523,0.6361) -- (0.1464,0.6387) -- (0.1406,0.6412) -- (0.1347,0.6437) -- (0.1288,0.646) -- (0.123,0.6484) -- (0.1171,0.6506) -- (0.1113,0.6529) -- (0.1054,0.655) -- (0.0996,0.6571) -- (0.0937,0.6592) -- (0.0879,0.6612) -- (0.082,0.6631) -- (0.0761,0.665) -- (0.0703,0.6669) -- (0.0644,0.6687) -- (0.0586,0.6704) -- (0.0527,0.6722) -- (0.0469,0.6738) -- (0.041,0.6755) -- (0.0351,0.6771) -- (0.0293,0.6787) -- (0.0234,0.6802) -- (0.0176,0.6817) -- (0.0117,0.6832) -- (0.0059,0.6847) -- (0,0.6861) -- cycle;
\end{scope}
\draw[line width=0.5pt] (0,0.5) -- (0.86,0.5) -- (0.86,0) -- (0.71,0) -- (0.7011,0.0279) -- (0.6922,0.0529) -- (0.6834,0.0754) -- (0.6745,0.0957) -- (0.6656,0.1141) -- (0.6567,0.131) -- (0.6479,0.1464) -- (0.639,0.1606) -- (0.6301,0.1737) -- (0.6212,0.1859) -- (0.6124,0.1971) -- (0.6035,0.2076) -- (0.5946,0.2174) -- (0.5857,0.2265) -- (0.5769,0.2351) -- (0.568,0.2432) -- (0.5591,0.2507) -- (0.5502,0.2579) -- (0.5414,0.2646) -- (0.5325,0.271) -- (0.5236,0.277) -- (0.5148,0.2828) -- (0.5059,0.2882) -- (0.497,0.2934) -- (0.4881,0.2983) -- (0.4792,0.303) -- (0.4704,0.3075) -- (0.4615,0.3118) -- (0.4526,0.3159) -- (0.4437,0.3198) -- (0.4349,0.3236) -- (0.426,0.3272) -- (0.4171,0.3307) -- (0.4082,0.334) -- (0.3994,0.3372) -- (0.3905,0.3403) -- (0.3816,0.3432) -- (0.3728,0.3461) -- (0.3639,0.3488) -- (0.355,0.3515) -- (0.3461,0.3541) -- (0.3372,0.3565) -- (0.3284,0.3589) -- (0.3195,0.3612) -- (0.3106,0.3635) -- (0.3018,0.3657) -- (0.2929,0.3678) -- (0.284,0.3698) -- (0.2751,0.3718) -- (0.2662,0.3737) -- (0.2574,0.3755) -- (0.2485,0.3774) -- (0.2396,0.3791) -- (0.2308,0.3808) -- (0.2219,0.3825) -- (0.213,0.3841) -- (0.2041,0.3856) -- (0.1952,0.3872) -- (0.1864,0.3887) -- (0.1775,0.3901) -- (0.1686,0.3915) -- (0.1598,0.3929) -- (0.1509,0.3942) -- (0.142,0.3955) -- (0.1331,0.3968) -- (0.1242,0.3981) -- (0.1154,0.3993) -- (0.1065,0.4005) -- (0.0976,0.4016) -- (0.0887,0.4028) -- (0.0799,0.4039) -- (0.071,0.4049) -- (0.0621,0.406) -- (0.0532,0.407) -- (0.0444,0.408) -- (0.0355,0.409) -- (0.0266,0.41) -- (0.0177,0.411) -- (0.0089,0.4119) -- (0,0.4128) -- cycle;
\draw[line width=0.5pt] (0,0.82) -- (0.56,0.82) -- (0.56,0) -- (0.4685,0) -- (0.4627,0.0493) -- (0.4568,0.0931) -- (0.451,0.1321) -- (0.4451,0.1672) -- (0.4393,0.1989) -- (0.4334,0.2276) -- (0.4275,0.2538) -- (0.4217,0.2778) -- (0.4158,0.2998) -- (0.41,0.3201) -- (0.4041,0.3389) -- (0.3983,0.3563) -- (0.3924,0.3725) -- (0.3865,0.3876) -- (0.3807,0.4017) -- (0.3748,0.415) -- (0.369,0.4274) -- (0.3631,0.4391) -- (0.3573,0.4501) -- (0.3514,0.4605) -- (0.3455,0.4703) -- (0.3397,0.4796) -- (0.3338,0.4884) -- (0.328,0.4968) -- (0.3221,0.5047) -- (0.3163,0.5123) -- (0.3104,0.5195) -- (0.3045,0.5264) -- (0.2987,0.533) -- (0.2928,0.5393) -- (0.287,0.5453) -- (0.2811,0.5511) -- (0.2753,0.5566) -- (0.2694,0.5619) -- (0.2636,0.567) -- (0.2577,0.5719) -- (0.2518,0.5766) -- (0.246,0.5812) -- (0.2401,0.5855) -- (0.2343,0.5897) -- (0.2284,0.5938) -- (0.2226,0.5977) -- (0.2167,0.6015) -- (0.2108,0.6052) -- (0.205,0.6087) -- (0.1991,0.6122) -- (0.1933,0.6155) -- (0.1874,0.6187) -- (0.1816,0.6218) -- (0.1757,0.6248) -- (0.1698,0.6278) -- (0.164,0.6306) -- (0.1581,0.6334) -- (0.1523,0.6361) -- (0.1464,0.6387) -- (0.1406,0.6412) -- (0.1347,0.6437) -- (0.1288,0.646) -- (0.123,0.6484) -- (0.1171,0.6506) -- (0.1113,0.6529) -- (0.1054,0.655) -- (0.0996,0.6571) -- (0.0937,0.6592) -- (0.0879,0.6612) -- (0.082,0.6631) -- (0.0761,0.665) -- (0.0703,0.6669) -- (0.0644,0.6687) -- (0.0586,0.6704) -- (0.0527,0.6722) -- (0.0469,0.6738) -- (0.041,0.6755) -- (0.0351,0.6771) -- (0.0293,0.6787) -- (0.0234,0.6802) -- (0.0176,0.6817) -- (0.0117,0.6832) -- (0.0059,0.6847) -- (0,0.6861) -- cycle;
\draw[densely dotted] (0.56,0) -- (0.56,0.5) -- (0,0.5);
\fill (0.86,0.5) circle (1.3pt);
\node[font=\scriptsize, right, inner sep=2pt] at (0.86,0.5) {$t$};
\fill (0.56,0.82) circle (1.3pt);
\node[font=\scriptsize, above, inner sep=2pt] at (0.56,0.82) {$v$};
\draw[|-|, line width=0.4pt] (0.44,0.4) -- (0.56,0.4);
\draw[|-|, line width=0.4pt] (0.44,0.4) -- (0.44,0.5);
\fill[white] (0.44,0.4) circle (1.1pt); \draw (0.44,0.4) circle (1.1pt);
\node[font=\scriptsize, below left, inner sep=1pt, fill=white] at (0.44,0.4) {$x$};
\node[font=\scriptsize, below=1pt, inner sep=1pt, fill=white] at (0.5,0.4) {$w_1$};
\node[font=\scriptsize, left=1pt, inner sep=1pt, fill=white] at (0.44,0.45) {$w_2$};
\node[font=\scriptsize, below] at (0,0) {$0$};
\node[font=\scriptsize, below] at (0.56,0) {$v_1$};
\node[font=\scriptsize, below] at (0.86,0) {$t_1$};
\node[font=\scriptsize, left] at (0,0.5) {$t_2$};
\node[font=\scriptsize, left] at (0,0.82) {$v_2$};
\draw[decorate, decoration={brace, mirror, amplitude=3pt}] (0.57,-0.075) -- (0.85,-0.075) node[midway, below=3pt, font=\scriptsize] {$\delta_1$};
\draw[decorate, decoration={brace, amplitude=3pt}] (-0.075,0.51) -- (-0.075,0.81) node[midway, left=3pt, font=\scriptsize] {$\delta_2$};
\node[font=\small] at (0.5,1.08) {(b)};
\end{scope}
\end{tikzpicture}
\caption{The fine part of an anchored box for $d=2$. (a) The dyadic boxes $I_{r_1}(t_1)\times I_{r_2}(t_2)$ partition $[0,t)$, and the numbers are their orders $r_1+r_2$. The incidence vectors of the light boxes, of order at most $h=5$, lie in $W$. The dark boxes, of order larger than $h$, form $S_t$, which lies between the dashed curve and the corner $t$, that is, in $E_t(2^d\varepsilon)$, as in Lemma~\ref{lem:coarse}. (b) Two regions $E_t(\eta)$ and $E_v(\eta)$ and their intersection, dark. A point $x$ of the intersection lies below and to the left of $\min(t,v)$; its coordinatewise distances to $v$ are $w_1$ and $w_2+\delta_2$, and those to $t$ are $w_1+\delta_1$ and $w_2$, as in the proof of Lemma~\ref{lem:overlap}.}
\label{fig:fine}
\end{figure}
\begin{lemma}\label{lem:coarse}
With $\eps=2^{-h}$ and a constant $C_d$ depending only on $d$,
\begin{equation}\label{eq:dimW}
\dim W\le C_d2^hh^q,
\end{equation}
and for every $t\in[0,1]^d$,
\begin{equation}\label{eq:layer}
S_t\subseteq E_t(2^d\eps),\qquad E_t(\eta)=\Big\{x\in[0,1]^d:\ x<t\text{ coordinatewise},\ \prod_{a=1}^d(t_a-x_a)\le\eta\Big\}.
\end{equation}
\end{lemma}

\begin{proof}
For a fixed vector of orders $(r_1,\dots,r_d)$ with sum $\ell$ there are $2^\ell$ dyadic boxes, and there are $\binom{\ell+q}q$ such vectors. Hence
\[
\dim W\le\sum_{\ell=0}^h2^\ell\binom{\ell+q}q\le2^{h+1}\binom{h+q}q\le2^{h+1}\frac{(h+q)^q}{q!}\le\frac{2(1+q)^q}{q!}\,2^hh^q.
\]
If $x$ lies in a product $\prod_aI_{r_a}(t_a)$ of order $\ell=\sum_ar_a>h$, then $x<t$ coordinatewise and, by \eqref{eq:prefix}, $\prod_a(t_a-x_a)\le\prod_a2^{1-r_a}=2^{d-\ell}\le2^d\eps$. Taking the union over all such products proves \eqref{eq:layer}.
\end{proof}

Fix $\alpha\in(0,1)$ and increase the constant $C_d$ in
\eqref{eq:dimW} so that $C_d\ge1$.
For sufficiently large $N$, choose $h\ge1$ maximal with
\begin{equation}\label{eq:hdef}
C_d2^hh^q\le\frac{\alpha N}8.
\end{equation}
Then
\begin{equation}\label{eq:dimWalpha}
\dim W\le\frac{\alpha N}8.
\end{equation}
Maximality gives $N<8C_d2^{h+1}(h+1)^q/\alpha$.
Together with $2^h\le N$, this implies
\begin{equation}\label{eq:hlog}
\tfrac12\log_2N\le h\le\log_2N
\qquad(N\ge N_{d,\alpha}).
\end{equation}
Write
\begin{equation}\label{eq:Leps}
L=h^q,\qquad\eps=2^{-h}.
\end{equation}
The choice of $h$ and its maximality also give
\begin{equation}\label{eq:Neps}
N\eps\asymp_{d,\alpha}L.
\end{equation}

\subsection{Overlap moments}\label{sec:overlap}

The geometric input is an estimate of the overlaps of the sets $S_t$ which is uniform in a growing moment order. Put
\begin{equation}\label{eq:Rj}
R_j=\int_{[0,1]^d}\int_{[0,1]^d}\vol(S_t\cap S_v)^j\,dt\,dv\qquad(j\ge1).
\end{equation}

\begin{lemma}\label{lem:product}
For $\eta>0$ let
\[
J_d(\eta)=\int_{\{z\in(0,\infty)^d:\ \prod_az_a\le\eta\}}\ \prod_{a=1}^d\min(1,z_a)\,dz.
\]
Then there is a constant $C_d$ with
\begin{equation}\label{eq:Jbound}
J_d(\eta)\le C_d\,\eta\min(1,\eta)\,(1+|\log\eta|)^q.
\end{equation}
\end{lemma}

\begin{proof}
Substitute $z_a=e^{y_a}$. Since
$\min(1,e^y)e^y=e^{3y/2-|y|/2}$,
\[
J_d(\eta)=4^d\int_{-\infty}^{\log\eta}e^{3v/2}\varphi_d(v)\,dv,
\]
where $\varphi_d$ is the density of the sum of $d$ independent variables
with density $\frac14e^{-|y|/2}$.
Each such variable is the difference of two independent exponential
variables of rate $1/2$, so $\varphi_d$ is the density of the difference
of two independent $\Gam(d,1/2)$ variables, with the second parameter
denoting the rate. It is even, and for $v\ge0$,
\[
\varphi_d(v)
=\frac{e^{-v/2}}{2^{2d}(q!)^2}
  \int_0^\infty e^{-r}r^q(r+v)^q\,dr
\le C_d e^{-v/2}(1+v)^q.
\]
Consequently,
\[
J_d(\eta)
\le C_d\int_{-\infty}^{\log\eta}
e^{v+\min(v,0)}(1+|v|)^q\,dv
\le C_d\,\eta\min(1,\eta)(1+|\log\eta|)^q.
\]
The last bound follows by integrating over the negative and positive
half-lines. For $\eta\le1$, the substitution $v=\log\eta-r$ bounds
the integral by
$\eta^2(1+|\log\eta|)^q\int_0^\infty e^{-2r}(1+r)^q\,dr$.
For $\eta>1$, the negative part is bounded by a constant and the positive
part is at most $\eta(1+\log\eta)^q$.
\end{proof}

\begin{lemma}\label{lem:overlap}
There are $h_d$ and a constant $C_d$ such that for all $h\ge h_d$
\begin{equation}\label{eq:R1}
R_1\le C_d\,\eps^2L^2,
\end{equation}
and, uniformly for integers $2\le j\le h$,
\begin{equation}\label{eq:Rj_bound}
R_j\le C_d^{\,j}\,\eps^{j+1}L\,j^{jq}.
\end{equation}
\end{lemma}

\begin{proof}
Put $\eta=2^d\eps$ and $H=\log(1/\eta)=(h-d)\log2$. For $h\ge h_d:=2d+4$ we have $0<\eta<1$, $H\ge(d+4)\log2\ge1$, and $\frac12h\log2\le H\le h\log2$. By \eqref{eq:layer} it suffices to estimate the overlaps of the sets $E_t(\eta)$.

For almost every pair $t,v\in[0,1]^d$, we have $\delta_a=|t_a-v_a|>0$ for every $a$. Put $c_\delta=\prod_a\delta_a$. If $x\in E_t(\eta)\cap E_v(\eta)$ then $x_a<\min(t_a,v_a)$ for every $a$, and with $w_a=\min(t_a,v_a)-x_a>0$ the two numbers $t_a-x_a$ and $v_a-x_a$ are $w_a$ and $w_a+\delta_a$ in some order. Multiplying the two product constraints gives
\begin{equation}\label{eq:wconstraint}
\prod_{a=1}^dw_a(w_a+\delta_a)\le\eta^2.
\end{equation}
This is a necessary condition for $x$ to lie in the intersection, and $x\mapsto w$ is a composition of a translation and a reflection, so
\[
\omega(t,v):=\vol\big(E_t(\eta)\cap E_v(\eta)\big)\le\vol\Big\{w\in(0,\infty)^d:\ \prod_aw_a(w_a+\delta_a)\le\eta^2\Big\}.
\]
For each $a$, the map $w_a\mapsto y_a=\sqrt{w_a(w_a+\delta_a)}$ is an increasing bijection of $(0,\infty)$ onto itself, with inverse $w_a=\frac12(\sqrt{\delta_a^2+4y_a^2}-\delta_a)$ and
\[
\frac{dw_a}{dy_a}=\frac{2y_a}{\sqrt{\delta_a^2+4y_a^2}}\le\min\Big(1,\frac{2y_a}{\delta_a}\Big).
\]
Hence, substituting also $y_a=\delta_az_a$ and using $\min(1,2z)\le2\min(1,z)$ for $z>0$,
\[
\omega(t,v)\le\int_{\{\prod_ay_a\le\eta\}}\prod_a\min\Big(1,\frac{2y_a}{\delta_a}\Big)dy=c_\delta\int_{\{\prod_az_a\le\eta/c_\delta\}}\prod_a\min(1,2z_a)\,dz\le2^dc_\delta J_d(\eta/c_\delta).
\]
By Lemma~\ref{lem:product},
\begin{equation}\label{eq:omegabound}
\omega(t,v)\le C_d\,\eta\min\Big(1,\frac\eta{c_\delta}\Big)\Big(1+\big|\log(c_\delta/\eta)\big|\Big)^q.
\end{equation}

Let $2\le j\le h$. For independent uniform $t,v$ the coordinates $\delta_a=|t_a-v_a|$ are independent with density $2(1-\delta_a)\le2$ on $(0,1)$. Since the right-hand side of \eqref{eq:omegabound} depends on $(t,v)$ only through $c_\delta$,
\[
R_j\le\iint\omega(t,v)^j\,dt\,dv\le2^d\int_{(0,1)^d}\Big[C_d\eta\min\Big(1,\frac\eta{c_\delta}\Big)\big(1+|\log(c_\delta/\eta)|\big)^q\Big]^j\,d\delta.
\]
Under Lebesgue measure on $(0,1)^d$, the variables $-\log\delta_a$ are independent exponential variables, their sum has the gamma density $\tau^qe^{-\tau}/q!$ on $(0,\infty)$, and therefore the product $c_\delta$ has density $(-\log c)^q/q!$ in the variable $c\in(0,1)$. Substituting $c=\eta e^{y}$, so that $dc=\eta e^{y}dy$, $y<H$, $\min(1,\eta/c)=e^{-\max(y,0)}$ and $-\log c=H-y$, and absorbing the factors $2^d$ and $1/q!$ into $C_d^{\,j}$, we obtain
\begin{equation}\label{eq:Rjint}
R_j\le C_d^{\,j}\,\eta^{j+1}\int_{-\infty}^He^{y-j\max(y,0)}(1+|y|)^{jq}(H-y)^q\,dy.
\end{equation}
On the positive half-axis,
$e^{-(j-1)y}\le e^{-y}$ because $j\ge2$, and
$(H-y)^q\le(H+y)^q$. Hence the integral in \eqref{eq:Rjint}
is at most
\begin{align*}
2\int_0^\infty e^{-r}(1+r)^{jq}(H+r)^q\,dr
&\le2^{q+1}e\big(H^q(jq)!+((j+1)q)!\big)\\
&\le C_dH^q(jq)!
\le C_d^{\,j}H^qj^{jq}.
\end{align*}
Here we used
$\int_0^\infty e^{-r}(1+r)^i\,dr\le e\,i!$,
$((j+1)q)!/(jq)!\le((j+1)q)^q\le C_dj^q$,
and $j\le h\le2H/\log2$.
Substituting this estimate into \eqref{eq:Rjint}, and using
$\eta=2^d\eps$ and $H^q\le L$, proves \eqref{eq:Rj_bound}.

For fixed $x\in[0,1]^d$ put $Q_x=\{t\in[0,1]^d:x\in S_t\}$. By \eqref{eq:layer}, $x\in S_t$ implies $t>x$ coordinatewise and $\prod_a(t_a-x_a)\le\eta$, so the substitution $y=t-x\in(0,1)^d$ gives
\[
\vol(Q_x)\le\vol\Big\{y\in(0,1)^d:\prod_ay_a\le\eta\Big\}=\frac1{q!}\int_H^\infty e^{-\tau}\tau^q\,d\tau=e^{-H}\sum_{i=0}^q\frac{H^i}{i!}\le\eta(1+H)^q,
\]
where we used the gamma density of $-\log\prod_ay_a$ and $1/i!\le\binom qi$. By Tonelli's theorem,
\[
R_1=\iint\int\one\{x\in S_t\}\one\{x\in S_v\}\,dx\,dt\,dv=\int_{[0,1]^d}\vol(Q_x)^2\,dx\le\eta^2(1+H)^{2q}\le C_d\eps^2L^2,
\]
which is \eqref{eq:R1}.
\end{proof}

\section{From overlap moments to discrepancy}\label{sec:disc}

\subsection{Summing the correlations}\label{sec:sum}

Fix $\alpha\in(0,1)$ and use $h,L,\eps,W$ from
Section~\ref{sec:fine}.
Let $\cS_N$ be the spectral event of Lemma~\ref{lem:sampled} with
$\delta=\alpha/8$. It is Borel, has probability at least $3/4$,
and on it there is an orthogonal projection $E$ of rank at least
$(1-\alpha/8)N$ with
\begin{equation}\label{eq:Eagain}
K_X\succeq b_{d,\alpha/8}\frac{(\log N)^{2q}}N\,E.
\end{equation}

Choose a sufficiently small fixed constant $c_{d,\alpha}>0$,
depending only on $d,\alpha$, and put
\begin{equation}\label{eq:sdef}
s=c_{d,\alpha}L.
\end{equation}
Fix $X\in\cS_N$ and $U\subseteq\R^N$ with $\dim U\ge\alpha N/2$.
Let $P$ project orthogonally onto $U\cap W^\perp$ and set
$\sigma(t)=\|Pa(t)\|_2$.
By \eqref{eq:dimWalpha}, $\rank P\ge3\alpha N/8$.
Since $0\preceq P\preceq I$,
$\tr(P(I-E))\le\tr(I-E)\le\alpha N/8$, and hence
$\tr(PE)\ge\alpha N/4$.
Taking the trace in \eqref{eq:Eagain} against $P$ gives
\begin{equation}\label{eq:S}
\int_{[0,1]^d}\sigma(t)^2\,dt
=\tr(PK_X)
\ge\frac\alpha4b_{d,\alpha/8}(\log N)^{2q}
\ge2s^2,
\end{equation}
where the last step follows from \eqref{eq:hlog} and the choice of
$c_{d,\alpha}$. Define
\begin{equation}\label{eq:mu}
G=\{t\in[0,1]^d:\sigma(t)\ge s\},\qquad
S_G=\int_G\sigma(t)^2\,dt,\qquad
d\mu(t)=S_G^{-1}\one_G(t)\sigma(t)^2\,dt.
\end{equation}
Since $\int_{G^c}\sigma(t)^2\,dt\le s^2$, \eqref{eq:S} gives
$S_G\ge s^2>0$. Thus $\mu$ is a probability measure and
$\sigma\ge s$ on $G$.
Put $u(t)=Pa(t)/\sigma(t)$ for $t\in G$ and
$r(t,v)=\inner{u(t)}{u(v)}$.
For fixed $X$ the maps $a$ and $z$ take finitely many values on
Borel sets, so all these integrals are well defined.

\begin{lemma}\label{lem:contraction}
For every even integer $j\ge2$, every $X\in\cS_N$ and every subspace $U$ with $\dim U\ge\alpha N/2$,
\begin{equation}\label{eq:contraction}
\iint r(t,v)^j\,d\mu(t)\,d\mu(v)\le s^{-2j}T_j(X),\qquad T_j(X)=\iint Z_X(t,v)^j\,dt\,dv,
\end{equation}
where $Z_X(t,v)=\sum_{i=1}^N\one\{x_i\in S_t\cap S_v\}=\inner{z(t)}{z(v)}$ does not depend on $U$.
\end{lemma}

\begin{proof}
The range of $P$ lies in $W^\perp$, so \eqref{eq:Paz} gives $Pa(t)=Pz(t)$. Write $k(t,v)=\inner{Pz(t)}{Pz(v)}$, so that $r(t,v)=k(t,v)/(\sigma(t)\sigma(v))$. By the definition of $\mu$,
\[
\iint r^j\,d\mu\,d\mu=S_G^{-2}\iint_{G\times G}k(t,v)^j\sigma(t)^{2-j}\sigma(v)^{2-j}\,dt\,dv\le s^{-2j}\iint_{G\times G}k(t,v)^j\,dt\,dv,
\]
Here $S_G^{-2}\le s^{-4}$ and $\sigma^{2-j}\le s^{2-j}$ on $G$, and $k^j\ge0$ since $j$ is even. In the $j$-fold tensor product $(\R^N)^{\otimes j}$ consider $M_j=\int_Gz(t)^{\otimes j}\,dt$. Since $\inner{Pz(t)}{Pz(v)}^j=\inner{P^{\otimes j}z(t)^{\otimes j}}{P^{\otimes j}z(v)^{\otimes j}}$ and $P^{\otimes j}$ is an orthogonal projection,
\begin{align*}
\iint_{G\times G}k(t,v)^j\,dt\,dv=\big\|P^{\otimes j}M_j\big\|^2\le\|M_j\|^2&=\iint_{G\times G}\inner{z(t)}{z(v)}^j\,dt\,dv\\
&\le\iint\inner{z(t)}{z(v)}^j\,dt\,dv=T_j(X).
\end{align*}
The last step enlarges the domain of integration, which is permitted because $z(t)$ has nonnegative coordinates. This proves \eqref{eq:contraction}.
\end{proof}

Choose $\varkappa=\varkappa_{d,\alpha}\in(0,1]$ sufficiently small
for Lemma~\ref{lem:expected}.
Let $p$ be the largest even integer with $p\le\varkappa h$. For $N$ so large that $\varkappa h\ge4$,
\begin{equation}\label{eq:pdef}
2\le\tfrac12\varkappa h\le p\le\varkappa h\le\log_2N.
\end{equation}
For every sample $X$, whether or not it lies in $\cS_N$, define
\begin{equation}\label{eq:Phi}
\Phi(X)=\sum_{j=2}^p\frac{x^j}{j!}T_j(X),\qquad x=\frac p{s^2}.
\end{equation}
By Lemma~\ref{lem:contraction}, for $X\in\cS_N$ and every admissible $U$, the left-hand side of \eqref{eq:varhyp} for the measure \eqref{eq:mu} satisfies
\begin{equation}\label{eq:VarPhi}
\sum_{\substack{2\le j\le p\\ j\text{ even}}}\frac{p^j}{j!}\iint r^j\,d\mu\,d\mu\le\Phi(X),
\end{equation}
where the odd terms in \eqref{eq:Phi} are nonnegative. The same statistic $\Phi(X)$ bounds the correlations for every subspace $U$ under consideration.

\begin{lemma}\label{lem:expected}
For sufficiently small $\varkappa=\varkappa_{d,\alpha}>0$,
there is a constant $C'_{d,\alpha}$ such that
\begin{equation}\label{eq:EPhi}
\E_X\Phi(X)\le C'_{d,\alpha}\frac{p^2}N
\end{equation}
for all sufficiently large $N$.
\end{lemma}

\begin{proof}
For $Z\sim\operatorname{Bin}(N,\pi)$,
\[
\E e^{xZ}=(1+\pi(e^x-1))^N
=\sum_{k=0}^N\binom Nk\pi^k(e^x-1)^k.
\]
The term indexed by $k$ starts in degree $k$. Truncating both sides
to degrees $2,\dots,p$ therefore removes every term with $k>p$.
Since all Taylor coefficients are nonnegative, for $x\ge0$ we get
\[
\sum_{j=2}^p\frac{x^j}{j!}\E Z^j
\le N\pi(e^x-1-x)
  +\sum_{k=2}^p\binom Nk\pi^k(e^x-1)^k.
\]
For fixed $t,v$, the variable $Z_X(t,v)$ is binomial with parameters
$N$ and $\vol(S_t\cap S_v)$. Integration over $t,v$ yields
\begin{equation}\label{eq:binommom}
\E_X\Phi(X)
\le NR_1(e^x-1-x)
  +\sum_{k=2}^p\binom NkR_k(e^x-1)^k.
\end{equation}
Since $s\asymp_{d,\alpha}L$, $p\le h$ and $q\ge1$,
$x=p/s^2\to0$. For sufficiently large $N$,
$e^x-1\le2x$ and $e^x-1-x\le x^2$.
By \eqref{eq:R1} and \eqref{eq:Neps},
\begin{equation}\label{eq:k1}
NR_1x^2
\le C N\eps^2L^2\frac{p^2}{L^4}
\le C\frac{p^2}N.
\end{equation}
Using \eqref{eq:Rj_bound} for $2\le k\le p\le h$ and
\eqref{eq:Neps}, we also obtain
\begin{equation}\label{eq:k2}
\sum_{k=2}^p\binom NkR_k(2x)^k
\le\frac{L^2}N
\sum_{k=2}^p\left(\frac{Cp}L\right)^k\frac{k^{kq}}{k!},
\end{equation}
where $C$ depends only on $d,\alpha$ and is uniform for
$0<\varkappa\le1$.
The ratio of consecutive summands on the right is at most
\[
\frac{Cp}L(k+1)^{q-1}\left(1+\frac1k\right)^{kq}
\le Ce^q\frac{p^q}L
\le Ce^q\varkappa^q.
\]
Choose $\varkappa$ so that $Ce^q\varkappa^q\le1/2$.
The sum in \eqref{eq:k2} is then at most twice its first term,
which is $O_{d,\alpha}(p^2/L^2)$.
Together with \eqref{eq:k1} and \eqref{eq:binommom}, this proves
\eqref{eq:EPhi}.
\end{proof}

\begin{proof}[Proof of Proposition~\ref{prop:gauss}]
Regard the sample as a point of the compact cube $[0,1]^{dN}$.
For fixed $N$, the dyadic definitions make
$(X,t,v)\mapsto Z_X(t,v)$ Borel, so $T_j$ and $\Phi$ are Borel
functions of $X$. The matrix $Q$ of Lemma~\ref{lem:sampled} is
continuous in $X$, so $\cS_N$ is also Borel.
Take $N$ sufficiently large for Lemmas~\ref{lem:sampled},
\ref{lem:overlap}, \ref{lem:expected} and \eqref{eq:pdef},
and let $C'=C'_{d,\alpha}$ be the constant in \eqref{eq:EPhi}.
Markov's inequality and Lemma~\ref{lem:sampled} give
\[
\Prob_X\Big\{X\in\cS_N,\ \Phi(X)\le4C'\frac{p^2}N\Big\}
\ge\frac12.
\]
Finite Borel measures on metric spaces are inner regular with
respect to closed sets \cite[Theorem 1.1]{Bil99}.
Since the cube is compact, this event contains a compact set $H$
with $\Prob_X\{H\}\ge1/4$. Choose $H$ independently of $B$.

Fix $X\in H$ and $U\subseteq\R^N$ with $\dim U\ge\alpha N/2$.
Let $P$ be the orthogonal projection onto $U\cap W^\perp$.
Since $PP_U=P$, Lemma~\ref{lem:criterion} applies to
$t\mapsto P_Ua_X(t)$, the projection $P$ and the measure $\mu$
from \eqref{eq:mu}, with $C=4C'$.
Indeed, $\sigma\ge s$ on $G$ by \eqref{eq:mu},
$p\le\log_2N$ by \eqref{eq:pdef}, and \eqref{eq:varhyp} follows
from \eqref{eq:VarPhi} and the choice of $H$.
For every $B>0$, we obtain
\[
\Prob_g\{\|A_XP_Ug\|_\infty\le c(C,B)s\sqrt p\}\le e^{-BN}.
\]
By \eqref{eq:sdef}, \eqref{eq:pdef} and \eqref{eq:hlog},
\[
s\sqrt p\asymp_{d,\alpha}h^{q+1/2}
\asymp_{d,\alpha}(\log N)^{d-1/2}.
\]
Thus there is $a_{d,\alpha,B}>0$ such that, with
$\varrho_N=a_{d,\alpha,B}(\log N)^{d-1/2}$,
\begin{equation}\label{eq:halfevent}
\Prob_g\{\|A_XP_Ug\|_\infty\le\varrho_N\}\le e^{-BN}
\qquad(X\in H,\ \dim U\ge\alpha N/2).
\end{equation}

For a sample $X$ let
\[
J(X)=\min\big\{|\{i\in[N]:x_i\ne y_i\}|:\ Y=(y_1,\dots,y_N)\in H\big\}
\]
be the least number of points of $X$ that must be changed to reach $H$. For $j\ge0$ the sublevel set $\{J\le j\}$ is the union, over the index sets $I\subseteq[N]$ with $|I|\ge N-j$, of the sets $\pi_I^{-1}(\pi_I(H))$, where $\pi_I$ projects onto the coordinates in $I$. Compactness of $\pi_I(H)$ makes these sets closed, so $J$ is Borel. Changing one point of $X$ changes $J(X)$ by at most one. By the bounded differences inequality of McDiarmid \cite[Lemma 1.2]{McD89}, applied to the independent points $x_1,\dots,x_N$, for every $u\ge0$
\[
\Prob_X\{J\le\E J-u\}\le e^{-2u^2/N},\qquad\Prob_X\{J\ge\E J+u\}\le e^{-2u^2/N}.
\]
Taking $u=\E J$ in the first inequality and using $\Prob_X\{J=0\}\ge\Prob_X\{H\}\ge1/4$ gives $\E J\le\sqrt{N\log(4)/2}$, which is at most $\alpha N/4$ for $N\ge8\log(4)/\alpha^2$. The second inequality with $u=\alpha N/4$ then gives
\begin{equation}\label{eq:Jtail}
\Prob_X\{J(X)>\alpha N/2\}\le e^{-\alpha^2N/8}.
\end{equation}

Let $X$ be a sample with $J(X)\le\alpha N/2$, and let $V\subseteq\R^N$ be a subspace with $\dim V\ge\alpha N$. There are $Y\in H$ and an index set $I$ with $|I|\ge(1-\alpha/2)N$ and $x_i=y_i$ for all $i\in I$. Let $\R^I$ be the subspace of vectors supported on $I$ and $U=V\cap\R^I$, so that $\dim U\ge\dim V+|I|-N\ge\alpha N/2$. For every $t$ and every $g$, the vector $P_Ug$ is supported on $I$, where the coordinates of $a_X(t)$ and $a_Y(t)$ agree, hence $\inner{a_X(t)}{P_Ug}=\inner{a_Y(t)}{P_Ug}$ and $\|A_XP_Ug\|_\infty=\|A_YP_Ug\|_\infty$. Since $Y\in H$, \eqref{eq:halfevent} gives
\[
\Prob_g\{\|A_XP_Ug\|_\infty\le\varrho_N\}=\Prob_g\{\|A_YP_Ug\|_\infty\le\varrho_N\}\le e^{-BN}.
\]
Since $U\subseteq V$, the vectors $A_XP_Ug$ and
$A_X(P_V-P_U)g$ are independent centered Gaussians.
Anderson's inequality \eqref{eq:anderson} gives
\[
\Prob_g\{\|A_XP_Vg\|_\infty\le\varrho_N\}\le\Prob_g\{\|A_XP_Ug\|_\infty\le\varrho_N\}\le e^{-BN}.
\]
Thus the assertion of the proposition holds on the Borel event $\{J(X)\le\alpha N/2\}$, whose probability is at least $1-e^{-b_{d,\alpha}N}$ with $b_{d,\alpha}=\alpha^2/8$ by \eqref{eq:Jtail}. The thresholds imposed on $N$ depend only on $d$ and $\alpha$.
\end{proof}

\subsection{Rounding random column restrictions}\label{sec:rounding}

For a real matrix $A$ with $N$ columns and $J\subseteq[N]$ write $A_J$ for the submatrix formed by the columns with indices in $J$, and
\[
\disc(A_J)=\min_{\chi\in\{-1,1\}^J}\|A_J\chi\|_\infty,
\]
with value zero for $J=\emptyset$, as in \cite{LSV86}. When $A=A_X$ and $J$ indexes a subset $X_J$ of the sample, the rows of $A_J$ are the restrictions of the vectors $a(t)$ to the coordinates in $J$, which give the anchored-box incidence vectors on $X_J$. Removing repeated rows gives
\begin{equation}\label{eq:restrdisc}
\disc((A_X)_J)=\disc(X_J,\cA_d).
\end{equation}
The following lemma uses the binary rounding of Lov\'asz, Spencer and Vesztergombi \cite{LSV86} on a random input. For a real matrix $A$ write $\|A\|_{\infty\to\infty}=\max_i\sum_j|A_{ij}|$ for its largest absolute row sum, so that $\|Au\|_\infty\le\|A\|_{\infty\to\infty}\|u\|_\infty$ for every vector $u$.

\begin{lemma}\label{lem:restriction}
Let $A$ be a real matrix with $N$ columns, and suppose that every column restriction has discrepancy at most $H$. Let $J$ be a uniformly random subset of $[N]$, each index belonging to $J$ independently with probability $1/2$, and for $D\ge0$ set $q_D=\Prob_J\{\disc(A_J)\le D\}$. Then for every integer $l\ge1$ and $g\sim N(0,I_N)$,
\begin{equation}\label{eq:restriction}
\Prob_g\big\{\|Ag\|_\infty\le(1-2^{-l})D+2^{-l}H\big\}\ge(2\pi e)^{-N/2}q_D^{\,l}.
\end{equation}
\end{lemma}

\begin{proof}
For every $S\subseteq[N]$ fix a vector $\chi_S\in\{-1,0,1\}^N$ with $\chi_{S,i}\ne0$ exactly for $i\in S$, attaining $\|A\chi_S\|_\infty=\disc(A_S)$, with $\chi_\emptyset=0$. Thus $\|A\chi_S\|_\infty\le H$ always, and $\|A\chi_S\|_\infty\le D$ whenever $\disc(A_S)\le D$. These choices are deterministic once $A$ is fixed.

Fix an integer $K\ge l$ and let $a^{(K)}$ be uniformly distributed on $\{0,1,\dots,2^K-1\}^N$. For $k=K,K-1,\dots,1$ define
\begin{equation}\label{eq:recursion}
J_k=\{i:\ a^{(k)}_i\text{ is odd}\},\qquad a^{(k-1)}=\frac{a^{(k)}+\chi_{J_k}}2.
\end{equation}
Induction gives $a^{(k)}\in\{0,\dots,2^k\}^N$. An even coordinate is halved, while an odd coordinate lies between $1$ and $2^k-1$ and changes by one before halving. Thus $z=a^{(0)}\in\{0,1\}^N$. Unwinding the recursion $a^{(k-1)}=2^{-1}a^{(k)}+2^{-1}\chi_{J_k}$ gives
\begin{equation}\label{eq:telescope}
z-2^{-K}a^{(K)}=\sum_{k=1}^K2^{-k}\chi_{J_k}.
\end{equation}

\smallskip
\noindent\emph{Independence of the parity sets.} We claim that, conditionally on $J_K,\dots,J_{k+1}$, the residue vector $a^{(k)}\bmod2^k$ is uniformly distributed on $(\Z/2^k\Z)^N$. This holds for $k=K$ by construction. Suppose it holds for some $k\ge1$, and write the residue of $a^{(k)}$ uniquely as $2b+\varsigma$ with $b\in\{0,\dots,2^{k-1}-1\}^N$ and $\varsigma\in\{0,1\}^N$. Conditionally on the earlier parity sets, $b$ and $\varsigma$ are independent and uniform. The vector $\varsigma$ determines $J_k=\{i:\varsigma_i=1\}$. Writing $a^{(k)}=2^kn+2b+\varsigma$ with an integer vector $n$, we get
\[
a^{(k-1)}=2^{k-1}n+b+\frac{\varsigma+\chi_{J_k}}2,
\]
and $(\varsigma+\chi_{J_k})/2$ is a vector in $\{0,1\}^N$ determined by $\varsigma$, since $\varsigma_i+\chi_{J_k,i}\in\{0,2\}$ for $i\in J_k$ and $=0$ otherwise. Hence $a^{(k-1)}\bmod2^{k-1}$ is a translate of the uniform vector $b$ by a vector depending only on $\varsigma$, so conditionally on $J_K,\dots,J_{k+1}$ and $J_k$ it is uniform on $(\Z/2^{k-1}\Z)^N$, and $J_k$ is independent of $J_K,\dots,J_{k+1}$ with the law of a uniformly random subset of $[N]$. By induction, $J_K,\dots,J_1$ are independent, each with this law.

By independence, with probability $q_D^{\,l}$ we have
$\disc(A_{J_k})\le D$ for every $1\le k\le l$.
On this event, \eqref{eq:telescope} gives
\[
\big\|A\big(z-2^{-K}a^{(K)}\big)\big\|_\infty
\le\sum_{k=1}^l2^{-k}D+\sum_{k=l+1}^K2^{-k}H
\le(1-2^{-l})D+2^{-l}H=:\varrho.
\]
Let $E_K$ be the union of the cubes
$x_0+[0,2^{-K})^N$, where $x_0=2^{-K}a^{(K)}$ ranges over
the grid points for which $\disc(A_{J_k})\le D$ for
$1\le k\le l$. These cubes are disjoint, lie in $[0,1)^N$,
and have total volume $q_D^{\,l}$.
For $x\in E_K$, with $z$ the rounded vector of the grid point of its cube, $x-z\in[-1,1]^N$ and
\[
\|A(x-z)\|_\infty\le\|A(x_0-z)\|_\infty+\|A(x-x_0)\|_\infty\le\varrho+2^{-K}\|A\|_{\infty\to\infty}.
\]
Set $C_K=\{u\in[-1,1]^N:\|Au\|_\infty\le\varrho+2^{-K}\|A\|_{\infty\to\infty}\}$, a closed set. Then $\sum_{z\in\Z^N}\one_{C_K}(x-z)\ge1$ for every $x\in E_K$. Integrating over $[0,1)^N$ and using Tonelli's theorem, the left-hand side becomes $\sum_{z\in\Z^N}\vol(C_K\cap([0,1)^N-z))=\vol(C_K)$, because the cubes $[0,1)^N-z$ tile $\R^N$. It follows that $\vol(C_K)\ge\vol(E_K)=q_D^{\,l}$. On $[-1,1]^N$ the standard Gaussian density is at least $(2\pi)^{-N/2}e^{-N/2}$, since $\|u\|_2^2\le N$ there. Therefore
\[
\Prob_g\{g\in C_K\}\ge(2\pi e)^{-N/2}q_D^{\,l}.
\]
The sets $C_K$ decrease as $K$ increases, and
$\bigcap_{K\ge l}C_K
=\{u\in[-1,1]^N:\|Au\|_\infty\le\varrho\}$.
By continuity from above,
\[
\Prob_g\{\|Ag\|_\infty\le\varrho\}
\ge\lim_{K\to\infty}\Prob_g\{g\in C_K\}
\ge(2\pi e)^{-N/2}q_D^{\,l}.
\]
This proves \eqref{eq:restriction}.
\end{proof}

Taking $D=H=\max_J\disc(A_J)$ recovers the Gaussian lower bound $(2\pi e)^{-N/2}$ for the body of the hereditary discrepancy, which goes back to \cite{LSV86}. For random point sets, the estimate will bound $q_D$ by requiring
discrepancy at most $D$ at a fixed number of rounding levels.

\subsection{Proof of Theorem \ref{thm:main}}

We need two elementary facts. First, for $0<\alpha\le1/2$,
\begin{equation}\label{eq:entropy}
\sum_{j\le\alpha N}\binom Nj\le e^{Nh_2(\alpha)},\qquad h_2(\alpha)=-\alpha\log\alpha-(1-\alpha)\log(1-\alpha).
\end{equation}
Indeed, $1=\sum_j\binom Nj\alpha^j(1-\alpha)^{N-j}\ge\sum_{j\le\alpha N}\binom Nj\alpha^j(1-\alpha)^{N-j}$, and since $\alpha/(1-\alpha)\le1$ the weight $\alpha^j(1-\alpha)^{N-j}$ is nonincreasing in $j$, so it is at least $\alpha^{\alpha N}(1-\alpha)^{(1-\alpha)N}=e^{-Nh_2(\alpha)}$ for $j\le\alpha N$. Second, if $J$ is a uniformly random subset of $[N]$ and $S\subseteq[N]$ is fixed, then the symmetric difference $J\triangle S$ is again a uniformly random subset of $[N]$, because $J\mapsto J\triangle S$ is a bijection of the subsets of $[N]$.

\begin{proof}[Proof of Theorem~\ref{thm:main}]
Put $N=2n$ and apply Proposition~\ref{prop:gauss} with $\alpha=1/2$ and $B=1+\frac12\log(2\pi e)$. Let $\cG_N$ be its sample event, of probability at least $1-e^{-b_dN}$ with $b_d=b_{d,1/2}$. For a parent sample $X\in\cG_N$ of size $N$ and $V=\R^N$ we have, for all large $N$,
\[
\Prob_g\{\|A_Xg\|_\infty\le r\}\le e^{-BN},\qquad r=a_{d,1/2,B}(\log N)^{d-1/2}.
\]
Put $D=r/2$. By Nikolov's theorem \eqref{eq:nikolov} and \eqref{eq:restrdisc}, every column restriction of $A_X$ has discrepancy at most $H=C_d(\log N)^{d-1/2}$, where we enlarge $C_d$ so that this includes the restrictions to at most one point, which have discrepancy at most one. Let $l=l_d$ be the least positive integer with $2^{-l}C_d\le a_{d,1/2,B}/2$, so that $2^{-l}H\le r/2$. This choice of $l$ depends only on $d$. Lemma~\ref{lem:restriction} gives, for $X\in\cG_N$,
\[
(2\pi e)^{-N/2}q_D^{\,l}\le\Prob_g\big\{\|A_Xg\|_\infty\le(1-2^{-l})D+2^{-l}H\big\}\le\Prob_g\{\|A_Xg\|_\infty\le r\}\le e^{-BN},
\]
hence $q_D^{\,l}\le e^{-N}$ and
\begin{equation}\label{eq:qD}
q_D\le e^{-N/l}.
\end{equation}

Choose $\alpha_d\in(0,1/2]$ so that
$h_2(\alpha_d)\le1/(2l)$, and let $J$ be a uniformly random subset
of $[N]$. For each fixed $S\subseteq[N]$, the set $J\triangle S$
is again uniformly random. A union bound and
\eqref{eq:entropy}, \eqref{eq:qD} give
\[
\Prob_J\Big\{\exists S\subseteq[N],\ |S|\le\alpha_dN,
  \disc((A_X)_{J\triangle S})\le D\Big\}
\le q_D\sum_{j\le\alpha_dN}\binom Nj
\le e^{-N/(2l)}.
\]
Since $N=2n$, the largest atom of $|J|\sim\operatorname{Bin}(N,1/2)$
is at $n$, so $\Prob\{|J|=n\}\ge1/(N+1)$.
Conditionally on $|J|=n$, the set $J$ is a uniform $n$-element
subset of $[N]$, and the preceding bound implies
\[
\Prob_J\Big\{\exists S\subseteq[N],\ |S|\le\alpha_dN,
  \disc((A_X)_{J\triangle S})\le D\ \Big|\ |J|=n\Big\}
\le(N+1)e^{-N/(2l)}.
\]
Deleting at most $\alpha_dn$ indices from $J$ produces
$J\triangle S$ for some $S\subseteq J$ with
$|S|\le\alpha_dn\le\alpha_dN$.
Thus, for every $X\in\cG_N$ and a uniform $n$-element subset $J$,
with probability at least $1-(N+1)e^{-N/(2l)}$, every such deletion
leaves anchored discrepancy greater than $D$, by
\eqref{eq:restrdisc}.

Now let the parent sample $X_N$ be random and let $J$ be a uniform $n$-element subset of $[N]$, independent of $X_N$. The points $x_i$, $i\in J$, are $n$ independent uniform points of $(0,1)^d$, so the subsample $X_J$ has the law of $X_n$. The event that some subset of $X_J$ obtained by deleting at most $\alpha_dn$ points has anchored discrepancy at most $D$ has probability at most
\[
\Prob_X\{X_N\notin\cG_N\}+(N+1)e^{-N/(2l)}\le e^{-b_dN}+(N+1)e^{-N/(2l)}\le e^{-b_d'n}
\]
for $n\ge n_d$. Since $D=\frac12a_{d,1/2,B}(\log2n)^{d-1/2}\ge c_d(\log n)^{d-1/2}$, it follows that for $n\ge n_d$, with probability at least $1-e^{-b'_dn}$, every $Y\subseteq X_n$ with $|X_n\setminus Y|\le\alpha_dn$ satisfies $\disc(Y,\cA_d)\ge c_d(\log n)^{d-1/2}$. Renaming $n$ as $N$ proves the second assertion of the theorem with $\eta_d=\alpha_d$.

Taking $Y=X_N$ and using \eqref{eq:anchored} gives $\Delta_d(N)\ge c_d(\log N)^{d-1/2}$ for $N\ge N_d$. Since $\Delta_d(N)\ge1$ for $N\ge2$, decreasing $c_d$ covers the remaining values of $N$. The upper bound is Nikolov's theorem \eqref{eq:nikolov}.
\end{proof}

\section{Lower bounds for every fixed deletion fraction}\label{sec:fixed}

The source fixes the sign of each retained observation. We first bound the number of ternary vectors with small discrepancy, then fold the observations to obtain independent fair signs.

\subsection{Counting ternary vectors}\label{sec:counting}

For subsets $K,L$ of a Euclidean space, let $N(K,L)$ be the least number of translates of $L$ covering $K$. Covers compose by adding their translation centers, so
\begin{equation}\label{eq:covercompose}
N(K,L_1)\le N(K,L_2)N(L_2,L_1),\qquad
N(tK,tL)=N(K,L)\quad(t>0).
\end{equation}
For a subspace $V$, let $B_2^V$ be its closed Euclidean unit ball and $\gamma_V$ the standard Gaussian measure on $V$. Pisier's regular ellipsoid theorem with regularity parameter one gives an absolute constant $C_0\ge1$ with the following property \cite[Corollary 0.2]{Pis89}. For every finite-dimensional Euclidean space $W$ and every compact convex set $K=-K\subset W$ with nonempty interior, there is an ellipsoid $E\subset W$ centered at the origin such that
\begin{equation}\label{eq:pisier}
\log N(K,tE)\le\frac{C_0\dim W}{t},\qquad
\log N(E,tK)\le\frac{C_0\dim W}{t}
\qquad(t\ge C_0).
\end{equation}

\begin{lemma}\label{lem:count}
For every $\tau>0$ there are constants $\kappa,C>0$, depending only on $\tau$, such that for every integer $N\ge1$ and every compact convex set $K=-K\subseteq[-1,1]^N$ containing at least $e^{\tau N}$ points of $\{-1,0,1\}^N$, some subspace $V\subseteq\R^N$ satisfies
\[
\dim V\ge\kappa N,\qquad \gamma_V(K\cap V)\ge e^{-CN}.
\]
\end{lemma}

\begin{proof}
Let $W$ be the span of $K$. The set $K$ has nonempty interior relative to $W$, and $\dim W\ge1$ because $K$ contains at least two ternary points. Put $R=\sqrt N$, so that $K\subseteq RB_2^W$.

Choose $u\in(0,1/2)$ with $h_2(u)+u\log2\le\tau/2$. By \eqref{eq:entropy}, the number of ternary points in a Hamming ball of radius $uN$ is at most
\[
\sum_{j\le uN}\binom Nj2^j\le e^{N(h_2(u)+u\log2)}\le e^{\tau N/2}.
\] Selecting a ternary point of $K$ and discarding its Hamming ball at each step gives at least $e^{\tau N/2}$ points with pairwise Euclidean distances greater than $\sqrt u\,R$. Set $\delta=\sqrt u/4$. Every ball of radius $\delta R$ contains at most one selected point, hence
\begin{equation}\label{eq:packing}
\log N(K,\delta RB_2^W)\ge\frac{\tau N}{2}.
\end{equation}

Choose $t\ge C_0$, depending only on $\tau$, with $C_0/t\le\tau/8$. For the ellipsoid $E$ supplied by \eqref{eq:pisier}, the bound $\dim W\le N$ gives
\begin{equation}\label{eq:pisierN}
\log N(K,tE)\le\frac{\tau N}{8},\qquad
\log N(E,tK)\le\frac{\tau N}{8}.
\end{equation}
Put $r=\delta R/t$. The covering inequalities give
\[
N(K,\delta RB_2^W)
\le N(K,tE)N(tE,\delta RB_2^W)
=N(K,tE)N(E,rB_2^W).
\]
Together with \eqref{eq:packing} and \eqref{eq:pisierN}, this yields
\begin{equation}\label{eq:Elower}
\log N(E,rB_2^W)\ge\frac{3\tau N}{8}.
\end{equation}
Also, $tK\subseteq tRB_2^W$ implies
\begin{equation}\label{eq:Eupper}
\log N(E,tRB_2^W)\le\log N(E,tK)\le\frac{\tau N}{8}.
\end{equation}

\smallskip
\noindent\emph{Counting long semiaxes.}
Let $a_1\ge\cdots\ge a_{\dim W}>0$ be the semiaxes of $E$, with orthonormal directions $e_i$. Let $j$ count the semiaxes greater than $r/2$ and put $V=\operatorname{span}\{e_1,\dots,e_j\}$. The bound \eqref{eq:Elower} ensures $j\ge1$, since otherwise $E\subseteq(r/2)B_2^W$.

The projection of $E$ onto $V$ equals $E\cap V$. Take a maximal family of points in $E\cap V$ at pairwise distances greater than $r/2$. The balls of radius $r/4$ centered at these points are disjoint. Since every semiaxis of $E\cap V$ exceeds $r/2$, these balls lie in
\[
(E\cap V)+\tfrac r4B_2^V\subseteq\tfrac32(E\cap V).
\]
Comparing volumes in $V$, the family has at most $\prod_{i\le j}(6a_i/r)$ points. By maximality it covers $E\cap V$ within distance $r/2$. For $x=\sum_i x_i e_i\in E$, the orthogonal component satisfies
\[
\sum_{i>j}x_i^2\le\frac{r^2}{4}\sum_{i>j}\frac{x_i^2}{a_i^2}\le\frac{r^2}{4}.
\]
Thus the same family covers $E$ within distance $r/\sqrt2\le r$, and
\begin{equation}\label{eq:semiaxes1}
\frac{3\tau N}{8}\le\log N(E,rB_2^W)
\le\sum_{i\le j}\log\frac{6a_i}{r}.
\end{equation}

Next project a cover of $E$ by balls of radius $tR$ onto the span of the axes with $a_i>tR$. The projected balls cover the ellipsoid with those semiaxes. Comparing volumes and using \eqref{eq:Eupper} gives
\begin{equation}\label{eq:semiaxes2}
\sum_{i\le\dim W}\log^+\frac{a_i}{tR}
\le\log N(E,tRB_2^W)\le\frac{\tau N}{8}.
\end{equation}
Since $tR/r=t^2/\delta$, the two estimates imply
\[
\frac{3\tau N}{8}
\le j\log\frac{6t^2}{\delta}
   +\sum_{i\le j}\log^+\frac{a_i}{tR}
\le j\log\frac{6t^2}{\delta}+\frac{\tau N}{8}.
\]
Hence $j\ge\kappa N$ for $\kappa=\tau/(4\log(6t^2/\delta))$.

\smallskip
\noindent\emph{Gaussian measure of the section.}
The section $E\cap V$ contains $(r/2)B_2^V$, hence the cube with coordinate bounds $\delta/(2t)$ in the basis $e_1,\dots,e_j$, because
$\sqrt j\,\delta/(2t)\le\sqrt N\,\delta/(2t)=r/2$.
Put
\[
p=\gamma_1\bigl([-\delta/(2t),\delta/(2t)]\bigr)>0.
\]
Since $j\le N$,
\[
\gamma_V(E\cap V)\ge p^j\ge p^N.
\]
Anderson's inequality \eqref{eq:anderson}, applied to a standard Gaussian supported on $V$ and the symmetric convex set $tK\subset W$, gives
\[
\gamma_V((y+tK)\cap V)\le\gamma_V(tK\cap V)
\qquad(y\in W).
\]
Covering $E$ by at most $e^{\tau N/8}$ translates of $tK$ therefore yields
\[
p^N\le\gamma_V(E\cap V)\le e^{\tau N/8}\gamma_V(tK\cap V).
\]
A change of variables, using $t\ge1$ and $j\le N$, gives
\[
\gamma_V(K\cap V)
=t^{-j}\int_{tK\cap V}
  \frac{e^{-\|x\|_2^2/(2t^2)}}{(2\pi)^{j/2}}\,dx
\ge t^{-N}\gamma_V(tK\cap V)
\ge(p/t)^N e^{-\tau N/8}.
\]
This proves the claim with $C=\log(t/p)+\tau/8$. Both $\kappa$ and $C$ depend only on $\tau$.
\end{proof}

\begin{proposition}\label{prop:count}
For each fixed $d\ge2$ and $\tau>0$ there are positive constants $a,b$, depending on $d$ and $\tau$, such that for all sufficiently large $N$, with probability at least $1-e^{-bN}$ over $N$ independent uniform points $X$ in $(0,1)^d$,
\[
\#\big\{z\in\{-1,0,1\}^N:\ \|A_Xz\|_\infty\le a(\log N)^{d-1/2}\big\}<e^{\tau N}.
\]
\end{proposition}

\begin{proof}
Take $\kappa\le1/2$ and $C$ from Lemma~\ref{lem:count}, and apply Proposition~\ref{prop:gauss} with $\alpha=\kappa$ and $B=C+1$. Let $a,b$ be the resulting constants. Suppose that $X$ lies in its sample event and that the asserted count fails. Then
\[
K=[-1,1]^N\cap\{x\in\R^N:\|A_Xx\|_\infty\le a(\log N)^{d-1/2}\}
\]
is compact, convex and symmetric and contains at least $e^{\tau N}$ ternary points. Lemma~\ref{lem:count} gives a subspace $V$ with $\dim V\ge\kappa N$ and $\gamma_V(K\cap V)\ge e^{-CN}$. Since $P_Vg$ is a standard Gaussian in $V$,
\[
\gamma_V(K\cap V)
=\Prob_g\{P_Vg\in K\}
\le\Prob_g\big\{\|A_XP_Vg\|_\infty\le a(\log N)^{d-1/2}\big\}
\le e^{-(C+1)N}<e^{-CN},
\]
a contradiction.
\end{proof}

\begin{corollary}\label{cor:compatible}
Fix $d\ge2$. Let $x_1,\dots,x_N$ be independent uniform points in $(0,1)^d$ and let $\xi_1,\dots,\xi_N$ be independent fair signs, independent of the points. For every $0<\beta\le1$ there are constants $c_{d,\beta},b_{d,\beta}>0$ such that, for all sufficiently large $N$,
\begin{equation}\label{eq:compatible}
\Prob\Big\{\exists z\in\prod_{i=1}^N\{0,\xi_i\}:\ 
|\{i:z_i\ne0\}|\ge\beta N,\ 
\|A_Xz\|_\infty\le c_{d,\beta}(\log N)^{d-1/2}\Big\}
\le e^{-b_{d,\beta}N}.
\end{equation}
\end{corollary}

\begin{proof}
Apply Proposition~\ref{prop:count} with $\tau=\frac12\beta\log2$, and let $a,b>0$ be its constants. Put $c_{d,\beta}=a$ and condition on the points. Outside a sample event of probability at most $e^{-bN}$, fewer than $e^{\tau N}$ ternary vectors satisfy $\|A_Xz\|_\infty\le a(\log N)^{d-1/2}$. A fixed vector with $s$ nonzero entries lies in $\prod_i\{0,\xi_i\}$ with conditional probability $2^{-s}$, since each sign on its support is prescribed. A union bound over the vectors with at least $\beta N$ nonzero entries gives
\[
e^{\tau N}2^{-\beta N}=e^{-(\beta\log2)N/2}.
\]
Adding the exceptional sample probability proves the claim.
\end{proof}

\subsection{Folding the first coordinate}\label{sec:folding}

Use the completed reveal order of Section~\ref{sec:model}, put $N=2n$, and let $a_i\in\{0,1\}$ indicate whether $x_i$ is kept. With the source signs $\eps_i$ of that model,
\[
S_j(t)=\sum_{i\le j}a_i\eps_i\one\{x_i\le t\},
\qquad t\in[0,1]^k,\quad 0\le j\le N.
\]

\begin{lemma}\label{lem:folding}
Consider the uniform law on $(0,1)^k$. There are independent identically distributed pairs $(y_i,\xi_i)$, $1\le i\le N$, such that $y_i$ is uniform on $(0,1)^k$, $\xi_i$ is a fair sign independent of $y_i$, and, simultaneously for $0\le j\le N$,
\begin{equation}\label{eq:fold}
\sup_{t\in[0,1]^k}\Big|\sum_{i\le j}a_i\xi_i\one\{y_i\le t\}\Big|
\le3\sup_{t\in[0,1]^k}|S_j(t)|.
\end{equation}
\end{lemma}

\begin{proof}
Let $\cF_{i-1}$ be generated by the seed and the sources and locations of the first $i-1$ observations. The source sign $\eps_i$ is $\cF_{i-1}$-measurable. On the event that source $s$ is selected and $m$ of its observations have been read, every event of $\cF_{i-1}$ restricted to this event is measurable with respect to the seed, the first $m$ observations of that source and the first $i-1-m$ observations of the other source. The next observation of source $s$ is independent of these variables and uniform. Summing over $s$ and $m$ shows that the conditional law of $x_i$ given $\cF_{i-1}$ is uniform on $(0,1)^k$.

Define
\[
b_i=\one\{x_{i1}>1/2\},\qquad
y_i=(2x_{i1}-b_i,x_{i2},\dots,x_{ik}),\qquad
\xi_i=\eps_i(-1)^{b_i}.
\]
Conditionally on $\cF_{i-1}$, the bit $b_i$ is fair and independent of the uniform point $y_i$. Since $\eps_i$ is $\cF_{i-1}$-measurable, the conditional law of $(y_i,\xi_i)$ is the product of the uniform law and a fair sign. The earlier pairs are $\cF_{i-1}$-measurable, so induction gives the asserted independence.

For $t=(t_1,t')$, with $t'$ omitted when $k=1$, splitting the first coordinate at $1/2$ gives
\[
\sum_{i\le j}a_i\xi_i\one\{y_i\le t\}
=S_j(t_1/2,t')-S_j((1+t_1)/2,t')+S_j(1/2,t').
\]
Indeed, when $b_i=0$ the last two indicators cancel, leaving $x_{i1}\le t_1/2$. When $b_i=1$, only the negative indicator remains. Taking the supremum and applying the triangle inequality proves \eqref{eq:fold}.
\end{proof}
\begin{figure}[t]
\centering
\begin{tikzpicture}[x=5.4cm,y=5.4cm]
\begin{scope}
\fill[black!13] (0,0) rectangle (0.31,0.68);
\fill[black!13] (0.5,0) rectangle (0.81,0.68);
\draw[thick] (0,0) rectangle (0.31,0.68);
\draw[thick] (0.5,0) rectangle (0.81,0.68);
\draw (0,0) rectangle (1,1);
\draw[densely dashed] (0.5,0) -- (0.5,1);
\fill (0.6419,0.1066) circle (1.5pt);
\fill (0.3737,0.0934) circle (1.5pt);
\fill (0.0652,0.439) circle (1.5pt);
\fill (0.2562,0.547) circle (1.5pt);
\fill (0.2398,0.6172) circle (1.5pt);
\fill (0.5725,0.4049) circle (1.5pt);
\draw[fill=white] (0.0738,0.8298) circle (1.5pt);
\draw[fill=white] (0.424,0.5374) circle (1.5pt);
\fill (0.1269,0.5655) circle (1.5pt);
\fill (0.3801,0.5439) circle (1.5pt);
\draw[fill=white] (0.5605,0.6095) circle (1.5pt);
\draw[fill=white] (0.6696,0.4334) circle (1.5pt);
\draw[fill=white] (0.4677,0.8896) circle (1.5pt);
\draw[fill=white] (0.6871,0.2646) circle (1.5pt);
\fill (0.7157,0.3049) circle (1.5pt);
\draw[fill=white] (0.141,0.4247) circle (1.5pt);
\fill (0.1729,0.4898) circle (1.5pt);
\draw[fill=white] (0.9343,0.1114) circle (1.5pt);
\draw[fill=white] (0.3497,0.3622) circle (1.5pt);
\fill (0.5751,0.4597) circle (1.5pt);
\fill (0.918,0.4762) circle (1.5pt);
\draw[fill=white] (0.6383,0.9536) circle (1.5pt);
\draw[fill=white] (0.6585,0.0608) circle (1.5pt);
\draw[fill=white] (0.3641,0.602) circle (1.5pt);
\fill (0.0854,0.7468) circle (1.5pt);
\draw[fill=white] (0.7241,0.4061) circle (1.5pt);
\draw[fill=white] (0.1057,0.4533) circle (1.5pt);
\draw[fill=white] (0.8604,0.7937) circle (1.5pt);
\draw[fill=white] (0.694,0.9475) circle (1.5pt);
\fill (0.9303,0.1788) circle (1.5pt);
\node[font=\small, fill=white, inner sep=1.2pt, circle] at (0.155,0.75) {$+$};
\node[font=\small, fill=white, inner sep=1.2pt, circle] at (0.655,0.75) {$-$};
\node[font=\scriptsize, above] at (0.25,1) {$b_i=0$};
\node[font=\scriptsize, above] at (0.75,1) {$b_i=1$};
\node[font=\scriptsize, below] at (0,0) {$0$};
\node[font=\scriptsize, below] at (0.31,0) {$\frac{t_1}2$};
\node[font=\scriptsize, below] at (0.5,0) {$\frac12$};
\node[font=\scriptsize, below] at (0.81,0) {$\frac{1+t_1}2$};
\node[font=\scriptsize, below] at (1,0) {$1$};
\node[font=\scriptsize, left] at (0,0.68) {$t_2$};
\node[font=\small] at (0.5,1.13) {(a)};
\end{scope}
\draw[->, thick] (1.08,0.5) -- (1.32,0.5);
\node[font=\scriptsize, above, align=center] at (1.2,0.52) {$y_{i1}=2x_{i1}-b_i$};
\begin{scope}[xshift=7.6cm]
\fill[black!13] (0,0) rectangle (0.62,0.68);
\draw[thick] (0,0) rectangle (0.62,0.68);
\draw (0,0) rectangle (1,1);
\draw[fill=white] (0.2838,0.1066) circle (1.5pt);
\fill (0.7475,0.0934) circle (1.5pt);
\fill (0.1305,0.439) circle (1.5pt);
\fill (0.5124,0.547) circle (1.5pt);
\fill (0.4797,0.6172) circle (1.5pt);
\draw[fill=white] (0.145,0.4049) circle (1.5pt);
\draw[fill=white] (0.1476,0.8298) circle (1.5pt);
\draw[fill=white] (0.848,0.5374) circle (1.5pt);
\fill (0.2537,0.5655) circle (1.5pt);
\fill (0.7601,0.5439) circle (1.5pt);
\fill (0.121,0.6095) circle (1.5pt);
\fill (0.3392,0.4334) circle (1.5pt);
\draw[fill=white] (0.9353,0.8896) circle (1.5pt);
\fill (0.3741,0.2646) circle (1.5pt);
\draw[fill=white] (0.4314,0.3049) circle (1.5pt);
\draw[fill=white] (0.282,0.4247) circle (1.5pt);
\fill (0.3457,0.4898) circle (1.5pt);
\fill (0.8686,0.1114) circle (1.5pt);
\draw[fill=white] (0.6994,0.3622) circle (1.5pt);
\draw[fill=white] (0.1502,0.4597) circle (1.5pt);
\draw[fill=white] (0.836,0.4762) circle (1.5pt);
\fill (0.2766,0.9536) circle (1.5pt);
\fill (0.3171,0.0608) circle (1.5pt);
\draw[fill=white] (0.7283,0.602) circle (1.5pt);
\fill (0.1708,0.7468) circle (1.5pt);
\fill (0.4481,0.4061) circle (1.5pt);
\draw[fill=white] (0.2115,0.4533) circle (1.5pt);
\fill (0.7208,0.7937) circle (1.5pt);
\fill (0.388,0.9475) circle (1.5pt);
\draw[fill=white] (0.8605,0.1788) circle (1.5pt);
\node[font=\scriptsize, below] at (0,0) {$0$};
\node[font=\scriptsize, below] at (0.62,0) {$t_1$};
\node[font=\scriptsize, below] at (1,0) {$1$};
\node[font=\scriptsize, left] at (0,0.68) {$t_2$};
\node[font=\small] at (0.5,1.13) {(b)};
\end{scope}
\end{tikzpicture}
\caption{Folding the first coordinate, for $k=2$. (a) Observations $x_i$ with their source signs $\varepsilon_i$, filled for $+1$ and open for $-1$; the bit $b_i$ records the half of the square containing $x_i$. (b) The points $y_i$, obtained by stretching each half to the whole square, with the signs $\xi_i=\varepsilon_i(-1)^{b_i}$, so that the points of the right half change their marks. The box $\{y\le t\}$ in (b) corresponds to the two shaded rectangles in (a), counted with the signs shown. The left rectangle is counted by $S_j(t_1/2,t_2)$ and the right one by $S_j((1+t_1)/2,t_2)-S_j(1/2,t_2)$, which gives the factor $3$ in \eqref{eq:fold}.}
\label{fig:fold}
\end{figure}
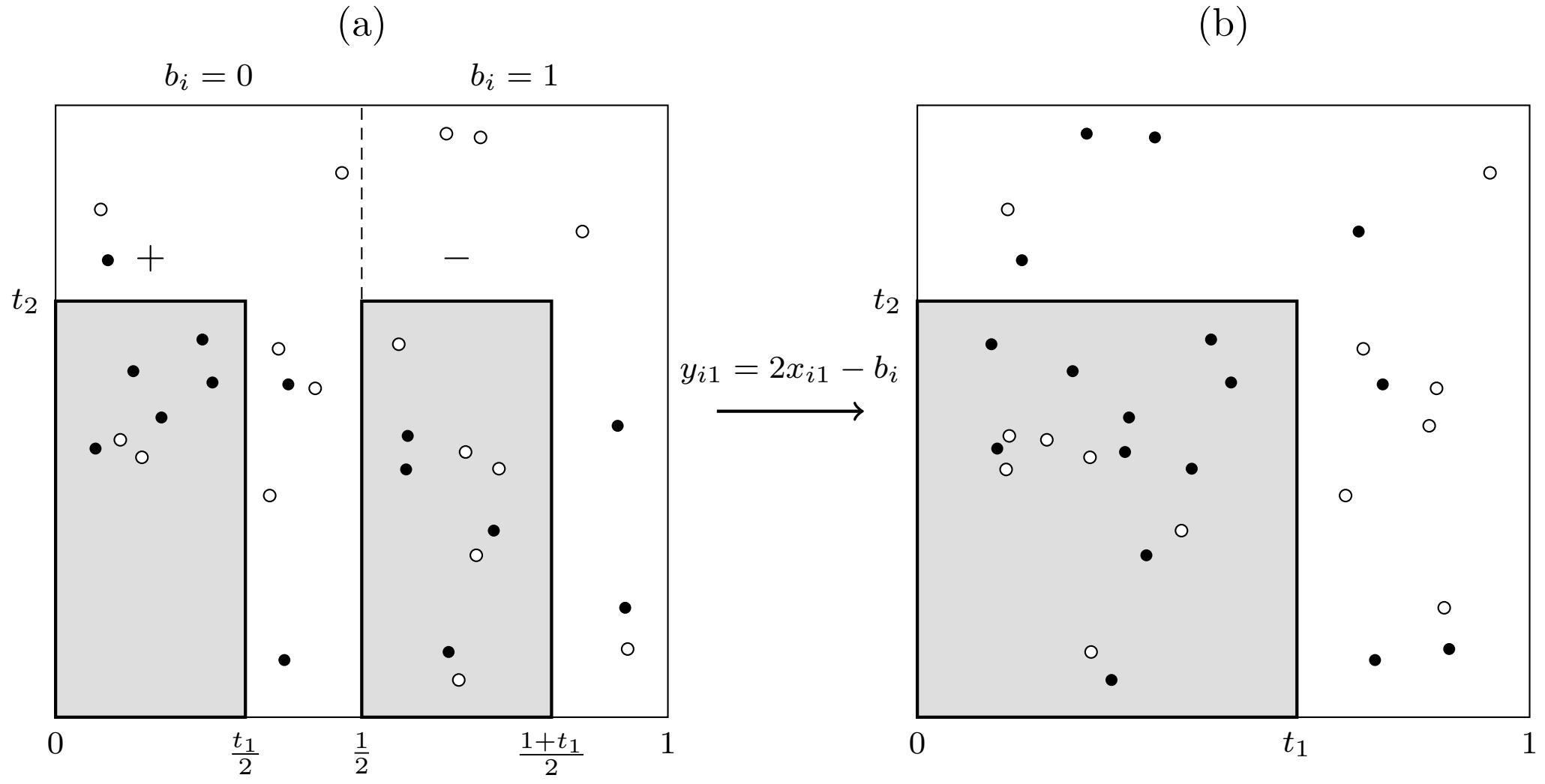
\begin{proof}[Proof of the lower bounds in Theorem~\ref{thm:multi}]
Let $\mu$ be uniform on $(0,1)^k$, put $N=2n$, and fix a retained fraction $\beta\in(0,1)$. On a finite sample in an open unit cube, strict and weak orthants cut out the same subsets by moving their boundaries between consecutive coordinate values. Thus Corollary~\ref{cor:compatible} applies to the weak inequalities in \eqref{eq:St}.

\smallskip
\noindent\emph{Terminal discrepancy.}
Let $k\ge2$. Apply Corollary~\ref{cor:compatible} in dimension $d=k$ to the folded points and signs from Lemma~\ref{lem:folding}. The output gives $z=(a_i\xi_i)_{i\le N}$, with $N-R$ nonzero entries. Thus, on $\{R\le(1-\beta)N\}$, the corollary and \eqref{eq:fold} give
\[
D\ge\tfrac13\sup_t\Big|\sum_i a_i\xi_i\one\{y_i\le t\}\Big|
\ge\tfrac13c_{k,\beta}(\log N)^{k-1/2}
\]
outside an event of probability at most $e^{-b_{k,\beta}N}$.

For an offline procedure, fix the reveal order in advance and choose the $a_i$ after seeing both samples. The folded pairs still have the law of Lemma~\ref{lem:folding}, and both Corollary~\ref{cor:compatible} and \eqref{eq:fold} apply to every choice of the $a_i$.

\smallskip
\noindent\emph{Prefix discrepancy.}
Let $k\ge1$. Independently of the run, take the order statistics
$u_{(1)}<\cdots<u_{(N)}$ of $N$ independent uniform variables on $(0,1)$, and put $\widehat y_i=(y_i,u_{(i)})$.
The marked sequence $(\widehat y_i,\xi_i)_{i\le N}$ has the law of independent uniform points in $(0,1)^{k+1}$ with independent fair marks, sorted by their last coordinate. The event in \eqref{eq:compatible} is invariant under a common permutation of the points and marks, so that corollary applies in dimension $k+1$.
An orthant with last threshold $\theta$ selects the arrival prefix $\{i:u_{(i)}\le\theta\}$, and every prefix occurs. Hence
\[
\sup_{(t,\theta)}\Big|\sum_i a_i\xi_i\one\{\widehat y_i\le(t,\theta)\}\Big|
=\max_{0\le j\le N}\sup_t\Big|\sum_{i\le j}a_i\xi_i\one\{y_i\le t\}\Big|.
\]
Corollary~\ref{cor:compatible} and \eqref{eq:fold} therefore give
\[
\Dpre\ge\tfrac13c_{k+1,\beta}(\log N)^{k+1/2}
\qquad\text{on }\{R\le(1-\beta)N\}
\]
outside an event of probability at most $e^{-b_{k+1,\beta}N}$.

\smallskip
\noindent\emph{Deletion budget.}
If $\E R\le\rho n$, Markov's inequality gives
\[
\Prob\{R>(1-\beta)N\}\le\frac{\rho n}{(1-\beta)N}
=\frac{\rho}{2(1-\beta)}.
\]
For $0<\rho\le1/8$, take $\beta=1/2$. Each lower bound then holds with probability at least $1-\rho-e^{-b_kn}$, with thresholds $c_k(\log n)^{k+1/2}$ for the prefix objective and $c_k(\log n)^{k-1/2}$ for the terminal objective when $k\ge2$. The constants $c_k,b_k$ and the large-sample threshold depend only on $k$.

For an arbitrary fixed $0<\rho<2$, take $\beta=(2-\rho)/4$. The probability is at least
\[
\frac{2-\rho}{2+\rho}-e^{-b_{k,\rho}n},
\]
which is bounded away from zero for large $n$. The uniform law is among the input laws in the minimax definitions. Taking expectations gives the asserted prefix lower bound for $k\ge1$ and terminal lower bound for $k\ge2$.
\end{proof}

\section{Uniform coordinates and thinning rules}\label{sec:tools}

We use the model of Section~\ref{sec:model}. All probabilities and expectations include the randomness of the samples and of the algorithm.

\subsection{Tags and uniform coordinates}\label{sec:tags}

We first give the reduction on the line. In dimension $k$, apply it coordinatewise using independent tags for each observation and coordinate. The argument uses only the marginal laws and therefore also applies to dependent coordinates.

Following the uniformization in Smirnov and Vershynin \cite[Section 3]{SV26a}, we use the distributional transform recorded by R\"uschendorf \cite[Section 2]{Rus09}. Write $F(b)=\mu((-\infty,b])$ and $F(b^-)=\mu((-\infty,b))$. Attach independent uniform tags $u_i\in(0,1)$ to the observations $x_i$, and put
\[
\widehat x_i=F(x_i^-)+u_i\big(F(x_i)-F(x_i^-)\big).
\]
Applied to independent observations from $\mu$, this transform produces independent uniform variables on $(0,1)$. Almost surely the transformed values are pairwise distinct and their order agrees with the lexicographic order of the tagged observations $(x_i,u_i)$, and every original half-line $(-\infty,b]$ cuts out an initial segment in this order.

The algorithms use tagged observations only through their order, so they make the same decisions in uniform coordinates. Since each original half-line is an initial segment of the tagged order, simultaneously for all $t$,
\begin{equation}\label{eq:tags-dominate}
\sup_{b\in\R}|S_t(b)|\le\sup_{x\in(0,1)}\Big|\sum_{\substack{i\le t\\ x_i\text{ kept}}}\eps_i\one\{\widehat x_i\le x\}\Big|.
\end{equation}
We analyze the upper-bound constructions in uniform coordinates. The reduction also applies to the offline comparison in Section~\ref{sec:model}.

\subsection{Source coins and controllers}\label{sec:general-model}

We use the independent-sign reduction of Smirnov and Vershynin \cite[Section 2]{SV26a}. Let $(V_i,\eps_i)_{i\ge1}$ be independent pairs, where $V_i\sim\mu$ and $\eps_i$ is an independent fair sign. The subsequences with signs $+1$ and $-1$ are independent samples from $\mu$, independent of the source coins. Take their first $n$ terms as the two samples.

For a construction that plans to reveal $2n-p$ points and leave $p$ unread, select each source by the next coin and stop if that source is exhausted. Let $E_p$ be the event that neither sign occurs more than $n$ times among the first $2n-p$ coins. On $E_p$, the sample run agrees with the first $2n-p$ points of the signed stream. More generally, the two runs agree until the sample run stops. For $1\le p<2n$, Hoeffding's inequality \cite{Hoe63} for the two coin counts gives
\begin{equation}\label{eq:source-failure}
\Prob(E_p^c)
\le2\exp\Big(-\frac{p^2}{2(2n-p)}\Big)
\le2\exp\Big(-\frac{p^2}{4n}\Big).
\end{equation}

For the one-dimensional lower bounds of Section~\ref{sec:lower}, let $\mu$ be uniform on $(0,1)$ and $N=2n$. Since the source of each point is selected before the point is seen, the locations $U_1,\dots,U_N$ in reveal order are independent uniform variables, independent of the seed $\omega$. Assign the coefficient $c_i=\eps_i$ to a kept point and $c_i=0$ to every discarded or previously unrevealed point. Then $c_i$ is a measurable function of $(\omega,U_1,\dots,U_i)$, and
\begin{equation}\label{eq:controller-identities}
R=\#\{i\le N:c_i=0\},\qquad S_t(b)=\sum_{i\le t}c_i\one\{U_i\le b\},\qquad0\le t\le N,
\end{equation}
while the discrepancies are still given by \eqref{eq:St}. A \emph{controller} may choose any $c_i\in\{-1,0,1\}$ after seeing $U_i$, using only the seed and the locations seen so far, and pays one deletion for $c_i=0$. This class contains the thinning algorithms, with the same discrepancies and deletion count, but it also allows the sign to be chosen after the location is seen. The lower bounds in Section~\ref{sec:lower} will be proved for controllers under the corresponding expected deletion budgets.

\subsection{One-dimensional pilot spacings}\label{sec:pilot}

A \emph{pilot sample} is an initial set of observations which is discarded and used only to choose interval boundaries. In uniform coordinates, write its order statistics as
\[
0=a_0<a_1<\cdots<a_p<a_{p+1}=1.
\]
The intervals $(a_{s-1},a_s)$, $1\le s\le p+1$, are called \emph{cells}. All later observations avoid their endpoints almost surely.

We use two facts about these cells. First, consecutive groups of sufficiently many cells have nearly equal mass. Suppose $p+1=m2^h$, where $2\le m\le n$ and $2^h\ge400\log(2n)$, and put
\[
I_j=(a_{(j-1)2^h},a_{j2^h}),\qquad1\le j\le m.
\]
By the normalized-exponential representation of uniform spacings \cite[Section 4.1]{Pyk65}, the $p+1$ spacings $a_s-a_{s-1}$ have the joint law of $(\xi_1,\dots,\xi_{p+1})/\sum_{s\le p+1}\xi_s$ for independent exponential variables $\xi_s$ of mean one. Grouping consecutive blocks of $2^h$ spacings, the vector of bin masses satisfies
\[
(|I_1|,\dots,|I_m|)\overset d=\frac{(\Gamma_1,\dots,\Gamma_m)}{\Gamma_1+\cdots+\Gamma_m},\qquad\Gamma_j\text{ independent with law }\Gam(2^h,1).
\]
For $\Gamma\sim\Gam(\nu,1)$ and $\theta<1$ we have $\E e^{\theta\Gamma}=(1-\theta)^{-\nu}$, so Markov's inequality gives $\Prob(\Gamma\ge5\nu/4)\le e^{-5\nu\theta/4}(1-\theta)^{-\nu}$ and $\Prob(\Gamma\le3\nu/4)\le e^{3\nu\theta/4}(1+\theta)^{-\nu}$ for $\theta>0$. Choosing $\theta=1/5$ and $\theta=1/3$, respectively, gives the exponents $\nu(\frac14-\log\frac54)$ and $\nu(\log\frac43-\frac14)$, both at least $\nu/40$. Hence $\Prob(\Gamma_j\notin[3\cdot2^h/4,5\cdot2^h/4])\le2e^{-2^h/40}\le2(2n)^{-10}$, and a union bound over the $m\le n$ bins gives
\begin{equation}\label{eq:pilot-bins}
\Prob\Big(\frac1{2m}\le|I_j|\le\frac2m\text{ for every }1\le j\le m\Big)\ge1-(2n)^{-9},
\end{equation}
since on the displayed Gamma event each normalized mass lies between $(3/4)/(5m/4)=3/(5m)$ and $(5/4)/(3m/4)=5/(3m)$.

Second, a fresh point is unlikely to fall in a cell which already contains many later observations. This estimate requires no high-probability condition on the pilot.

\begin{lemma}\label{lem:cell}
Let $p\ge1$ pilot points, $t\ge0$ further points, and a fresh point be independent and uniform on $(0,1)$. Let $K$ count the further points in the pilot cell containing the fresh point. For every integer $r\ge1$,
\begin{equation}\label{eq:cell}
\Prob(K\ge r)\le2\Big(\frac t{p+t}\Big)^{\lceil r/2\rceil}.
\end{equation}
\end{lemma}

\begin{proof}
Order all $p+t+1$ points, distinguishing the fresh point. Its cell is bounded by the nearest pilot points on either side, or by the endpoints of $(0,1)$. If the cell contains at least $r$ further points, then on at least one side of the fresh point the first $\lceil r/2\rceil$ points in rank order are further points. Conditional on the rank of the fresh point, the labels of all other ranks form a uniformly random arrangement of $p$ pilot labels and $t$ further-point labels, by exchangeability. The probability that the first $a$ ranks on one specified side all carry further-point labels is at most
\[
\frac{t(t-1)\cdots(t-a+1)}{(p+t)(p+t-1)\cdots(p+t-a+1)}\le\Big(\frac t{p+t}\Big)^a,
\]
and it is zero if that side has fewer than $a$ ranks or if $t<a$. Take $a=\lceil r/2\rceil$ and sum over the two sides.
\end{proof}

For the upper bounds in dimension $k$, use the cell limit
\begin{equation}\label{eq:cell-limit}
r=\Big\lceil\frac{84}{\rho}\log\frac{200k}{\rho}\Big\rceil,
\qquad 0<\rho<2.
\end{equation}
In the signed stream, suppose that $p\ge\rho n/20$ pilot observations are followed by a stage of at most $2n$ further observations. Apply the limit to each marginal cell, discarding a point if one of its cells already contains $r$ earlier arrivals from this stage. Count every arrival, including discarded points. Lemma~\ref{lem:cell} and a union bound over the $k$ coordinates give, at each step,
\begin{equation}\label{eq:cell-deletion}
\Prob(\text{cell-limit deletion})
\le2k\Big(1+\frac\rho{40}\Big)^{-r/2}
\le2k e^{-\rho r/84}\le\frac\rho{100}.
\end{equation}
Here $\log(1+\rho/40)\ge\rho/(40+\rho)\ge\rho/42$. Thus the expected number of cell-limit deletions is at most $\rho n/50$. Each marginal cell supplies at most $r$ kept points during this stage.

\subsection{The Metropolis rule}\label{sec:metropolis}

For $A>0$ and $d\ge1$ put
\[
f_A(z)=\frac1A\cos^2\Big(\frac{\pi z}{2A}\Big)\one\{|z|<A\},\qquad\varphi_A=f_A^{\otimes d},
\]
so that $f_A$ is a probability density on $\R$, continuously differentiable, even, with $f_A'(z)=-\frac\pi{2A^2}\sin(\pi z/A)$ for $|z|<A$ and $f_A'(z)=0$ for $|z|\ge A$, and $\varphi_A$ is a probability density on $\R^d$ supported in the cube $(-A,A)^d$. Smirnov and Vershynin \cite[Theorem 1.1, Corollary 1.3 and Section 2]{SV26b} use this stationary density to confine a signed walk to a cube while charging each proposal according to its Euclidean length. The following lemma states their rule and its rejection estimate in the form used here.

\begin{lemma}[{\cite[Corollary 1.3 and Section 2]{SV26b}}]\label{lem:metropolis}
Let $d\ge1$, $A>0$, and let $v_1,\dots,v_M\in\R^d$ be fixed vectors. Let $\eps_1,\dots,\eps_M$ be fair signs and $\eta_1,\dots,\eta_M$ uniform variables on $(0,1)$, with all these variables mutually independent. Independently of them, let $w_0$ have density $\varphi_A$. For $t=1,\dots,M$ put $w_t=w_{t-1}+\eps_tv_t$ if
\[
\eta_t\le\frac{\varphi_A(w_{t-1}+\eps_tv_t)}{\varphi_A(w_{t-1})},
\]
in which case step $t$ is called accepted, and $w_t=w_{t-1}$ otherwise. Then the following hold.

(i) $w_t$ has density $\varphi_A$ for every $0\le t\le M$.

(ii) $w_t-w_0=\sum_{s\le t\ \mathrm{accepted}}\eps_sv_s$, and $\|w_t-w_0\|_\infty<2A$ almost surely, for every $0\le t\le M$.

(iii) For every $1\le t\le M$, the probability that step $t$ is rejected is at most $\dfrac{\pi}{2A}\|v_t\|_2$.
\end{lemma}

Only the coordinates in the support of $v_t$ enter the acceptance ratio, since $\varphi_A$ is a product density. In the applications below the vectors $v_t$ are indicator vectors of small sets of coordinates, and they become fixed after conditioning on all locations, while the source signs remain independent and fair. Under this conditioning the rejection probability in (iii) is averaged over the stationary state as well as over the sign and the acceptance variable. Since the coordinate $w_0(\ell)$ enters the computation only when $\ell$ lies in the support of some proposal, it may be sampled from $f_A$ at that moment, which gives the same joint law as sampling $w_0$ at the start.

\section{The one-dimensional terminal upper bound}\label{sec:terminal-upper}

Throughout this section $k=1$ and $0<\rho<1/8$. The construction uses only intervals fixed by a pilot sample. The pilot fixes the bins, the main stage thins inside them, and the last arrivals cancel the signed total of each bin.

\subsection{Fixed dyadic blocks}\label{sec:blocks}

Suppose that a bin consists of $2^h$ consecutive pilot cells, numbered $1,\dots,2^h$ from left to right, with $h\ge1$. For $0\le i\le h$ and $1\le b\le2^{h-i}$ the set of cells $\{(b-1)2^i+1,\dots,b2^i\}$ is called a \emph{dyadic block} of level $i$. Level $0$ consists of the cells, and level $h$ consists of the whole bin. Each cell lies in exactly one block of each level. The pilot ranks determine all blocks before the main stage.

Give each block a coordinate of the Metropolis state of Lemma~\ref{lem:metropolis}. When a point arrives during the main stage, first count the earlier main-stage arrivals, kept or not, in its cell. Discard the point if this count is at least $r$. Otherwise the point proposes the indicator vector of the $h+1$ blocks containing its cell, and the rule of Lemma~\ref{lem:metropolis} is applied with the sign of the point. Keep the point exactly when the proposal is accepted. Since the cell limit counts all main-stage arrivals, the sequence of proposal vectors is determined by the locations alone. The state is initialized with independent $f_A$ coordinates after the pilot stage.

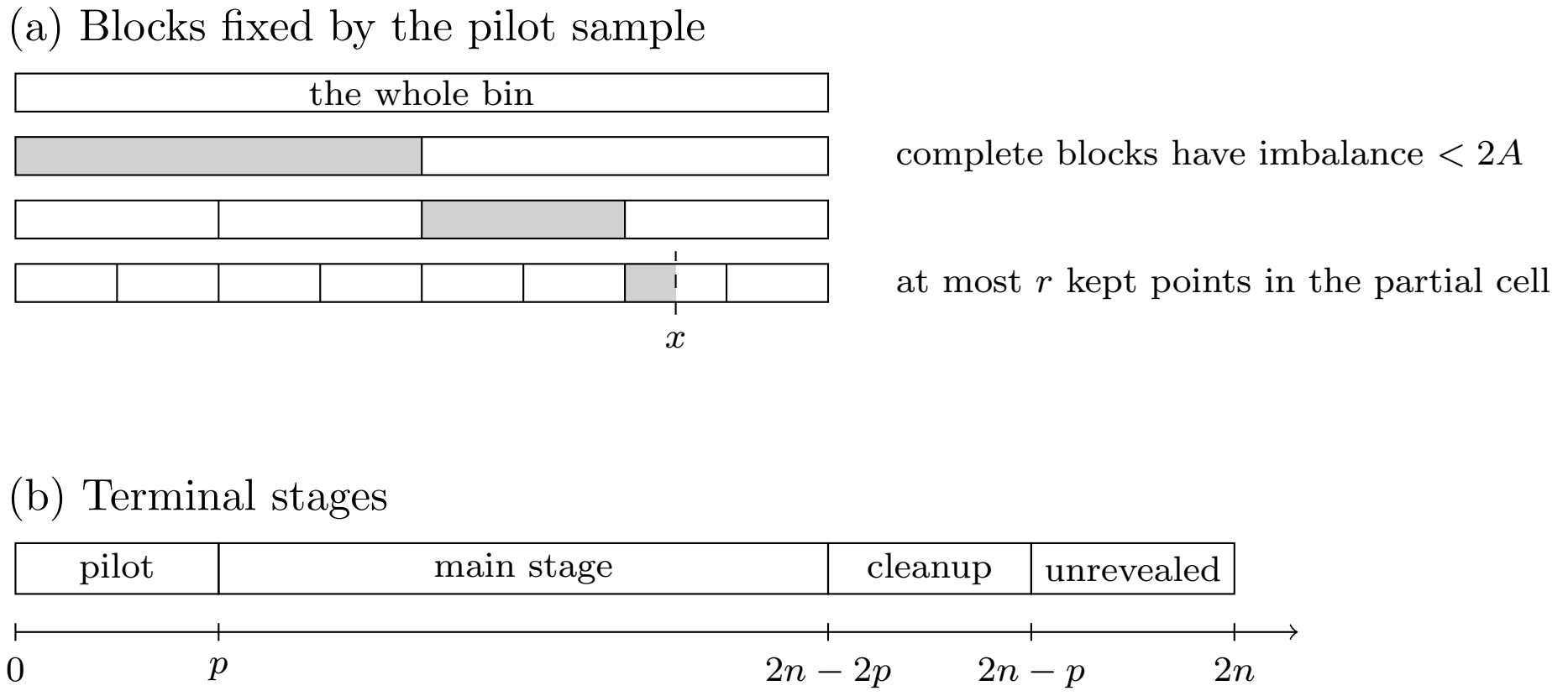
\begin{figure}[t]
\centering
\begin{tikzpicture}[x=1cm,y=1cm]
\begin{scope}
\node[font=\small,anchor=west] at (-0.2,2.55) {(a) Blocks fixed by the pilot sample};
\draw (0,1.9) rectangle (6.4,2.2); \node[font=\scriptsize] at (3.2,2.05) {the whole bin};
\fill[black!18] (0,1.4) rectangle (3.2,1.7);
\draw (0,1.4) rectangle (6.4,1.7); \draw (3.2,1.4) -- (3.2,1.7);
\fill[black!18] (3.2,0.9) rectangle (4.8,1.2);
\draw (0,0.9) rectangle (6.4,1.2); \foreach \x in {1.6,3.2,4.8} \draw (\x,0.9) -- (\x,1.2);
\fill[black!18] (4.8,0.4) rectangle (5.2,0.7);
\draw (0,0.4) rectangle (6.4,0.7); \foreach \x in {0.8,1.6,2.4,3.2,4.0,4.8,5.6} \draw (\x,0.4) -- (\x,0.7);
\draw[dashed] (5.2,0.3) -- (5.2,0.8); \node[font=\scriptsize,below] at (5.2,0.3) {$x$};
\node[font=\scriptsize,anchor=west,align=left] at (6.8,1.55) {complete blocks have imbalance $<2A$};
\node[font=\scriptsize,anchor=west,align=left] at (6.8,0.55) {at most $r$ kept points in the partial cell};
\end{scope}
\begin{scope}[yshift=-1.9cm]
\node[font=\small,anchor=west] at (-0.2,0.75) {(b) Terminal stages};
\draw (0,0) rectangle (1.6,0.4); \node[font=\scriptsize] at (0.8,0.2) {pilot};
\draw (1.6,0) rectangle (6.4,0.4); \node[font=\scriptsize] at (4.0,0.2) {main stage};
\draw (6.4,0) rectangle (8.0,0.4); \node[font=\scriptsize] at (7.2,0.2) {cleanup};
\draw (8.0,0) rectangle (9.6,0.4); \node[font=\scriptsize] at (8.8,0.2) {unrevealed};
\draw[->] (0,-0.3) -- (10.1,-0.3);
\foreach \x/\l in {0/$0$,1.6/$p$,6.4/$2n-2p$,8.0/$2n-p$,9.6/$2n$} {\draw (\x,-0.37) -- (\x,-0.23); \node[font=\scriptsize,below] at (\x,-0.37) {\l};}
\end{scope}
\end{tikzpicture}
\caption{Fixed blocks and terminal cancellation. (a) A bin of $2^h=8$ cells and its blocks of levels $3,2,1,0$. A threshold $x$ leaves complete dyadic blocks, shown shaded, and at most one partial cell. (b) The pilot, cleanup and unread tail each contain $p$ observations. Stage lengths are schematic. When the cleanup balances every bin, all complete bins preceding a threshold contribute zero.}
\label{fig:terminal}
\end{figure}

This applies the dyadic representation of \cite{SV26a} to empirical quantiles. The cell limit controls thresholds that fall between two consecutive pilot points. Figure~\ref{fig:terminal} shows the blocks of one bin and the stages of the terminal algorithm.
\begin{algorithm}[t]
\caption{One step of the cell limit and the Metropolis rule}
\label{alg:step}
\begin{algorithmic}[1]
\Procedure{Step}{$x,\varepsilon$}
\Comment{arrival counts $K_a(\cdot)$ of the marginal cells; state $w$}
  \State $c\gets(c_1(x),\dots,c_k(x))$, the cell indices of the tagged point $x$
  \If{$K_a(c_a)\ge r$ for some $a$}
  \Comment{cell limit}
    \State $K_a(c_a)\gets K_a(c_a)+1$ for $a=1,\dots,k$
    \State \Return \Discard
    \Comment{the state is unchanged}
  \EndIf
  \State $K_a(c_a)\gets K_a(c_a)+1$ for $a=1,\dots,k$
  \State $\mathcal P_x\gets$ the dyadic boxes containing $c$
  \Comment{the $h+1$ blocks of its bin when $k=1$}
  \State draw $w(B)$ from $f_A$ for every $B\in\mathcal P_x$ not used before
  \State draw $\eta$ uniformly from $(0,1)$
  \If{$\eta\prod_{B\in\mathcal P_x}f_A\big(w(B)\big)\le\prod_{B\in\mathcal P_x}f_A\big(w(B)+\varepsilon\big)$}
    \State $w(B)\gets w(B)+\varepsilon$ for every $B\in\mathcal P_x$
    \State \Return \Keep
  \Else
    \State \Return \Discard
  \EndIf
\EndProcedure
\end{algorithmic}
\end{algorithm}

\begin{lemma}\label{lem:interval}
At every time during the main stage, every dyadic block, including the whole bin, has signed retained count of absolute value less than $2A$, and every initial part of the bin has signed retained count of absolute value at most
\begin{equation}\label{eq:interval}
2Ah+r.
\end{equation}
\end{lemma}

\begin{proof}
A kept point changes the coordinate of a block $B$ by its sign exactly when its cell lies in $B$, and no other point changes this coordinate. Hence the signed retained count of $B$ at time $t$ is $w_t(B)-w_0(B)$, whose absolute value is less than $2A$ by Lemma~\ref{lem:metropolis}(ii).

Fix a threshold $x$ inside the bin. The binary expansion of the number of complete cells preceding $x$ partitions them into at most $h$ dyadic blocks. Each contributes less than $2A$ in absolute value. The partial cell contributes at most $r$, since at most $r$ main-stage arrivals in a cell are offered to the Metropolis rule. This proves \eqref{eq:interval}.
\end{proof}

\subsection{Parameters and algorithm}\label{sec:algorithm}

Let
\begin{equation}\label{eq:par1}
L=\log(2n),\qquad m=\Big\lfloor\frac{\rho^2n}{10^4L}\Big\rfloor.
\end{equation}
If $m<2$, reveal the samples in a fixed order and keep all points. Otherwise put
\begin{equation}\label{eq:par2}
h=\Big\lceil\log_2\frac{1200L}\rho\Big\rceil,\qquad p=m2^h-1,\qquad A=\frac{8\pi\sqrt{h+1}}\rho.
\end{equation}
Use $r$ from \eqref{eq:cell-limit} with $k=1$. There are $m$ bins, each containing $2^h$ pilot cells. The pilot, cleanup and unread tail each have length $p$. The Metropolis scale is $A$, and each cell supplies at most $r$ proposals.

The algorithm plans to reveal $2n-p$ points using independent fair source coins as in Section~\ref{sec:general-model}, and stops if a coin selects an exhausted sample. Its stages are as follows.

\smallskip
\noindent\emph{Pilot, steps $1,\dots,p$.} Discard these points and sort them. Their $p+1=m2^h$ cells form $m$ bins $I_1,\dots,I_m$ of $2^h$ consecutive cells each, as in Section~\ref{sec:pilot}. Build the fixed dyadic blocks of every bin.

\smallskip
\noindent\emph{Main stage, steps $p+1,\dots,2n-2p$.} Apply the cell limit and the Metropolis rule of Section~\ref{sec:blocks}, with the state initialized independently across the block coordinates of all bins. Record the signed retained count $d_j$ of each bin.

\smallskip
\noindent\emph{Cleanup, steps $2n-2p+1,\dots,2n-p$.} For a point in bin $j$ with sign $\eps$, keep it exactly when $\eps d_j<0$, and then replace $d_j$ by $d_j+\eps$. Once $d_j=0$, keep no more points in that bin. All unrevealed observations count as discarded.
\begin{algorithm}[t]
\caption{Terminal thinning on the line}
\label{alg:terminal}
\begin{algorithmic}[1]
\Require $n\ge2$ and $0<\rho<1/8$; online access to $X_1,\dots,X_n$ and $Y_1,\dots,Y_n$
\State $L\gets\log(2n)$, \ $m\gets\lfloor\rho^2n/(10^4L)\rfloor$
\If{$m<2$}
  \State reveal $X_1,\dots,X_n$, then $Y_1,\dots,Y_n$, keep every point, and stop
\EndIf
\State $h\gets\lceil\log_2(1200L/\rho)\rceil$, \ $p\gets m2^h-1$, \ $A\gets8\pi\sqrt{h+1}/\rho$, \ $r\gets\lceil(84/\rho)\log(200/\rho)\rceil$
\For{$i=1,\dots,2n-p$}
  \State toss a fair sign $\varepsilon_i$; select the $X$-sample if $\varepsilon_i=+1$ and the $Y$-sample otherwise
  \If{the selected sample is exhausted}
    \State \textbf{stop}
  \EndIf
  \State reveal its next observation $x_i$ and draw a tag $u_i$ uniformly from $(0,1)$
  \Statex \hspace{\algorithmicindent}\textit{comparisons below use the tagged points $(x_i,u_i)$ in lexicographic order}
  \If{$i\le p$}
  \Comment{pilot}
    \State discard $x_i$
    \If{$i=p$}
      \State sort the pilot; form the $m$ bins of $2^h$ consecutive cells and their dyadic blocks
      \State $d_j\gets0$ for $j=1,\dots,m$
    \EndIf
  \ElsIf{$i\le2n-2p$}
  \Comment{main stage}
    \State let $j$ be the bin containing $x_i$
    \If{\Call{Step}{$x_i,\varepsilon_i$} $=$ \Keep}
      \State keep $x_i$ and set $d_j\gets d_j+\varepsilon_i$
    \Else
      \State discard $x_i$
    \EndIf
  \Else
  \Comment{cleanup}
    \State let $j$ be the bin containing $x_i$
    \If{$\varepsilon_id_j<0$}
      \State keep $x_i$ and set $d_j\gets d_j+\varepsilon_i$
    \Else
      \State discard $x_i$
    \EndIf
  \EndIf
\EndFor
\State discard the observations that were never revealed
\end{algorithmic}
\end{algorithm}

\smallskip
For general inputs the comparisons use the independent tags of Section~\ref{sec:tags}. We record the parameter estimates needed in the proof. When $m\ge2$,
\begin{equation}\label{eq:parfacts}
\begin{gathered}
\frac{\rho n}{20}\le p\le\frac{\rho n}4,\qquad
\rho^2n\ge2\cdot10^4L,\qquad L\ge14,\\
h+1\le4L,\qquad
4A\le\frac{150L}\rho\le\frac p{4m}.
\end{gathered}
\end{equation}
Indeed, $m\ge2$ gives $\rho^2n\ge2\cdot10^4L$, and the floor and ceiling in the definitions imply
\[
\frac{\rho^2n}{2\cdot10^4L}\le m\le\frac{\rho^2n}{10^4L},\qquad
\frac{1200L}{\rho}\le2^h<\frac{2400L}{\rho}.
\]
Consequently $3\rho n/50-1\le p<6\rho n/25$. Since $\rho<1/8$, we also have $\rho n\ge16\cdot10^4L\ge100$, which gives the asserted bounds on $p$. The same assumptions yield $2n\ge2.56\cdot10^6L$, so $L\ge14$, and
\[
\frac1\rho\le\frac{e^{L/2}}{200\sqrt L},\qquad
h+1\le\frac{L/2+(\log L)/2+\log6}{\log2}+2\le4L.
\]
The last inequality follows from $L\ge14$. Finally,
\[
4A\le\frac{64\pi\sqrt L}\rho\le\frac{150L}\rho
\le\frac{2^h}{8}\le\frac p{4m}.
\]
Since $h+1\le C(\log L+\log(2/\rho))$ and $r\le84\rho^{-1}\log(200/\rho)+1$, the definitions give
\begin{equation}\label{eq:explicit-constant}
2A(h+1)+r+1\le\frac C\rho\Big(\log L+\log\frac2\rho\Big)^{3/2}
\end{equation}
for an absolute constant $C$.

\subsection{Expected deletions}\label{sec:deletions}

Assume $m\ge2$ and first run the procedure on the independent signed stream for the planned $2n-p$ points. The main stage has at most $2n$ arrivals. Since $p\ge\rho n/20$, \eqref{eq:cell-deletion} bounds the expected number of cell-limit deletions by $\rho n/50$.

Now condition on all locations of the stream. The pilot partition, the cell-limit decisions and the proposal vectors at the eligible steps are then fixed, every proposal has exactly $h+1$ unit coordinates, and the source signs remain independent and fair. Lemma~\ref{lem:metropolis}(iii) bounds the rejection probability at every eligible step by
\[
\frac{\pi\sqrt{h+1}}{2A}=\frac\rho{16},
\]
so the expected number of Metropolis rejections is at most $\rho n/8$. The estimate averages over the stationary state and the independent sign and acceptance variable.

The pilot and the cleanup together discard at most $2p$ points, and the planned unrevealed tail has length $p$. By \eqref{eq:source-failure} and \eqref{eq:parfacts},
\begin{equation}\label{eq:stream-failure}
\Prob(E_p^c)\le2e^{-p^2/(4n)}\le2e^{-\rho^2n/1600}\le2e^{-12.5L}\le2(2n)^{-12}.
\end{equation}
On $E_p$ the run on the samples coincides with the stream run, and the number of unrevealed observations is exactly $p$. Charging $2n$ deletions on the complementary event, we obtain
\begin{equation}\label{eq:deletions}
\E R\le3p+\frac{\rho n}{50}+\frac{\rho n}8+4n(2n)^{-12}\le\frac34\rho n+\frac{\rho n}{50}+\frac{\rho n}8+1\le\frac9{10}\rho n+1\le\rho n,
\end{equation}
using $p\le\rho n/4$ and $\rho n\ge100$. If $m<2$, all points are kept and the budget holds trivially.

\subsection{Terminal cancellation}\label{sec:cancellation}

In the stream run, let $G$ be the event in \eqref{eq:pilot-bins}. Since $2^h\ge1200L/\rho\ge400L$ and $m\le n$, we have $\Prob(G^c)\le(2n)^{-9}$.

At the start of the cleanup, Lemma~\ref{lem:interval} applied to the whole bin gives $|d_j|<2A$ for each bin. Each cleanup point kept in bin $j$ has sign $-\sgn(d_j)$ and reduces $|d_j|$ by one until it reaches zero. Thus fewer than $2A$ cleanup points are kept in any bin.

Condition on the complete history of the stream run through the main stage, including the pilot points, the locations and source signs of the main-stage points and the Metropolis randomness. The bins and the counts $d_j$ are then fixed, while the $p$ cleanup points of the stream have independent uniform locations and independent fair signs. For a bin with $d_j\ne0$, let $N_j$ count the cleanup arrivals in $I_j$ with sign $-\sgn(d_j)$. These signs remain useful until balance is reached, so $N_j\ge|d_j|$ implies that bin $j$ is balanced at the end of the cleanup, and
\[
N_j\mid\text{main-stage history}\ \sim\ \Bin(p,|I_j|/2).
\]
On $G$ its mean is at least $p/(4m)\ge4A$ by \eqref{eq:parfacts}. The Chernoff bound $\Prob(X\le\E X-u)\le\exp(-u^2/(2\E X))$ for a binomial variable $X$, see \cite[Theorem 2.1]{JLR00}, applied with $u=\E X/2$, gives
\[
\Prob(N_j<2A\mid\text{main-stage history})\le e^{-\E N_j/8}\le e^{-p|I_j|/16}\le e^{-p/(32m)}\le e^{-150L/(8\rho)}\le(2n)^{-10},
\]
where we used $2A\le\frac12\E N_j$ in the first step and $p/(32m)\ge150L/(8\rho)\ge150L$ in the last two. Bins with $d_j=0$ stay balanced. A union bound over the at most $n$ bins with $d_j\ne0$, followed by averaging over the history on $G$, shows that the probability that some bin is unbalanced at the end of the cleanup, and $G$ occurs, is at most $n(2n)^{-10}\le(2n)^{-9}$.

Suppose that all bins are balanced at the end of the cleanup, and let $x\in I_j$ be a threshold. Pilot points are discarded, and later arrivals avoid them almost surely. The kept points with location at most $x$ consist of the kept points of the complete bins $I_1,\dots,I_{j-1}$, whose signed counts are zero, and of the kept points of the initial part of $I_j$ at $x$. The main-stage points of this initial part contribute at most $2Ah+r$ in absolute value by Lemma~\ref{lem:interval}, applied at the end of the main stage, and the cleanup points of bin $j$ contribute less than $2A$ in absolute value. Therefore, in uniform coordinates and hence by \eqref{eq:tags-dominate} for the original half-lines,
\begin{equation}\label{eq:det-bound}
D\le2A(h+1)+r.
\end{equation}
On $E_p$ the run on the samples keeps the same points as the stream run. Combining the two exceptional events above with \eqref{eq:stream-failure} yields, for the run on the samples,
\begin{equation}\label{eq:D-whp}
\Prob\big(D>2A(h+1)+r\big)\le2(2n)^{-12}+2(2n)^{-9}\le4(2n)^{-9}.
\end{equation}
Since always $D\le2n$, it follows that
\begin{equation}\label{eq:D-mean}
\E D\le2A(h+1)+r+2n\cdot4(2n)^{-9}\le2A(h+1)+r+1.
\end{equation}

\begin{proof}[Proof of Theorem~\ref{thm:terminal}(i)]
The algorithm uses only comparisons of tagged observations and the parameters $n$ and $\rho$. The deletion bound is \eqref{eq:deletions}. When $m\ge2$, the estimates \eqref{eq:explicit-constant}, \eqref{eq:D-whp} and \eqref{eq:D-mean} prove both discrepancy assertions, with an absolute constant $C$.

When $m<2$, all points are kept. Let $F_n^X$ and $F_n^Y$ be the empirical distribution functions of the two samples and $F$ that of $\mu$. The Dvoretzky--Kiefer--Wolfowitz inequality with Massart's constant \cite{Mas90} gives $\Prob(\sqrt n\sup_b|F_n^X(b)-F(b)|>\lambda)\le2e^{-2\lambda^2}$ for all $\lambda>0$, hence $\E\sup_b|F_n^X-F|\le n^{-1/2}\int_0^\infty2e^{-2\lambda^2}\,d\lambda=\sqrt{\pi/(2n)}$, and the triangle inequality for $D=n\sup_b|F_n^X(b)-F_n^Y(b)|$ gives
\[
\E D\le\sqrt{2\pi n}.
\]
The inequality $\rho^2n<2\cdot10^4L$ also implies $L<11\log(2/\rho)$. Indeed, put $u=\log(2/\rho)\ge\log16$. Then $e^L=2n<4\cdot10^4L/\rho^2=10^4Le^{2u}$. If we had $L\ge11u$, then, since $x\mapsto e^x/x$ increases for $x>1$, we would get $e^{11u}/(11u)\le e^L/L<10^4e^{2u}$, that is $e^{9u}<11\cdot10^4u$, which fails for $u\ge\log16$ because $9u-\log(11\cdot10^4u)$ is increasing and positive at $u=\log16$. Consequently
\[
\E D\le\sqrt{2\pi n}<\frac{200\sqrt\pi}\rho\sqrt L<\frac{200\sqrt{11\pi}}\rho\Big(\log\frac2\rho\Big)^{1/2}\le\frac C\rho\Big(\log L+\log\frac2\rho\Big)^{3/2},
\]
using $2n<4\cdot10^4L/\rho^2$ in the second step and $\log(2/\rho)\ge1$ in the last. This proves the expectation bound for all $n\ge2$ with a single absolute constant. The condition in the high-probability assertion is exactly $m\ge2$.
\end{proof}

Sorting the pilot takes $O(p\log p)$ comparisons. Each later point uses $O(\log n)$ comparisons to locate its cell and $O(h)$ arithmetic operations for its proposal. There are fewer than $2(p+1)$ dyadic blocks, so storage is $O(n)$. These counts assume exact real arithmetic, unit-cost evaluations of $f_A$, and independent sampling from $f_A$.

\section{Prefix upper bounds}\label{sec:prefix-upper}

Fix $k\ge1$, $0<\rho<2$ and $n\ge2$. If $\rho^2n<2\cdot10^4\log(2n)$, reveal the samples in a fixed order and keep every point. Otherwise put
\begin{equation}\label{eq:multipar}
p=\Big\lfloor\frac{\rho n}8\Big\rfloor,\qquad
h=\lceil\log_2(p+1)\rceil,\qquad
s=(h+1)^k,\qquad A=\frac{8\pi\sqrt s}{\rho},
\end{equation}
and take $r$ from \eqref{eq:cell-limit}. All comparisons use the coordinatewise tags of Section~\ref{sec:tags}.

Plan to reveal $2n-p$ points using the source coins of Section~\ref{sec:general-model}, and stop if a selected source is exhausted. Discard the first $p$ observations as the pilot. In each coordinate, pad the $p+1$ pilot cells with empty cells to obtain a list of length $2^h$, and use the dyadic blocks of Section~\ref{sec:blocks} on this list. A \emph{dyadic box} is a product of one block from each coordinate. Let $\cP$ be the set of these boxes and index a Metropolis state by $\cP$, with independent initial coordinates of density $f_A$.

For a point $x$, let $c(x)=(c_1(x),\dots,c_k(x))$ record its cell indices. Number the arrivals after the pilot by $t$ and write $\eps_t$ for their source signs. The proposal vector $v_t\in\{0,1\}^{\cP}$ indicates the $s$ dyadic boxes containing $c(x)$. Apply the cell limit of Section~\ref{sec:pilot} in every coordinate, counting all earlier arrivals after the pilot. If the point is discarded by the cell limit, leave the state unchanged, $w_t=w_{t-1}$. Otherwise apply the Metropolis rule of Lemma~\ref{lem:metropolis} to $\eps_tv_t$ and keep the point exactly when the proposal is accepted. All unread observations are discarded.
\begin{algorithm}[t]
\caption{Prefix thinning in $\mathbb{R}^k$}
\label{alg:prefix}
\begin{algorithmic}[1]
\Require $k\ge1$, $n\ge2$ and $0<\rho<2$; online access to $X_1,\dots,X_n$ and $Y_1,\dots,Y_n$ in $\mathbb{R}^k$
\If{$\rho^2n<2\cdot10^4\log(2n)$}
  \State reveal $X_1,\dots,X_n$, then $Y_1,\dots,Y_n$, keep every point, and stop
\EndIf
\State $p\gets\lfloor\rho n/8\rfloor$, \ $h\gets\lceil\log_2(p+1)\rceil$, \ $s\gets(h+1)^k$, \ $A\gets8\pi\sqrt s/\rho$, \ $r\gets\lceil(84/\rho)\log(200k/\rho)\rceil$
\For{$i=1,\dots,2n-p$}
  \State toss a fair sign $\varepsilon_i$; select the $X$-sample if $\varepsilon_i=+1$ and the $Y$-sample otherwise
  \If{the selected sample is exhausted}
    \State \textbf{stop}
  \EndIf
  \State reveal its next observation $x_i$; draw independent uniform tags $u_{i1},\dots,u_{ik}$
  \Statex \hspace{\algorithmicindent}\textit{coordinate $a$ is compared through the tagged values $(x_{ia},u_{ia})$}
  \If{$i\le p$}
  \Comment{pilot}
    \State discard $x_i$
    \If{$i=p$}
      \State sort the pilot in each coordinate; pad its $p+1$ cells with empty cells to $2^h$
      \State form the dyadic blocks of each coordinate; the dyadic boxes are their products
    \EndIf
  \ElsIf{\Call{Step}{$x_i,\varepsilon_i$} $=$ \Keep}
    \State keep $x_i$
  \Else
    \State discard $x_i$
  \EndIf
\EndFor
\State discard the observations that were never revealed
\end{algorithmic}
\end{algorithm}
\begin{figure}[t]
\centering
\begin{tikzpicture}[x=6.4cm,y=6.4cm]
\fill[black!7] (0.7951,0) rectangle (0.9092,1);
\fill[black!7] (0,0.6609) rectangle (1,0.7806);
\fill[black!16] (0,0) rectangle (0.5657,0.5483);
\fill[black!28] (0,0.5483) rectangle (0.5657,0.6609);
\fill[black!28] (0.5657,0) rectangle (0.7951,0.5483);
\fill[black!40] (0.5657,0.5483) rectangle (0.7951,0.6609);
\draw[black!45, line width=0.3pt] (0.1542,0) -- (0.1542,1);
\draw[black!45, line width=0.3pt] (0.2797,0) -- (0.2797,1);
\draw[black!45, line width=0.3pt] (0.4791,0) -- (0.4791,1);
\draw[black!45, line width=0.3pt] (0.5657,0) -- (0.5657,1);
\draw[black!45, line width=0.3pt] (0.7098,0) -- (0.7098,1);
\draw[black!45, line width=0.3pt] (0.7951,0) -- (0.7951,1);
\draw[black!45, line width=0.3pt] (0.9092,0) -- (0.9092,1);
\draw[black!45, line width=0.3pt] (0,0.1535) -- (1,0.1535);
\draw[black!45, line width=0.3pt] (0,0.3333) -- (1,0.3333);
\draw[black!45, line width=0.3pt] (0,0.4327) -- (1,0.4327);
\draw[black!45, line width=0.3pt] (0,0.5483) -- (1,0.5483);
\draw[black!45, line width=0.3pt] (0,0.6609) -- (1,0.6609);
\draw[black!45, line width=0.3pt] (0,0.7806) -- (1,0.7806);
\draw[black!45, line width=0.3pt] (0,0.8668) -- (1,0.8668);
\draw[thick] (0,0) rectangle (0.5657,0.5483);
\draw[thick] (0,0.5483) rectangle (0.5657,0.6609);
\draw[thick] (0.5657,0) rectangle (0.7951,0.5483);
\draw[thick] (0.5657,0.5483) rectangle (0.7951,0.6609);
\draw (0,0) rectangle (1,1);
\fill (0.3987,0.25) circle (0.9pt);
\fill (0.423,0.3002) circle (0.9pt);
\fill (0.8325,0.1907) circle (0.9pt);
\fill (0.0643,0.7592) circle (0.9pt);
\fill (0.6004,0.6245) circle (0.9pt);
\fill (0.2399,0.5918) circle (0.9pt);
\fill (0.4416,0.4722) circle (0.9pt);
\fill (0.2043,0.1031) circle (0.9pt);
\fill (0.7562,0.078) circle (0.9pt);
\fill (0.395,0.9646) circle (0.9pt);
\fill (0.2382,0.0933) circle (0.9pt);
\fill (0.4321,0.8059) circle (0.9pt);
\fill (0.7398,0.0895) circle (0.9pt);
\fill (0.0496,0.9578) circle (0.9pt);
\fill (0.4336,0.6821) circle (0.9pt);
\fill (0.6544,0.9663) circle (0.9pt);
\fill (0.5366,0.0939) circle (0.9pt);
\fill (0.5432,0.6278) circle (0.9pt);
\fill (0.6277,0.0272) circle (0.9pt);
\fill (0.9416,0.4027) circle (0.9pt);
\fill (0.3866,0.2441) circle (0.9pt);
\fill (0.0365,0.9443) circle (0.9pt);
\fill (0.9285,0.3873) circle (0.9pt);
\fill (0.3316,0.6125) circle (0.9pt);
\fill (0.6319,0.9158) circle (0.9pt);
\fill (0.419,0.2505) circle (0.9pt);
\fill (0.8832,0.4086) circle (0.9pt);
\fill (0.1169,0.0461) circle (0.9pt);
\fill (0.0496,0.7046) circle (0.9pt);
\fill (0.2485,0.9668) circle (0.9pt);
\fill (0.8149,0.5272) circle (0.9pt);
\fill (0.0466,0.1881) circle (0.9pt);
\fill (0.1103,0.1034) circle (0.9pt);
\fill (0.2221,0.6267) circle (0.9pt);
\fill (0.7666,0.4002) circle (0.9pt);
\fill (0.4591,0.7543) circle (0.9pt);
\fill (0.8198,0.5046) circle (0.9pt);
\fill (0.8383,0.919) circle (0.9pt);
\fill (0.0675,0.1107) circle (0.9pt);
\fill (0.5111,0.8943) circle (0.9pt);
\fill (0.6689,0.2361) circle (0.9pt);
\fill (0.0451,0.5081) circle (0.9pt);
\fill (0.6373,0.9084) circle (0.9pt);
\fill (0.3594,0.9631) circle (0.9pt);
\fill (0.1261,0.7489) circle (0.9pt);
\fill (0.945,0.75) circle (0.9pt);
\draw[line width=0.7pt] (0.8962,0.6479) -- (0.9222,0.6739) (0.8962,0.6739) -- (0.9222,0.6479);
\draw[line width=0.7pt] (0.1412,0.5353) -- (0.1672,0.5613) (0.1412,0.5613) -- (0.1672,0.5353);
\draw[line width=0.7pt] (0.2667,0.1405) -- (0.2927,0.1665) (0.2667,0.1665) -- (0.2927,0.1405);
\draw[line width=0.7pt] (0.4661,0.4197) -- (0.4921,0.4457) (0.4661,0.4457) -- (0.4921,0.4197);
\draw[line width=0.7pt] (0.5527,0.8538) -- (0.5787,0.8798) (0.5527,0.8798) -- (0.5787,0.8538);
\draw[line width=0.7pt] (0.6968,0.7676) -- (0.7228,0.7936) (0.6968,0.7936) -- (0.7228,0.7676);
\draw[line width=0.7pt] (0.7821,0.3203) -- (0.8081,0.3463) (0.7821,0.3463) -- (0.8081,0.3203);
\draw[very thick, densely dashed] (0.8621,0) -- (0.8621,0.7257) -- (0,0.7257);
\node[star, star points=5, fill=black, inner sep=1.3pt] at (0.8621,0.7257) {};
\node[font=\scriptsize, above right, inner sep=2pt, fill=white, fill opacity=0.85, text opacity=1] at (0.8621,0.7257) {$\tau$};
\node[font=\tiny, below] at (0.0771,0) {$1$};
\node[font=\tiny, left] at (0,0.0767) {$1$};
\node[font=\tiny, below] at (0.2169,0) {$2$};
\node[font=\tiny, left] at (0,0.2434) {$2$};
\node[font=\tiny, below] at (0.3794,0) {$3$};
\node[font=\tiny, left] at (0,0.383) {$3$};
\node[font=\tiny, below] at (0.5224,0) {$4$};
\node[font=\tiny, left] at (0,0.4905) {$4$};
\node[font=\tiny, below] at (0.6377,0) {$5$};
\node[font=\tiny, left] at (0,0.6046) {$5$};
\node[font=\tiny, below] at (0.7524,0) {$6$};
\node[font=\tiny, left] at (0,0.7207) {$6$};
\node[font=\tiny, below] at (0.8521,0) {$7$};
\node[font=\tiny, left] at (0,0.8237) {$7$};
\node[font=\tiny, below] at (0.9546,0) {$8$};
\node[font=\tiny, left] at (0,0.9334) {$8$};
\draw[decorate, decoration={brace, mirror, amplitude=3pt}] (0.005,-0.07) -- (0.5607,-0.07) node[midway, below=3pt, font=\scriptsize] {$\{1,\dots,4\}$};
\draw[decorate, decoration={brace, mirror, amplitude=3pt}] (0.5707,-0.07) -- (0.7901,-0.07) node[midway, below=3pt, font=\scriptsize] {$\{5,6\}$};
\draw[decorate, decoration={brace, amplitude=3pt}] (-0.07,0.005) -- (-0.07,0.5433) node[midway, left=3pt, font=\scriptsize] {$\{1,\dots,4\}$};
\draw[decorate, decoration={brace, amplitude=3pt}] (-0.07,0.5533) -- (-0.07,0.6559) node[midway, left=3pt, font=\scriptsize] {$\{5\}$};
\end{tikzpicture}
\caption{The comparison in Section~\ref{sec:prefix-upper} for $k=2$ and $p=7$, so that $h=3$ and each coordinate has $p+1=2^h=8$ cells, with no padding; the picture is the two-dimensional analogue of Figure~\ref{fig:terminal}(a). The cells are cut by the coordinates of the pilot points, shown as crosses, and the later arrivals are shown as dots. For the threshold $\tau$, the orthant $\{\widehat x\le\tau\}$, bounded by the dashed lines, contains the grid query $Q(\ell)$ with $\ell=(6,5)$, obtained by rounding each coordinate of $\tau$ down to a cell boundary. This query is the union of the four dyadic boxes with thick boundaries, built from the blocks $\{1,\dots,4\}$ and $\{5,6\}$ of the first coordinate and $\{1,\dots,4\}$ and $\{5\}$ of the second. The rest of the orthant lies in the two light strips formed by the marginal cells that contain $\tau_1$ and $\tau_2$, and the cell limit leaves at most $r$ kept points in each strip.}
\label{fig:grid}
\end{figure}
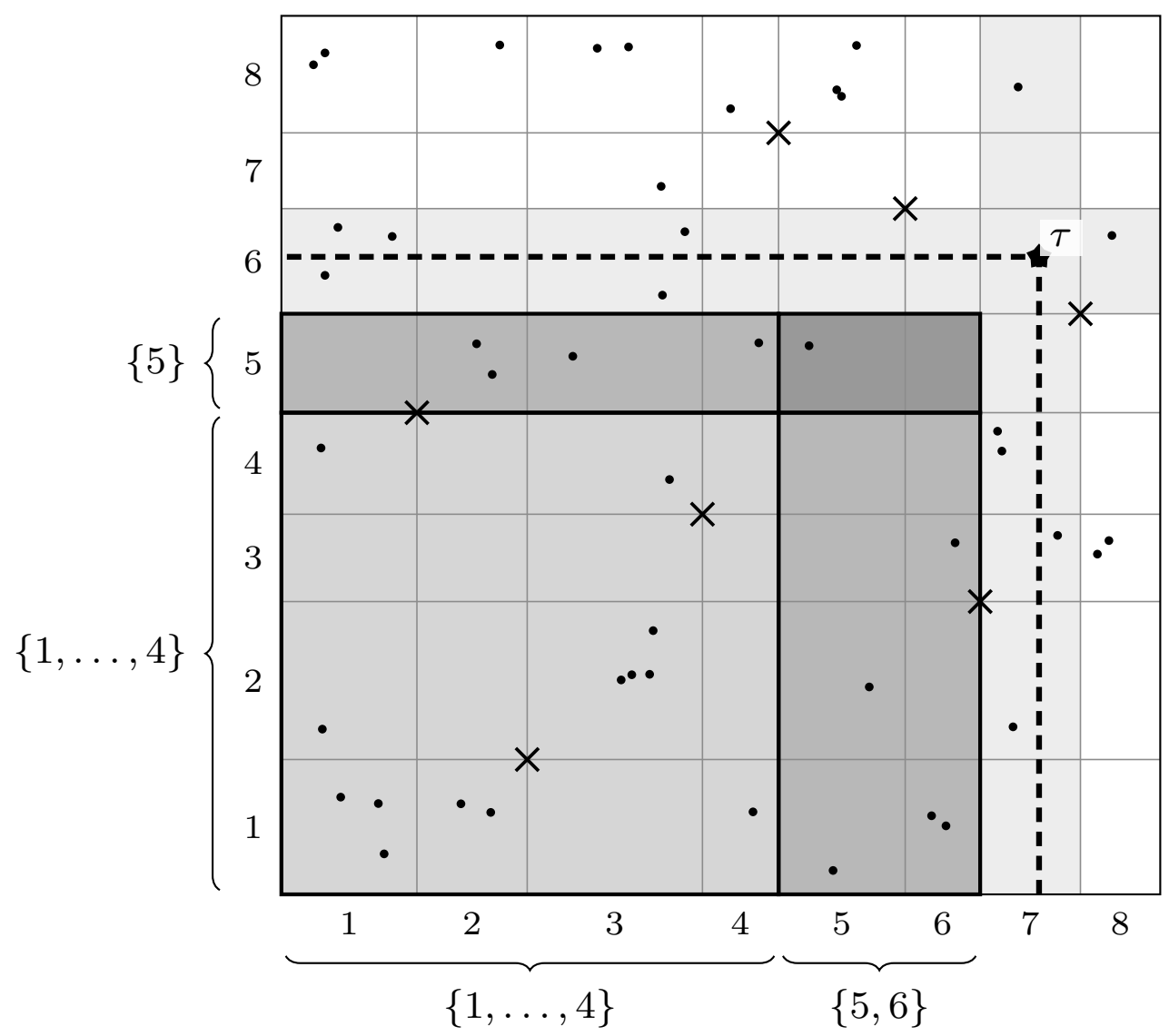

\begin{proof}[Proof of the upper bounds in Theorem~\ref{thm:multi}]
In the nontrivial branch, the threshold on $n$ and $\rho<2$ give $\rho n\ge10^4\log(2n)\ge100$. Hence
\[
\frac{\rho n}{16}\le p\le\frac{\rho n}8.
\]
Work first in the signed stream. The cell-limit deletions have expectation at most $\rho n/50$ by \eqref{eq:cell-deletion}. Conditional on all tagged locations, the cells, proposal vectors and the set of arrivals submitted to the Metropolis rule are fixed. The submitted source signs remain independent and fair. Lemma~\ref{lem:metropolis} applies to this subsequence, and each submitted proposal is rejected with probability at most $\pi\sqrt s/(2A)=\rho/16$. Holding the state fixed at every unsubmitted arrival shows that $w_j$ has density $\varphi_A$ at every post-pilot time $j$.

The pilot and unread tail contribute $2p$ deletions. On $E_p$ the sample run agrees with the stream run, and on $E_p^c$ we charge $2n$ deletions. Since
\[
\frac{p^2}{4n}\ge\frac{\rho^2n}{1024}\ge12\log(2n),
\]
the source bound \eqref{eq:source-failure} gives
\begin{equation}\label{eq:multideletions}
\E R\le2p+\frac{\rho n}{50}+\frac{\rho n}8+4n(2n)^{-12}
\le\frac25\rho n+1\le\rho n.
\end{equation}
The keep-all branch has $R=0$.

For $\ell\in\{0,\dots,2^h\}^k$, let $Q(\ell)$ be the grid query with cell indices $c_a\le\ell_a$ in every coordinate. Each initial cell segment is a disjoint union of at most $h+1$ dyadic blocks. For a segment shorter than $2^h$, use the binary decomposition in the proof of Lemma~\ref{lem:interval}. The whole list is itself one block. Taking products gives a family $\cD(\ell)\subseteq\cP$ of at most $s$ disjoint boxes with union $Q(\ell)$. Every kept point in the query contributes to exactly one of these boxes. Lemma~\ref{lem:metropolis}(ii) therefore identifies its signed count after $j$ post-pilot arrivals as
\begin{equation}\label{eq:grid-count}
\sum_{B\in\cD(\ell)}\big(w_j(B)-w_0(B)\big).
\end{equation}
Let
\[
\Xi=\max_{0\le j\le2n-2p}\ \max_{\ell\in\{0,\dots,2^h\}^k}
\Big|\sum_{B\in\cD(\ell)}\big(w_j(B)-w_0(B)\big)\Big|.
\]
Conditional on the locations, the coordinates of each fixed state $w_j$ are independent, centered and bounded by $A$ in absolute value. For every family $\cD$ of at most $s$ coordinates, Hoeffding's inequality \cite{Hoe63} gives
\[
\Prob\Big(\Big|\sum_{B\in\cD}w_j(B)\Big|\ge u\ \Big|\ \text{locations}\Big)
\le2\exp\Big(-\frac{u^2}{2A^2s}\Big).
\]
Apply this to $w_j$ and $w_0$, and take a union bound over the at most $(2^h+1)^k\le(2n+2)^k$ grid queries and $2n+1$ times. The triangle inequality gives
\[
\Prob(\Xi\ge2u\mid\text{locations})
\le4(2n+2)^{k+1}\exp\Big(-\frac{u^2}{2A^2s}\Big).
\]
Integrating yields
\begin{equation}\label{eq:Xi}
\E\,\Xi\le C_kA\sqrt{s\log(2n)}.
\end{equation}

Round each coordinate of a tagged threshold downward to the nearest cell boundary, using $0$ and $1$ as outer boundaries. The resulting grid query differs from the orthant only within at most one marginal cell in each coordinate. Each such cell contains at most $r$ kept points. The pilot is discarded, so \eqref{eq:grid-count} bounds every signed orthant count by $\Xi+kr$. The sample run is an initial part of the stream run, including when a source is exhausted, and its signed counts stay fixed after it stops. The tagged-order reduction therefore gives
\begin{equation}\label{eq:prefix-cell-bound}
\Dpre\le\Xi+kr.
\end{equation}
Since $h+1\le\log_2(p+1)+2\le\log_2n+2$, we obtain
\[
\E\Dpre\le C_kA\sqrt{s\log(2n)}+kr
\le C_{k,\rho}(1+\log n)^{k+1/2}.
\]
The bounds are uniform over all Borel input laws. This proves the upper bounds of Theorem~\ref{thm:multi}, with $D\le\Dpre$ for the terminal objective.

\end{proof}

\paragraph{A pathwise prefix bound on the line.}
For $k=1$, the coordinate bound $|w_j(B)-w_0(B)|<2A$ also gives $\Xi\le2As$ on every realization. Thus
\[
\Dpre\le2A(h+1)+r
=\frac{16\pi}{\rho}(h+1)^{3/2}+r
\le C_\rho(1+\log n)^{3/2}.
\]
For each fixed $\rho$, only finitely many $n$ use the keep-all branch, where $\Dpre\le2n$. Enlarging $C_\rho$ covers these values. Hence the same algorithm satisfies $\E R\le\rho n$ and this almost-sure bound for every $n\ge2$ and every Borel input law.

\section{One-dimensional thinning lower bounds}\label{sec:lower}

Throughout this section the input law is uniform on $[0,1]$ and $N=2n$. We use the adapted coefficients of Section~\ref{sec:general-model}. For terminal discrepancy, Wei's planar estimate forces a large bin count, and a Poisson clock and a long interarrival gap preserve one such count until the final observation.

\subsection{Random planar point sets}\label{sec:planar}

Following Wei \cite{Wei26}, for points $P_1,\dots,P_q$ of the unit square and coefficients $c\in\{-1,0,1\}^q$ put
\[
F_c(x,y)=\sum_{i=1}^qc_i\one\{P_i\in[0,x]\times[0,y]\},\qquad x,y\in[0,1].
\]
Here $\|F_c\|_\infty$ is the anchored-rectangle discrepancy, with zero coefficients representing deletions.

Wei \cite[Theorem 1.1]{Wei26} proves that for every $a>0$ there is $c_a>0$ such that $q$ independent uniform points in $[0,1]^2$ satisfy
\[
\Prob\big(\|F_\chi\|_\infty\ge c_a(\log_2q)^{3/2}\text{ for all }\chi\in\{-1,1\}^q\big)\ge1-e^{-aq}.
\]
A union bound over retained subsets gives the following estimate.

\begin{lemma}\label{lem:robust}
For every fixed $0<\beta<1$ there is $\kappa_\beta>0$ such that, for every $q\ge4$ and independent uniform points $P_1,\dots,P_q$ in $[0,1]^2$, with probability at least $1-e^{-3q}$ every $c\in\{-1,0,1\}^q$ with at least $\beta q$ nonzero entries satisfies
\[
\|F_c\|_\infty\ge\kappa_\beta(\log q)^{3/2}.
\]
\end{lemma}

\begin{proof}
First suppose $q\ge\beta^{-2}$. For each fixed subset $T\subseteq\{1,\dots,q\}$ with $|T|\ge\beta q$, apply Wei's theorem with exponential rate $4/\beta$ to the points indexed by $T$. Its failure probability is at most $e^{-4q}$, so a union bound over at most $2^q$ subsets gives total failure probability at most $2^qe^{-4q}\le e^{-3q}$. The inequality $\log(\beta q)\ge\tfrac12\log q$ gives the claimed threshold after changing the constant to use natural logarithms. For the finitely many smaller integers $q\ge4$, every nonempty support almost surely has distinct first coordinates. A rectangle $[0,x]\times[0,1]$ isolating its leftmost point has imbalance one. Decreasing $\kappa_\beta$ therefore handles these remaining sizes with probability one.
\end{proof}

\subsection{Poissonization}\label{sec:poisson}

Fix $0<\rho<2$ and put
\[
\beta=\frac{2-\rho}{8}\in(0,1/4).
\]
Let a controller have expected deletion count $\E R\le\rho n=(1-4\beta)N$. Its locations $U_1,U_2,\dots$ are independent and uniform, and its coefficients have the form $c_i=c_i(\omega,U_1,\dots,U_i)$ for $i\le N$. Extend them by $c_i=0$ for $i>N$.

Attach an independent Poisson clock of rate $(1-\beta)N$, with arrival times $0<T_1<T_2<\cdots$. By the marking theorem \cite[Chapter 5]{Kin93}, the pairs $(T_i,U_i)_{i\ge1}$ form a Poisson point process on $[0,\infty)\times(0,1)$ with intensity $(1-\beta)N\,dt\,du$. Let $K(t)=\#\{i:T_i\le t\}$ and $R(t)=\#\{i\le K(t):c_i=0\}$. Then
\begin{equation}\label{eq:R1poisson}
R(1)\le R+(K(1)-N)_+,\qquad
\E R(1)\le(1-4\beta)N+(1-\beta)N e^{-\beta^2N/2}.
\end{equation}
Indeed, for $K=K(1)\sim\Poi((1-\beta)N)$, the Poisson identity and Chernoff bound give
\begin{align*}
\E(K-N)_+
&\le\E[K\one\{K>N\}]
=(1-\beta)N\Prob(K\ge N)\\
&\le(1-\beta)N e^{-N[-\beta-\log(1-\beta)]}
\le(1-\beta)N e^{-\beta^2N/2}.
\end{align*}
Since $K(4)\sim\Poi(4(1-\beta)N)$ and $4(1-\beta)>3$, a Poisson lower-tail bound also gives
\begin{equation}\label{eq:TN}
\Prob(T_N>4)=\Prob(K(4)<N)\le e^{-(2-\log3)N}.
\end{equation}

Fix an integer $m\ge1$ and let $I_j=[(j-1)/m,j/m)$ for $1\le j\le m$. The arrivals in these bins form independent Poisson processes of rate
\[
\lambda=\frac{(1-\beta)N}{m},
\]
with independent marks uniform in their respective bins. Let $Q_j$ count the arrivals in $I_j$ before time $1$. Then $Q_j\sim\Poi(\lambda)$, and conditional on $Q_j=q$ the unordered set
\begin{equation}\label{eq:planar-points}
P_i=\big(mU_i-(j-1),T_i\big),\qquad U_i\in I_j,\quad T_i\le1,
\end{equation}
has the law of $q$ independent uniform points of $[0,1]^2$. Lemma~\ref{lem:robust} is invariant under relabeling, so it applies to these points with their arrival indices.

\begin{figure}[t]
\centering
\begin{tikzpicture}
\begin{scope}
\fill[black!10] (0,0) rectangle (3.472,3.192);
\draw (0,0) rectangle (5.6,5.6);
\draw[thick] (0,0) rectangle (3.472,3.192);
\fill (4.124,1.200) circle (1.6pt);
\fill (3.302,1.304) circle (1.6pt);
\draw[gray] (4.728,3.213) -- (4.848,3.333) (4.728,3.333) -- (4.848,3.213);
\fill (3.862,0.455) circle (1.6pt);
\draw[gray] (0.146,5.392) -- (0.266,5.512) (0.146,5.512) -- (0.266,5.392);
\draw[fill=white] (0.931,3.060) circle (1.6pt);
\fill (4.645,4.523) circle (1.6pt);
\fill (1.622,1.651) circle (1.6pt);
\draw[fill=white] (4.285,1.808) circle (1.6pt);
\draw[fill=white] (2.789,2.942) circle (1.6pt);
\fill (2.975,5.428) circle (1.6pt);
\fill (2.339,1.211) circle (1.6pt);
\draw[gray] (1.202,0.189) -- (1.322,0.309) (1.202,0.309) -- (1.322,0.189);
\draw[fill=white] (2.979,4.024) circle (1.6pt);
\fill (4.526,1.983) circle (1.6pt);
\draw[gray] (1.056,4.542) -- (1.176,4.662) (1.056,4.662) -- (1.176,4.542);
\fill (5.477,4.390) circle (1.6pt);
\draw[fill=white] (4.166,2.987) circle (1.6pt);
\draw[fill=white] (4.259,5.157) circle (1.6pt);
\fill (0.270,2.843) circle (1.6pt);
\draw[fill=white] (2.437,4.137) circle (1.6pt);
\draw[gray] (2.942,0.017) -- (3.062,0.137) (2.942,0.137) -- (3.062,0.017);
\fill (1.468,5.091) circle (1.6pt);
\fill (0.487,3.913) circle (1.6pt);
\fill (1.762,5.506) circle (1.6pt);
\fill (2.502,0.960) circle (1.6pt);
\draw[fill=white] (3.810,4.240) circle (1.6pt);
\draw[fill=white] (0.551,3.322) circle (1.6pt);
\draw[fill=white] (4.026,3.020) circle (1.6pt);
\draw[fill=white] (3.144,3.115) circle (1.6pt);
\draw[fill=white] (4.271,4.330) circle (1.6pt);
\draw[fill=white] (3.331,3.387) circle (1.6pt);
\draw[fill=white] (1.715,3.030) circle (1.6pt);
\draw[fill=white] (4.032,2.621) circle (1.6pt);
\fill (3.534,0.528) circle (1.6pt);
\draw[fill=white] (5.456,1.715) circle (1.6pt);
\draw[fill=white] (0.706,5.117) circle (1.6pt);
\draw[gray] (0.893,0.740) -- (1.013,0.860) (0.893,0.860) -- (1.013,0.740);
\draw[fill=white] (0.385,3.067) circle (1.6pt);
\draw[fill=white] (1.498,2.648) circle (1.6pt);
\fill (2.517,4.925) circle (1.6pt);
\fill (4.648,4.638) circle (1.6pt);
\fill (3.817,0.191) circle (1.6pt);
\draw[gray] (1.061,1.433) -- (1.181,1.553) (1.061,1.553) -- (1.181,1.433);
\draw[gray] (1.911,1.066) -- (2.031,1.186) (1.911,1.186) -- (2.031,1.066);
\draw[gray] (0.040,1.741) -- (0.160,1.861) (0.040,1.861) -- (0.160,1.741);
\fill (4.465,1.096) circle (1.6pt);
\draw[fill=white] (1.915,2.537) circle (1.6pt);
\node[font=\scriptsize, below] at (2.8,0) {location $U_i$};
\node[font=\scriptsize, rotate=90, above] at (0,2.8) {time coordinate $\vartheta_{(i)}$};
\node[font=\scriptsize, right] at (3.472,3.192) {$(x,\vartheta_{(t)})$};
\node[font=\small] at (2.8,6.0) {(a)};
\end{scope}
\begin{scope}[xshift=7.6cm]
\draw[black!30] (0.000,0) -- (0.000,5.6);
\draw[black!30] (1.233,0) -- (1.233,5.6);
\draw[black!30] (2.467,0) -- (2.467,5.6);
\draw[black!30] (3.700,0) -- (3.700,5.6);
\draw[black!30] (4.933,0) -- (4.933,5.6);
\draw[black!30] (6.167,0) -- (6.167,5.6);
\draw[black!30] (7.400,0) -- (7.400,5.6);
\draw (0,0) rectangle (7.4,5.6);
\fill[black!10] (2.467,0.806) rectangle (3.700,5.286);
\fill (0.300,4.912) circle (0.9pt);
\fill (0.158,2.688) circle (0.9pt);
\fill (0.465,5.373) circle (0.9pt);
\fill (0.524,0.540) circle (0.9pt);
\fill (0.780,1.901) circle (0.9pt);
\fill (0.510,3.103) circle (0.9pt);
\fill (0.152,3.085) circle (0.9pt);
\fill (0.706,3.860) circle (0.9pt);
\fill (0.785,2.786) circle (0.9pt);
\fill (1.022,0.859) circle (0.9pt);
\fill (1.091,0.175) circle (0.9pt);
\fill (0.350,1.784) circle (0.9pt);
\fill (0.935,2.172) circle (0.9pt);
\fill (0.939,0.078) circle (0.9pt);
\fill (0.266,4.733) circle (0.9pt);
\fill (0.109,5.570) circle (0.9pt);
\fill (0.965,0.115) circle (0.9pt);
\fill (0.940,2.091) circle (0.9pt);
\fill (0.362,3.705) circle (0.9pt);
\fill (0.140,0.942) circle (0.9pt);
\fill (0.270,2.610) circle (0.9pt);
\fill (0.131,1.913) circle (0.9pt);
\fill (0.506,0.022) circle (0.9pt);
\fill (0.322,2.670) circle (0.9pt);
\fill (0.393,4.283) circle (0.9pt);
\fill (1.066,0.268) circle (0.9pt);
\fill (0.600,1.743) circle (0.9pt);
\fill (0.346,4.981) circle (0.9pt);
\fill (1.039,5.417) circle (0.9pt);
\fill (0.551,4.216) circle (0.9pt);
\fill (1.032,1.471) circle (0.9pt);
\fill (0.956,5.497) circle (0.9pt);
\fill (0.336,5.185) circle (0.9pt);
\fill (0.137,0.328) circle (0.9pt);
\fill (1.113,0.964) circle (0.9pt);
\fill (0.972,2.302) circle (0.9pt);
\fill (0.140,1.375) circle (0.9pt);
\fill (0.600,1.506) circle (0.9pt);
\fill (0.880,0.549) circle (0.9pt);
\fill (0.873,0.205) circle (0.9pt);
\fill (1.106,1.955) circle (0.9pt);
\fill (0.217,5.253) circle (0.9pt);
\fill (1.089,2.858) circle (0.9pt);
\fill (0.335,5.105) circle (0.9pt);
\fill (0.359,5.404) circle (0.9pt);
\fill (0.296,3.885) circle (0.9pt);
\fill (0.176,2.694) circle (0.9pt);
\fill (0.346,4.964) circle (0.9pt);
\fill (0.857,3.010) circle (0.9pt);
\fill (0.473,1.585) circle (0.9pt);
\fill (1.987,2.836) circle (0.9pt);
\fill (2.178,3.855) circle (0.9pt);
\fill (1.450,0.798) circle (0.9pt);
\fill (2.319,3.675) circle (0.9pt);
\fill (2.259,3.397) circle (0.9pt);
\fill (1.923,3.406) circle (0.9pt);
\fill (2.317,0.706) circle (0.9pt);
\fill (1.605,5.242) circle (0.9pt);
\fill (2.218,5.003) circle (0.9pt);
\fill (1.583,3.533) circle (0.9pt);
\fill (1.349,0.219) circle (0.9pt);
\fill (1.543,0.051) circle (0.9pt);
\fill (1.801,5.243) circle (0.9pt);
\fill (2.312,0.407) circle (0.9pt);
\fill (2.347,4.500) circle (0.9pt);
\fill (1.865,5.364) circle (0.9pt);
\fill (1.769,5.512) circle (0.9pt);
\fill (2.225,5.541) circle (0.9pt);
\fill (1.978,5.215) circle (0.9pt);
\fill (1.545,0.071) circle (0.9pt);
\fill (1.790,1.098) circle (0.9pt);
\fill (2.092,4.680) circle (0.9pt);
\fill (2.036,1.829) circle (0.9pt);
\fill (2.178,3.259) circle (0.9pt);
\fill (1.713,1.232) circle (0.9pt);
\fill (1.418,5.067) circle (0.9pt);
\fill (1.453,1.665) circle (0.9pt);
\fill (1.537,2.467) circle (0.9pt);
\fill (1.553,5.461) circle (0.9pt);
\fill (1.582,3.982) circle (0.9pt);
\fill (1.712,2.048) circle (0.9pt);
\fill (1.991,2.807) circle (0.9pt);
\fill (1.601,1.127) circle (0.9pt);
\fill (1.732,2.050) circle (0.9pt);
\fill (1.700,5.369) circle (0.9pt);
\fill (1.575,1.733) circle (0.9pt);
\fill (2.158,2.587) circle (0.9pt);
\fill (1.391,3.187) circle (0.9pt);
\fill (1.531,0.651) circle (0.9pt);
\fill (1.392,0.069) circle (0.9pt);
\fill (1.626,0.165) circle (0.9pt);
\fill (1.659,1.777) circle (0.9pt);
\fill (1.476,1.094) circle (0.9pt);
\fill (1.564,0.328) circle (0.9pt);
\fill (1.768,1.343) circle (0.9pt);
\fill (1.873,1.187) circle (0.9pt);
\fill (2.360,2.269) circle (0.9pt);
\fill (1.990,0.747) circle (0.9pt);
\fill (1.517,5.320) circle (0.9pt);
\fill (1.816,0.501) circle (0.9pt);
\fill (3.014,0.575) circle (0.9pt);
\fill (3.522,0.150) circle (0.9pt);
\fill (3.270,5.300) circle (0.9pt);
\fill (2.757,0.030) circle (0.9pt);
\fill (3.297,5.371) circle (0.9pt);
\fill (2.765,0.415) circle (0.9pt);
\fill (2.656,0.426) circle (0.9pt);
\fill (2.641,0.345) circle (0.9pt);
\fill (3.078,0.599) circle (0.9pt);
\fill (3.841,2.120) circle (0.9pt);
\fill (4.224,1.569) circle (0.9pt);
\fill (4.682,0.954) circle (0.9pt);
\fill (4.400,0.730) circle (0.9pt);
\fill (4.275,0.154) circle (0.9pt);
\fill (3.818,2.597) circle (0.9pt);
\fill (4.354,1.883) circle (0.9pt);
\fill (4.391,4.596) circle (0.9pt);
\fill (3.796,0.613) circle (0.9pt);
\fill (4.689,0.213) circle (0.9pt);
\fill (4.542,4.311) circle (0.9pt);
\fill (4.703,0.489) circle (0.9pt);
\fill (3.917,4.965) circle (0.9pt);
\fill (4.319,2.045) circle (0.9pt);
\fill (4.669,4.882) circle (0.9pt);
\fill (4.369,5.560) circle (0.9pt);
\fill (4.852,4.484) circle (0.9pt);
\fill (4.397,2.990) circle (0.9pt);
\fill (4.230,4.868) circle (0.9pt);
\fill (4.763,1.775) circle (0.9pt);
\fill (4.142,3.863) circle (0.9pt);
\fill (4.377,3.801) circle (0.9pt);
\fill (4.627,5.058) circle (0.9pt);
\fill (4.688,4.952) circle (0.9pt);
\fill (4.284,2.391) circle (0.9pt);
\fill (4.701,1.806) circle (0.9pt);
\fill (4.312,5.512) circle (0.9pt);
\fill (4.211,1.135) circle (0.9pt);
\fill (4.208,2.335) circle (0.9pt);
\fill (4.013,4.193) circle (0.9pt);
\fill (3.935,4.361) circle (0.9pt);
\fill (3.961,1.791) circle (0.9pt);
\fill (4.677,5.232) circle (0.9pt);
\fill (4.827,0.845) circle (0.9pt);
\fill (4.417,4.152) circle (0.9pt);
\fill (4.310,2.313) circle (0.9pt);
\fill (3.991,1.645) circle (0.9pt);
\fill (4.520,4.208) circle (0.9pt);
\fill (3.842,3.394) circle (0.9pt);
\fill (3.860,3.692) circle (0.9pt);
\fill (4.802,2.986) circle (0.9pt);
\fill (4.442,5.257) circle (0.9pt);
\fill (3.902,3.148) circle (0.9pt);
\fill (4.212,3.305) circle (0.9pt);
\fill (4.848,2.209) circle (0.9pt);
\fill (4.575,5.170) circle (0.9pt);
\fill (4.529,1.083) circle (0.9pt);
\fill (3.959,3.622) circle (0.9pt);
\fill (4.359,3.924) circle (0.9pt);
\fill (4.169,1.461) circle (0.9pt);
\fill (6.078,2.494) circle (0.9pt);
\fill (5.745,2.769) circle (0.9pt);
\fill (5.644,4.525) circle (0.9pt);
\fill (5.471,0.910) circle (0.9pt);
\fill (5.343,5.270) circle (0.9pt);
\fill (6.013,4.352) circle (0.9pt);
\fill (5.347,4.314) circle (0.9pt);
\fill (5.614,0.495) circle (0.9pt);
\fill (5.321,1.265) circle (0.9pt);
\fill (5.443,0.933) circle (0.9pt);
\fill (5.186,1.528) circle (0.9pt);
\fill (5.935,4.630) circle (0.9pt);
\fill (5.363,4.024) circle (0.9pt);
\fill (5.178,2.782) circle (0.9pt);
\fill (5.611,0.448) circle (0.9pt);
\fill (5.814,4.206) circle (0.9pt);
\fill (5.732,0.355) circle (0.9pt);
\fill (5.022,5.569) circle (0.9pt);
\fill (5.596,1.004) circle (0.9pt);
\fill (5.633,0.479) circle (0.9pt);
\fill (5.610,5.538) circle (0.9pt);
\fill (5.183,0.286) circle (0.9pt);
\fill (5.700,4.929) circle (0.9pt);
\fill (6.016,2.811) circle (0.9pt);
\fill (5.685,0.603) circle (0.9pt);
\fill (5.723,3.566) circle (0.9pt);
\fill (5.784,5.499) circle (0.9pt);
\fill (5.218,1.931) circle (0.9pt);
\fill (5.167,0.357) circle (0.9pt);
\fill (6.024,1.457) circle (0.9pt);
\fill (5.151,4.952) circle (0.9pt);
\fill (5.058,3.599) circle (0.9pt);
\fill (5.676,4.678) circle (0.9pt);
\fill (5.449,0.883) circle (0.9pt);
\fill (5.642,1.275) circle (0.9pt);
\fill (5.337,4.669) circle (0.9pt);
\fill (6.021,2.962) circle (0.9pt);
\fill (5.500,4.119) circle (0.9pt);
\fill (5.223,4.261) circle (0.9pt);
\fill (5.284,1.247) circle (0.9pt);
\fill (5.184,3.907) circle (0.9pt);
\fill (5.896,3.999) circle (0.9pt);
\fill (5.583,4.408) circle (0.9pt);
\fill (5.587,4.126) circle (0.9pt);
\fill (5.182,5.012) circle (0.9pt);
\fill (5.877,3.218) circle (0.9pt);
\fill (5.767,5.476) circle (0.9pt);
\fill (6.036,1.719) circle (0.9pt);
\fill (5.132,4.603) circle (0.9pt);
\fill (5.527,0.100) circle (0.9pt);
\fill (6.504,2.769) circle (0.9pt);
\fill (7.235,2.962) circle (0.9pt);
\fill (6.539,4.386) circle (0.9pt);
\fill (6.264,4.086) circle (0.9pt);
\fill (7.282,0.903) circle (0.9pt);
\fill (6.546,5.269) circle (0.9pt);
\fill (6.294,1.321) circle (0.9pt);
\fill (6.553,0.572) circle (0.9pt);
\fill (6.855,1.638) circle (0.9pt);
\fill (7.016,0.509) circle (0.9pt);
\fill (6.560,2.871) circle (0.9pt);
\fill (6.631,3.777) circle (0.9pt);
\fill (6.569,1.944) circle (0.9pt);
\fill (7.211,1.644) circle (0.9pt);
\fill (6.526,4.874) circle (0.9pt);
\fill (6.952,0.916) circle (0.9pt);
\fill (6.900,3.541) circle (0.9pt);
\fill (7.237,1.191) circle (0.9pt);
\fill (6.247,5.148) circle (0.9pt);
\fill (7.302,1.194) circle (0.9pt);
\fill (6.300,2.125) circle (0.9pt);
\fill (6.422,2.171) circle (0.9pt);
\fill (7.215,2.892) circle (0.9pt);
\fill (6.252,1.048) circle (0.9pt);
\fill (6.393,2.676) circle (0.9pt);
\fill (7.160,2.035) circle (0.9pt);
\fill (6.523,5.269) circle (0.9pt);
\fill (6.527,1.709) circle (0.9pt);
\fill (6.789,2.140) circle (0.9pt);
\fill (6.305,4.283) circle (0.9pt);
\fill (7.038,1.066) circle (0.9pt);
\fill (6.994,5.344) circle (0.9pt);
\fill (6.656,1.813) circle (0.9pt);
\fill (7.216,0.031) circle (0.9pt);
\fill (7.051,5.514) circle (0.9pt);
\fill (7.289,0.138) circle (0.9pt);
\fill (6.763,0.263) circle (0.9pt);
\fill (6.560,4.208) circle (0.9pt);
\fill (7.145,2.019) circle (0.9pt);
\fill (6.632,1.257) circle (0.9pt);
\fill (7.096,1.883) circle (0.9pt);
\fill (6.775,5.277) circle (0.9pt);
\fill (6.650,0.853) circle (0.9pt);
\fill (7.011,3.805) circle (0.9pt);
\fill (7.086,4.920) circle (0.9pt);
\fill (6.295,1.601) circle (0.9pt);
\fill (6.645,0.678) circle (0.9pt);
\fill (7.232,3.144) circle (0.9pt);
\fill (6.303,0.021) circle (0.9pt);
\fill (6.794,3.924) circle (0.9pt);
\node[star, star points=5, fill=black, inner sep=1.4pt] at (3.083,0.806) {};
\draw[dashed] (0,1.120) -- (7.4,1.120); \node[font=\scriptsize, left] at (0,1.120) {$1$};
\draw[thick, densely dashdotted] (0,2.408) -- (7.4,2.408); \node[font=\scriptsize, left] at (0,2.408) {$T_N$};
\node[font=\scriptsize, left] at (0,5.600) {$5$}; \node[font=\scriptsize, left] at (0,0) {$0$};
\draw[thin] (2.467,5.286) -- (3.700,5.286); \node[font=\scriptsize, right] at (3.700,5.286) {\colorbox{white}{$\tau+4$}};
\node[font=\scriptsize, fill=white, inner sep=1pt, align=center] at (3.083,3.646) {no arrival\\ in $(\tau,\tau+4]$};
\node[font=\scriptsize, right] at (3.203,0.756) {$\tau$};
\node[font=\scriptsize, below] at (0.617,0) {$I_{1}$};
\node[font=\scriptsize, below] at (1.850,0) {$I_{2}$};
\node[font=\scriptsize, below] at (3.083,0) {$I_{3}$};
\node[font=\scriptsize, below] at (4.317,0) {$I_{4}$};
\node[font=\scriptsize, below] at (5.550,0) {$I_{5}$};
\node[font=\scriptsize, below] at (6.783,0) {$I_{6}$};
\node[font=\scriptsize, rotate=90, above] at (-0.55,2.8) {Poisson time};
\node[font=\small] at (3.7,6.0) {(b)};
\end{scope}
\end{tikzpicture}
\caption{The two lower bounds. (a) A prefix at threshold $x$ is the signed count in $[0,x]\times[0,\vartheta_{(t)}]$. Filled discs, open discs and crosses have coefficients $+1$, $-1$ and $0$. (b) Arrivals of the Poisson clock in the bins $I_1,\dots,I_m$. After a bin is selected at time $\tau\le1$, its next interarrival gap exceeds $4$ with probability $e^{-4\lambda}$, and on the event $T_N\le4$ such a gap preserves its local discrepancy until the terminal time.}
\label{fig:lower}
\end{figure}
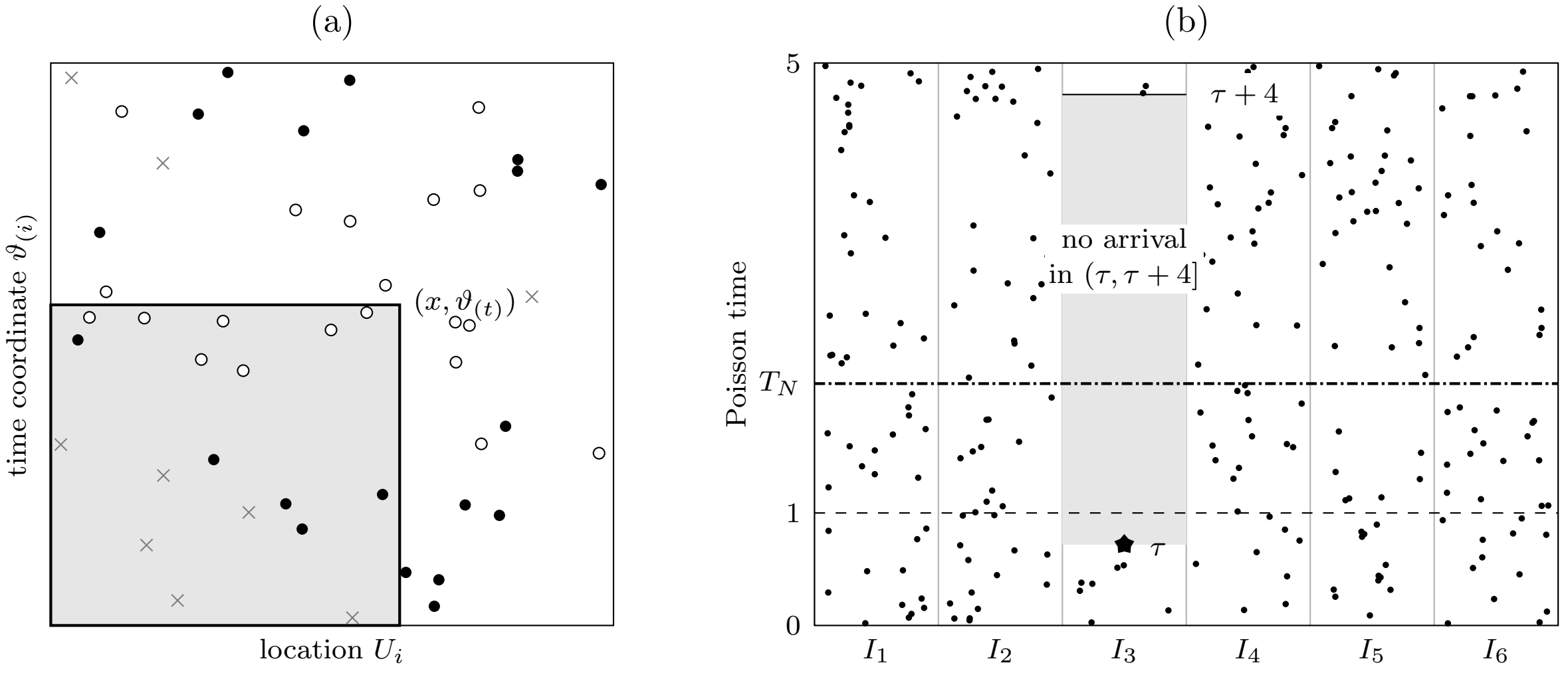

The local discrepancy of bin $j$ at time $t$ is
\begin{equation}\label{eq:local-disc}
\Delta_j(t)=\sup_{b\in I_j}\Big|\sum_{i:\ T_i\le t,\ U_i\in I_j,\ U_i\le b}c_i\Big|.
\end{equation}
It is determined by the joint past and changes only at an arrival in $I_j$ with a nonzero coefficient. For the arrivals of bin $j$ before time $1$, a planar rectangle selects exactly the points to the left of a spatial threshold and before a time threshold. Hence
\begin{equation}\label{eq:planar-local}
\|F_c\|_\infty=\sup_{0\le t\le1}\Delta_j(t).
\end{equation}
The controller's decisions remain measurable with respect to the joint past after the independent clock is added. Figure~\ref{fig:lower}(b) illustrates this embedding, while panel (a) depicts the construction used for the prefix lower bound.

\subsection{Long interarrival gaps}\label{sec:freezing}

Consider $m$ independent Poisson processes of rate $\lambda$, carrying independent marks whose law may depend on the process. For process $i$, write $\zeta_{i,\ell}$ for its independent exponential gaps of mean $1/\lambda$ and $\theta_{i,\ell}=\zeta_{i,1}+\cdots+\zeta_{i,\ell}$ for its arrival times. All gaps and marks are independent, and a seed $\omega$ is independent of them. The past at time $t$ consists of $\omega$ and of the process indices, times and marks of all arrivals up to $t$.

A rule selects at most one arrival from each process, using only the joint past at that arrival. Let $H>0$. A selected arrival $\theta_{i,\ell}$ succeeds if $\zeta_{i,\ell+1}>H$, so that its process has no arrival in the next $H$ time units, and fails otherwise. Let $\tau_j$ be the time of the $j$-th selected arrival, with $\tau_j=\infty$ if fewer than $j$ are selected.

\begin{lemma}\label{lem:freezing}
For every integer $j\ge1$,
\[
\Prob\big(\tau_j\le1\text{ and the first }j\text{ selections all fail}\big)\le\big(1-e^{-\lambda H}\big)^j.
\]
\end{lemma}

\begin{proof}
Put $\delta=e^{-\lambda H}$. For $j\ge1$, $1\le i\le m$ and $\ell\ge1$, let $A_{j,i,\ell}$ be the event that the $\ell$-th arrival of process $i$ is the $j$-th selected arrival, that $\theta_{i,\ell}\le1$, and that the first $j-1$ selections fail. Let $\cF_{i,\ell}$ be the $\sigma$-field generated by $\omega$, by the first $\ell$ arrival times and their marks in process $i$, and by the complete marked paths of all other processes. The time $\theta_{i,\ell}$ and the first $\ell$ arrivals of process $i$ are $\cF_{i,\ell}$-measurable. Together with the complete paths of the other processes, they determine the joint past at $\theta_{i,\ell}$. Hence all selection decisions up to time $\theta_{i,\ell}$ are $\cF_{i,\ell}$-measurable, in particular whether the $\ell$-th arrival of process $i$ is selected and whether it is the $j$-th selected arrival. If it is, the first $j-1$ selections lie on processes $i'\ne i$, because a process is selected at most once, and whether each of them fails is determined by the complete path of its process. Therefore $A_{j,i,\ell}\in\cF_{i,\ell}$. The gap $\zeta_{i,\ell+1}$ is independent of $\cF_{i,\ell}$, so
\[
\Prob\big(A_{j,i,\ell}\cap\{\zeta_{i,\ell+1}\le H\}\big)=(1-\delta)\,\Prob(A_{j,i,\ell}).
\]
The $j$-th selected arrival, if it occurs by time $1$, is exactly one arrival $(i,\ell)$, and it fails if and only if $\zeta_{i,\ell+1}\le H$. Hence
\begin{align*}
\Prob(\tau_j\le1,\ \text{first }j\text{ fail})&=\sum_{i,\ell}\Prob\big(A_{j,i,\ell}\cap\{\zeta_{i,\ell+1}\le H\}\big)\\
&=(1-\delta)\sum_{i,\ell}\Prob(A_{j,i,\ell})=(1-\delta)\,\Prob(\tau_j\le1,\ \text{first }j-1\text{ fail}).
\end{align*}
Since $\tau_{j-1}\le\tau_j$, the last probability is at most $\Prob(\tau_{j-1}\le1,\ \text{first }j-1\text{ fail})$ for $j\ge2$, and at most $1$ for $j=1$. Induction on $j$ proves the lemma.
\end{proof}

\subsection{The terminal lower bounds}\label{sec:terminal-allrho}

\begin{proof}[Proof of Theorem~\ref{thm:terminal}(ii) and Corollary~\ref{cor:terminal-allrho}]
The upper bound in the corollary follows by applying Theorem~\ref{thm:terminal}(i) with deletion parameter $\min\{\rho,1/16\}$.

For the lower bounds, use the Poisson embedding of Section~\ref{sec:poisson} and put
\begin{equation}\label{eq:lower-params}
m=\Big\lceil\frac{8(1-\beta)N}{\log N}\Big\rceil,\qquad
m_0=\Big\lfloor\frac{\beta m}{4}\Big\rfloor.
\end{equation}
For $N\ge8$, we have $(1-\beta)N\ge\log N$, so
\begin{equation}\label{eq:lambda-bounds}
\frac{\log N}{9}\le\lambda\le\frac{\log N}{8}.
\end{equation}
Take $N$ large enough that $\lambda\ge8$ and $m_0\ge1$. Call bin $j$ good if $Q_j\ge(1-\beta)\lambda$ and its points \eqref{eq:planar-points} satisfy Lemma~\ref{lem:robust} with retained fraction $\beta$. A Poisson lower-tail bound and the lemma give
\[
\Prob(\text{bin }j\text{ is not good})
\le e^{-\beta^2\lambda/2}+e^{-3(1-\beta)\lambda}
=:\varepsilon_N,
\]
where the planar estimate applies because $(1-\beta)\lambda\ge6$.

Set
\[
B=\frac{\kappa_\beta}{2}\bigl(\log((1-\beta)\lambda)\bigr)^{3/2}.
\]
In each bin select the first arrival at which $\Delta_j$ exceeds $B$, if such an arrival occurs. This rule uses only the joint past and selects at most one arrival per bin, so Lemma~\ref{lem:freezing} applies. Let $J$ be the number of bins selected by time $1$.

If a good bin has no selection by time $1$, fewer than $\beta Q_j$ of its arrivals before that time can have been retained. Otherwise Lemma~\ref{lem:robust} and \eqref{eq:planar-local} would give
\[
\sup_{0\le t\le1}\Delta_j(t)
\ge\kappa_\beta(\log Q_j)^{3/2}\ge2B.
\]
Such a bin therefore contributes more than $(1-\beta)Q_j\ge(1-\beta)^2\lambda$ deletions to $R(1)$. Consequently,
\[
J\ge\sum_{j=1}^m\one\{\text{bin }j\text{ is good}\}
-\frac{R(1)}{(1-\beta)^2\lambda}.
\]
Taking expectations and using $m\lambda=(1-\beta)N$ and \eqref{eq:R1poisson}, we obtain
\[
\frac{\E J}{m}
\ge1-\frac{1-4\beta}{(1-\beta)^3}
-\varepsilon_N-(1-\beta)^{-2}e^{-\beta^2N/2}.
\]
Since $0\le J\le m$ and $m_0/m\le\beta/4$, the bound
$\E J\le m_0+m\Prob(J\ge m_0)$ implies
\[
\Prob(\tau_{m_0}\le1)
=\Prob(J\ge m_0)
\ge1-\frac{1-4\beta}{(1-\beta)^3}-\frac\beta4
-\varepsilon_N-(1-\beta)^{-2}e^{-\beta^2N/2}.
\]

Apply Lemma~\ref{lem:freezing} with $H=4$. By \eqref{eq:lambda-bounds} and \eqref{eq:lower-params},
\[
e^{-4\lambda}\ge N^{-1/2},\qquad
m_0\ge\frac{2\beta(1-\beta)N}{\log N}-1.
\]
Thus
\begin{equation}\label{eq:all-fail}
\Prob(\tau_{m_0}\le1,\ \text{first }m_0\text{ selections fail})
\le e^{-m_0e^{-4\lambda}}=o_\rho(1).
\end{equation}
If a selection succeeds at time $\tau\le1$ and $T_N\le4$, its bin has no arrival in $(\tau,\tau+4]$. A selection requires a nonzero coefficient, hence $\tau\le T_N$. Each local half-line count is the difference of two global half-line counts. Hence
\[
2D\ge\Delta_j(T_N)=\Delta_j(\tau)>B.
\]
Combining the crossing probability with \eqref{eq:all-fail} and \eqref{eq:TN} gives
\begin{equation}\label{eq:terminal-probability}
\begin{split}
\Prob(D>B/2)
\ge{}&1-\frac{1-4\beta}{(1-\beta)^3}-\frac\beta4
-\varepsilon_N-(1-\beta)^{-2}e^{-\beta^2N/2}\\
&{}-e^{-m_0e^{-4\lambda}}-e^{-(2-\log3)N}.
\end{split}
\end{equation}
All error terms tend to zero for each fixed $\rho<2$, uniformly over controllers with $\E R\le\rho n$.

Bernoulli's inequality $(1-\beta)^4\ge1-4\beta$ gives
\[
1-\frac{1-4\beta}{(1-\beta)^3}-\frac\beta4\ge\frac{3\beta}{4}.
\]
Taking $N$ sufficiently large that the error sum in \eqref{eq:terminal-probability} is at most $\beta/4$ therefore gives
\[
\Prob(D>B/2)\ge\frac\beta2=\frac{2-\rho}{16}.
\]
Moreover, $(1-\beta)\lambda\ge(\log N)/12$, so
$B/2\ge c_\rho(\log\log n)^{3/2}$ for all sufficiently large $n$. This proves the probability bound in the corollary. Its expectation consequence, together with the upper bound, gives the claimed minimax order.

For Theorem~\ref{thm:terminal}(ii), take $\rho=1/8$, so that $\beta=15/64$ and the controller budget is $\E R\le N/16=n/8$. For all sufficiently large $N$, \eqref{eq:terminal-probability} gives \eqref{eq:terminal-lower}, since
\[
1-\frac{1-4\beta}{(1-\beta)^3}-\frac\beta4
\ge1-\frac{4}{27}-\frac1{16}>\frac14.\qedhere
\]
\end{proof}

\section*{AI disclosure}

ChatGPT was used to generate ideas, assist with the proof, review draft, check main argument, revise the manuscript and preparation of the mathematical exposition. Claude was used to assist with preparing the illustrations and algorithms for the manuscript.

\end{document}